\documentclass[BCOR8mm,DIV11,cleardoubleempty,openright,a4paper,
bibtotoc,idxtotoc,twoside,final]{scrbook}
\pdfoutput=1
\usepackage[automark]{scrpage2}
\usepackage[ngerman, english]{babel}
\usepackage[backref=page,final]{hyperref}
\usepackage[T1]{fontenc}
\usepackage{lmodern}
\usepackage[leqno]{amsmath}
\usepackage{amssymb,amscd}
\usepackage{amsthm}
\usepackage{amsopn}
\usepackage[all]{xypic}
\usepackage{mathtools}
\usepackage{bbm}
\usepackage{listings}
\usepackage{verbatim}
\usepackage{braket}
\usepackage{algorithmic}
\ifnum\pdfoutput=0
\usepackage[final]{graphicx}
\else
\usepackage[pdftex,final]{graphicx}
\usepackage{epstopdf}
\fi
\renewcommand{\emph}[1]{{\bf #1}}

\newcommand{\F}{\mathbb{F}}

\newcommand{\C}{\mathbb{C}}
\newcommand{\N}{\mathbb{N}}
\newcommand{\Z}{\mathbb{Z}}
\newcommand{\Q}{\mathbb{Q}}

\newcommand{\Proj}{\mathbb{P}}
\newcommand{\Segre}{\Sigma}
\newcommand{\SegreK}{\mathcal{S}}

\newcommand{\Orbit}{\mathcal{O}}

\newtheorem{defn}{Definition}[chapter]
\newtheorem{lmm}[defn]{Lemma}
\newtheorem{cory}[defn]{Corollary}
\newtheorem{theom}[defn]{Theorem}
\newtheorem{mtheom}{Theorem}
\newtheorem{propn}[defn]{Proposition}
\theoremstyle{remark}

\newtheorem{beispiel}[defn]{Example}
\newtheorem{remark}[defn]{Remark}

\DeclareMathOperator{\rank}{rank}
\DeclareMathOperator{\supp}{supp}

\DeclareMathOperator{\GL}{GL}
\DeclareMathOperator{\Tr}{Tr}

\DeclareMathOperator{\End}{End}
\DeclareMathOperator{\Aut}{Aut}
\DeclareMathOperator{\Hom}{Hom}
\DeclareMathOperator{\Ext}{Ext}
\DeclareMathOperator{\ind}{ind}
\DeclareMathOperator{\add}{add}
\DeclareMathOperator{\pd}{pd}
\DeclareMathOperator{\gldim}{gldim}

\DeclareMathOperator{\ext}{ext}
\DeclareMathOperator{\Ker}{Ker}
\DeclareMathOperator{\Bild}{Im}
\DeclareMathOperator{\Mod}{mod}
\DeclareMathOperator{\MMod}{Mod}

\DeclareMathOperator{\rad}{rad}
\DeclareMathOperator{\id}{id}

\DeclareMathOperator{\Mat}{Mat}
\DeclareMathOperator{\soc}{soc}
\DeclareMathOperator{\codim}{codim}
\DeclareMathOperator{\chark}{char}
\DeclareMathOperator{\Spec}{Spec}
\DeclareMathOperator{\ProjS}{Proj}

\DeclareMathOperator{\OGr}{Gr}
\DeclareMathOperator{\OFl}{Fl}
\DeclareMathOperator{\ORep}{Rep}
\DeclareMathOperator{\Orep}{rep}

\DeclareMathOperator{\ORepFl}{RepFl}
\DeclareMathOperator{\ORepHom}{RepHom}
\newcommand{\Sch}[1]{\mathrm{#1}}
\newcommand{\RMon}{\mathcal{SR}}
\newcommand{\EMon}{\mathcal{ER}}
\newcommand{\TEMon}{\widetilde{\mathcal{ER}}}
\newcommand{\SRoot}[1]{\{#1\}}
\newcommand{\flvec}[1]{\seqv{\dvec{#1}}}
\newcommand{\tensor}{\otimes}
\newcommand{\qbinom}[2]{\begin{bmatrix}#1\\#2\end{bmatrix}}
\newcommand{\qnum}[1]{\left[ #1 \right]}
\newcommand{\rep}[1]{\ORep(#1)}
\newcommand{\repK}[2]{\ensuremath{\Orep(#1, #2)}}

\newcommand{\Rep}[2][]{\ensuremath{\ORep_{#1}({#2})}}

\newcommand{\RepK}[3][]{\ensuremath{\ORep_{#1}({#2}, {#3})}}
\newcommand{\Gr}[3][]{\ensuremath{\OGr_{#1}\binom{#3}{#2}}}

\newcommand{\RepFl}[2][]{\ensuremath{\ORepFl_{#1}\left({#2}\right)}}
\newcommand{\RepHom}[2][]{\ensuremath{\ORepHom_{#1}\left({#2}\right)}}
\newcommand{\Fl}[3][]{\ensuremath{\OFl_{#1}\binom{#3}{#2}}}
\newcommand{\HFl}[3][]{\ensuremath{\hat \OFl_{#1}\binom{#3}{#2}}}

\newcommand{\euf}[1]{\left\langle {#1} \right\rangle}
\newcommand{\lrang}[1]{\left\langle {#1} \right\rangle}
\newcommand{\bform}[2]{\lrang{\,{#1}\, {,}\, {#2}\,}}
\newcommand{\sbform}[2]{\left({#1}\, {,}\, {#2}\right)}
\newcommand{\dvec}[1]{{\underline{#1}}}
\newcommand{\ledeg}{\le_{\mathrm{deg}}}

\newcommand{\degr}[1]{\lvert {#1} \rvert}
\newcommand{\dimv}[1]{\underline{\dim}\left(#1\right)}
\newcommand{\dimve}[1][]{\underline{\dim}\;{#1}}
\newcommand{\seqv}[1]{{\boldsymbol{#1}}}
\newcommand{\rever}[1]{\ensuremath\stackrel{\leftarrow}{#1}}

\def\clap#1{\hbox to 0pt{\hss#1\hss}}
\def\mathclap{\mathpalette\mathclapinternal}
\def\mathclapinternal#1#2{%
	\clap{$\mathsurround=0pt#1{#2}$}%
}

\usepackage[figuresright]{rotating}
\usepackage{makeidx}
\usepackage{lscape}
\usepackage{setspace}
\usepackage{tikz}
\usetikzlibrary{matrix,arrows,decorations.pathmorphing,calc}
\author{\textbf{Stefan Wolf}}
\date{2009}
\titlehead{
\resizebox{5cm}{!}{\includegraphics{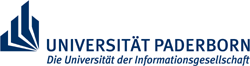}}\hfill
\begin{tabular}{r}
Fakult\"at f\"ur\\
Elektrotechnik, Informatik und Mathematik\\
Institut f\"ur Mathematik
\end{tabular}\hfill
}
\subject{Dissertation}
\title{\fontsize{35}{35}\selectfont {\Huge The}\\[1ex] Hall Algebra\\[1ex] {\Huge and the}\\[1ex] Composition Monoid}
\publishers{%
\begin{tabular}{rl}
\large Betreuer:& \large Andrew Hubery, Ph.D.\\
&\large Prof. Dr. Henning Krause
\end{tabular}
}
\makeindex
\newbox\ASYbox
\newdimen\ASYdimen
\long\def\ASYbase#1#2{\leavevmode\setbox\ASYbox=\hbox{#1}\ASYdimen=\ht\ASYbox%
\setbox\ASYbox=\hbox{#2}\lower\ASYdimen\box\ASYbox}
\long\def\ASYaligned(#1,#2)(#3,#4)#5#6#7{\leavevmode%
\setbox\ASYbox=\hbox{#7}%
\setbox\ASYbox\hbox{\ASYdimen=\ht\ASYbox%
\advance\ASYdimen by\dp\ASYbox\kern#3\wd\ASYbox\raise#4\ASYdimen\box\ASYbox}%
\put(#1,#2){#5\wd\ASYbox 0pt\dp\ASYbox 0pt\ht\ASYbox 0pt\box\ASYbox#6}}%
\long\def\ASYalignT(#1,#2)(#3,#4)#5#6{%
\ASYaligned(#1,#2)(#3,#4){%
}{%
}{#6}}
\long\def\ASYalign(#1,#2)(#3,#4)#5{\ASYaligned(#1,#2)(#3,#4){}{}{#5}}

\begin{document}
\frontmatter
\selectlanguage{ngerman}
\bibliographystyle{amsalphaurl}
\maketitle
%
\selectlanguage{english}
\pagestyle{scrheadings}
\cleardoublepage
\chapter*{}
\section*{Abstract}
Let $Q$ be a quiver.
M. Reineke and A. Hubery investigated the connection between the composition monoid
$\mathcal{CM}(Q)$, as introduced by M. Reineke,
and the generic composition algebra $\mathcal{C}_q(Q)$, as introduced by C. M. Ringel, specialised at $q=0$.
In this thesis we continue their work. We show that if $Q$ is a Dynkin quiver or an oriented cycle, then
$\mathcal{C}_0(Q)$ is isomorphic to the monoid algebra $\Q\mathcal{CM}(Q)$.
Moreover, if $Q$ is an acyclic, extended Dynkin quiver, we show that there exists
a surjective homomorphism
$\Phi\colon \mathcal{C}_0(Q) \rightarrow \Q\mathcal{CM}(Q)$, and we describe
its non-trivial kernel.

Our main tool is a geometric version of BGP reflection functors
on quiver Grassmannians and quiver flags, that is varieties consisting of
filtrations of a fixed representation by subrepresentations
of fixed dimension vectors.
These functors enable us to calculate various structure constants of the composition algebra.

Moreover, we investigate geometric properties of
quiver flags  and quiver Grassmannians, and show that under certain conditions,
quiver flags are irreducible and smooth.  If, in addition, we have a
counting polynomial, these properties imply the
positivity of the Euler characteristic of the quiver flag.

\section*{Zusammenfassung}
\selectlanguage{ngerman}
Sei $Q$ ein K\"ocher.
M. Reineke und A. Hubery untersuchten den Zusammenhang zwischen dem von M. Reineke eingef"uhrten Kompositionsmonoid
$\mathcal{CM}(Q)$ und der bei $q=0$ spezialisierten Kompositionsalgebra $\mathcal{C}_q(Q)$, die von C. M. Ringel definiert wurde.
Diese Dissertation f"uhrt diese Arbeit fort. Wir zeigen, dass $\mathcal{C}_0(Q)$ isomorph zu der Monoidalgebra $\Q\mathcal{CM}(Q)$
ist, wenn $Q$ ein Dynkin
K"ocher oder ein orientierter Zykel ist. Wenn $Q$ ein azyklischer erweiterter Dynkin K"ocher ist, so zeigen wir, dass es einen surjektiven Homomorphismus
$\Phi \colon \mathcal{C}_0(Q) \rightarrow \Q\mathcal{CM}(Q)$ gibt, und beschreiben dessen
nicht trivialen Kern.

Um dies zu beweisen, m"ussen wir viele Strukturkonstanten der Kompositionsalgebra berechnen.
Dazu f"uhren wir eine geometrische Version von BGP-Spiegelungsfunktoren
auf K"ochergrassmannschen und K"ocherfahnen ein. Dies sind Variet\"aten bestehend aus
Filtrierungen einer festen Darstellung durch Unterdarstellungen fester Dimensionvektoren.

Au"serdem untersuchen wir die geometrischen Eigenschaften von K"ochergrassmannschen und K"ocherfahnen. Wir
erhalten ein Kriterium, das uns erlaubt festzustellen, wann eine K"ocherfahne
irreduzibel und glatt ist. Wenn man zus"atzlich noch ein Z"ahlpolynom hat, so folgt aus diesen Eigenschaften die Positivit"at der
Euler-Charakteristik.
\selectlanguage{english}

\tableofcontents
\mainmatter
\chapter{Introduction}

The representation theory of quivers has its origin in a paper by P. Gabriel
\cite{Gabriel_unzdst}, who showed that a connected quiver is of finite
representation type if and only if the underlying graph is a Dynkin diagram of type
$A$, $D$ or $E$. In doing so, he observed that there is a strong connection
to the theory of Lie algebras. Namely, there is a bijection
between isomorphism classes of indecomposable representations of the quiver
and the set of positive roots of the associated complex Lie algebra. In
\cite{BernsteinGelfandPonomarev_coxeter} I. Bern{\v s}te{\u\i}n, I. Gel{$'$}fand and V. Ponomarev
gave a more direct proof of this result. More precisely, they used sequences of reflection
functors to obtain all indecomposable representations in a similar way
to how all roots are obtained by applying reflections in the Weyl group
to the simple roots.

P. Donovan and M. Freislich \cite{DonovanFreislich_repsoffingraphs} and independently
L. Nazarova \cite{Nazarova_repsquivinf} extended this work
to cover quivers of extended Dynkin, or affine, type. They described
the set of isomorphism classes of indecomposables, therefore showing that
these quivers are of tame representation type. A unified approach,
which can also be used for species of Dynkin or affine type, can be found in
\cite{DlabRingel_repsgraphsalgs}. Again there is a connection with the root
systems of affine Kac-Moody Lie algebras, namely the dimension vectors of the indecomposables
are exactly the positive roots. This does not extend to a bijection with
the isomorphism classes, since for each imaginary root there is a continuous family
of indecomposables.

V. Kac \cite{Kac_infiniterootsystems} proved that this correspondence holds in general. That is,
for a finite quiver without vertex loops the set of dimension vectors
of indecomposables over an algebraically closed field is precisely the
set of positive roots of the corresponding (symmetric) Kac-Moody Lie
algebra. Finally, A. Hubery \cite{Hubery_quiverauto} established this correspondence
in the case of species and general Kac-Moody Lie algebras.

These results hint towards a deep connection between the category of representations
of a quiver and the corresponding Kac-Moody Lie algebra.
This was further strengthened by the theory of Ringel-Hall algebras.
The Hall algebra appeared at first in the work of E. Steinitz \cite{Steinitz_abelsch} and afterwards
in the work of P. Hall \cite{Hall_part}. C. M. Ringel \cite{Ringel_hallalgsandquantumgroups} generalised this
construction to obtain an associative algebra structure
on the $\Q$-vector space $\mathcal{H}_K(Q)$ with
basis the isomorphism classes of representations of $Q$ over a finite field $K$, the Ringel-Hall algebra.
The structure constants are given basically by counting numbers of extensions.

The whole Ringel-Hall algebra is generally too complicated, and therefore one introduces the $\Q$-subalgebra
$\mathcal{C}_K(Q)$ generated by the isomorphism classes of simple representations without self-extensions.
There is a generic version $\mathcal{C}_q(Q)$, a $\Q(q)$-algebra, such that
specialising $q$ to $|K|$ recovers $\mathcal{C}_K(Q)$.
C.M. Ringel \cite{Ringel_hallalgsandquantumgroups}, for the Dynkin case, and later
J. Green \cite{Green_hallalgs}, for the general case, showed that, after twisting the multiplication with
the Euler form of $Q$, the generic composition algebra $\mathcal{C}_{q}(Q)$
is isomorphic
to the
positive part of the quantised enveloping algebra of the Kac-Moody
Lie algebra corresponding to $Q$.

When doing calculations in the Hall algebra one often has to decide
which representations are an extension of two other fixed representations.
If we fix two representations, then there is the easiest
extension, the direct sum. For the Dynkin case there is also the most complicated extension,
the generic extension. If we work over an algebraically closed field,
then the closure of the orbit of the generic extension contains all other extensions.
M. Reineke \cite{Reineke_genericexts} was the first to notice that the
multiplication by taking generic extensions is associative.
We therefore obtain a monoid structure on the set of isomorphism classes of
representations.

If we work over a non-Dynkin quiver, then, in general, generic extensions do not exist anymore.
M. Reineke \cite{Reineke_monoid} showed that one can solve this problem by, instead of taking
individual representations, taking irreducible closed subvarieties of the representation
variety to obtain a similar result. The multiplication is then
given by taking all possible extensions. Again, this yields a monoid,
the generic extension monoid $\mathcal{M}(Q)$. Similarly to the Hall algebra,
it is in general too complicated, so one restricts itself to the
submonoid generated by the orbits of simple representations without self-extensions,
the composition monoid $\mathcal{CM}(Q)$.
The elements of the composition monoid are the varieties consisting of representations
having a composition series with prescribed composition factors in prescribed order.

For the generic composition algebra viewed as a $\Q(v)$-algebra with $v^2 =q$
the (twisted) quantum Serre relations are defining as shown in \cite{Ringel_greenstheom}. M. Reineke
showed that the quantum Serre relations specialised to
$q=0$ hold in the composition monoid. They are in general not defining any more if we
specialise $q$ to $0$ in the composition algebra. But nonetheless one can conjecture,
as M. Reineke did in \cite{Reineke_genericexts} and \cite{Reineke_monoid}, that there is
a homomorphism of $\Q$-algebras
 \[\Phi \colon \mathcal{C}_0(Q) \rightarrow \Q\mathcal{CM}(Q)\]
sending simples to simples and therefore being automatically surjective.

The first step in this direction was done by A. Hubery \cite{Hubery_kronecker} showing that for
the Kronecker quiver
\[K =
\begin{array}{@{}c@{}}
\makeatletter%
\let\ASYencoding\f@encoding%
\let\ASYfamily\f@family%
\let\ASYseries\f@series%
\let\ASYshape\f@shape%
\makeatother%
\setlength{\unitlength}{1pt}
\includegraphics{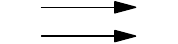}%
\definecolor{ASYcolor}{gray}{0.000000}\color{ASYcolor}
\fontsize{12.000000}{14.400000}\selectfont
\usefont{\ASYencoding}{\ASYfamily}{\ASYseries}{\ASYshape}%
\ASYalign(-46.424798,6.262465)(-0.500000,-0.500000){$\bullet$}
\ASYalign(-4.790163,6.262465)(-0.500000,-0.500000){$\bullet$}
\end{array}
\]
$\Phi$ is a homomorphism of $\Q$-algebras with non-trivial kernel. He did this by
calculating defining relations for $\mathcal{C}_0(K)$ and $\mathcal{CM}(K)$. He also was
able to give generators for the kernel of $\Phi$.

The main aim of this thesis is to extend this result to the Dynkin and extended Dynkin case.
More precisely we show that if $Q$ is a Dynkin quiver or an oriented cycle, then
$\mathcal{C}_0(Q)$ and $\Q\mathcal{CM}(Q)$ are isomorphic. Then we
prove that, if $Q$ is an acyclic, extended Dynkin quiver, there is such a morphism $\Phi$ and
the kernel of this morphism is given by the same relations as A. Hubery gave for
the Kronecker quiver.

In order to do this we need to calculate many of the structure
coefficients in the composition algebra. These are given by counting
points of quiver flags of a representation
$M$, i.e. increasing sequences of subrepresentations of fixed dimension vectors,
over finite fields.

For the oriented cycle we show directly that the coefficients at $q=0$ are
only one or zero. As an immediate consequence we obtain the desired isomorphism.

To obtain results for the Dynkin and extended Dynkin case we first
develop a framework for applying reflection functors
to quiver flags. Having done so, we can immediately show that for the Dynkin
case the coefficients at $q=0$ are also one or zero and by this we obtain
the desired isomorphism.

In the extended Dynkin case there will be coefficients
which are not equal to one or zero. Therefore, the proof becomes more
involved, but by using the framework of reflection functors on flags
it is possible to obtain the morphism. Along the way we obtain a
basis of PBW-type for $\mathcal{C}_0(Q)$, a normal form
for elements of $\mathcal{CM}(Q)$ and show that the second one does not depend
on the choice of the algebraically closed field, making geometric
arguments possible.

The varieties which appear, namely quiver Grassmannians and quiver flags, are of interest of their own.
They are
interesting projective varieties to study and their Euler
characteristics give coefficients in the cluster algebra as shown in
\cite{CalderoChapoton_clusterhall}. In general, they are neither smooth nor irreducible or even
reduced. We show that under some extra conditions they are smooth and irreducible. We
then use this information to calculate the constant coefficient of the counting polynomial
if it exists. Moreover, if we have a counting polynomial, we show that
the Euler characteristic is positive for rigid modules, and this implies
a certain positivity result for the associated cluster algebra. This works in the Dynkin and
the extended Dynkin case without using G. Lusztig's interpretation in terms of perverse sheaves,
and therefore gives an independent proof of positivity.

\subsection*{Outline}
This thesis is organised as follows:
In chapter 2 we recall basic facts about the representation theory of quivers and related topics. In chapter 3
we show, by a direct calculation of some Hall polynomials, that the Hall algebra of a cyclic quiver at
$q=0$ is isomorphic to the generic extension monoid. Then, in chapter 4, we develop normal
forms in terms of Schur roots for the composition monoid of an extended Dynkin quiver,
showing that it is independent of the choice of the algebraically closed field. These were the direct
results. In chapter 5 we focus our attention on the geometry of quiver flags. We prove
a dimension estimate, generalising the one of M. Reineke, and then show that in certain cases
the quiver flag varieties are smooth and irreducible. If one additionally has a counting polynomial,
we use this to deduce that the constant coefficient is one modulo $q$ and the Euler characteristic is
positive. In chapter 6 we develop the aforementioned calculus of reflections on quiver flags, culminating
in the proof that the composition algebra of a Dynkin quiver at $q=0$ is isomorphic to the composition
monoid. Finally, in chapter 7 we show the result on the composition algebra at $q=0$ and
the composition monoid for the extended Dynkin case.
\cleardoublepage
\subsection*{Acknowledgements}
First of all, I wish to express my sincere thanks to my supervisors Andrew Hubery and Henning Krause.
In particular, I would like to thank Andrew for coming up with the topic of the thesis, for always answering my questions
and for being open for my ideas and Henning for his continuing support and his encouragement.
I also thank the representation theory group in Paderborn and
their various guests for their mathematical and non-mathematical support.
In particular, I am grateful to Karsten Dietrich, Claudia K\"ohler and Torsten Wedhorn for
proofreading parts of this thesis. I thank the IRTG ``Geometry and Analysis of Symmetries'' for
providing financial support during these three years and the DAAD for funding my research stay in Lyon
which Philippe Caldero made so enjoyable. Last but not least, I would like to thank my parents
and Marta for their unconditional support.

\chapter{Preliminaries}
\label{chap:prelims}
\dictum[Douglas Adams]{In the beginning the Universe was created. This has made a lot
of people very angry and has been widely regarded as a bad move.}
\section{Quivers, Path Algebras and Root Systems}
For standard notations and results about quivers we refer the reader to
\cite{Ringel_integral}.
A \emph{quiver}\index{quiver} $Q=(Q_0, Q_1, s, t)$ is a directed graph with a set of \emph{vertices}\index{quiver!vertices} $Q_0$, a set of
\emph{arrows}\index{quiver!arrow} $Q_1$ and maps
$s,t \colon Q_1 \rightarrow Q_0$, sending an arrow to its starting respectively terminating vertex.
In particular, we write $\alpha \colon s(\alpha) \rightarrow t(\alpha)$ for an $\alpha \in Q_1$. A quiver $Q$ is finite
if $Q_0$ and $Q_1$ are finite sets. In the following, all quivers will be finite.
A quiver
is connected if its underlying graph is connected.
For a quiver $Q=(Q_0, Q_1, s ,t)$ we define $Q^{op} := (Q_0, Q_1, t, s)$ as the quiver with all arrows reversed.

A \emph{path} of length $r\ge0$ in $Q$ is a sequence of arrows $\xi = \alpha_1 \alpha_2 \cdots \alpha_r$ such that
$t(\alpha_i) = s(\alpha_{i+1})$ for all $1 \le i < r$. We write $t(\xi) := t(\alpha_r)$ and $s(\xi):= s(\alpha_1)$. Pictorially,
if $i_j = s(\alpha_j) = t(\alpha_{j-1})$, then
\[
\makeatletter%
\let\ASYencoding\f@encoding%
\let\ASYfamily\f@family%
\let\ASYseries\f@series%
\let\ASYshape\f@shape%
\makeatother%
\setlength{\unitlength}{1pt}
\includegraphics{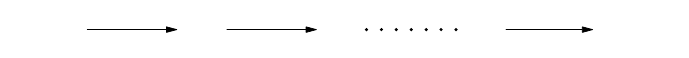}%
\definecolor{ASYcolor}{gray}{0.000000}\color{ASYcolor}
\fontsize{12.000000}{14.400000}\selectfont
\usefont{\ASYencoding}{\ASYfamily}{\ASYseries}{\ASYshape}%
\ASYalign(-181.157906,8.072525)(-0.500000,-0.406901){$i_1$}
\ASYalign(-140.807001,8.072525)(-0.500000,-0.406901){$i_2$}
\ASYalign(-100.456095,8.072525)(-0.500000,-0.406901){$i_3$}
\ASYalign(-60.105189,8.072525)(-0.500000,-0.406901){$i_{r}$}
\ASYalign(-24.210965,8.072525)(0.000000,-0.374607){$i_{r+1}.$}
\ASYalign(-184.771406,8.072525)(-1.000000,-0.390626){$\vphantom{i_r}\xi\colon\ $}
\ASYalign(-160.982454,11.686025)(-0.500000,0.000000){$\scriptstyle \alpha_1$}
\ASYalign(-120.631548,11.686025)(-0.500000,0.000000){$\scriptstyle \alpha_2$}
\ASYalign(-41.409697,10.627655)(0.000000,0.000000){$\scriptstyle \alpha_{r}$}
\]
Clearly, the paths of length one are exactly the arrows of $Q$. For each $i \in Q_0$ there is the
trivial path $\epsilon_i$ of length zero starting and terminating in the vertex $i$. We denote by $Q(i,j)$ the
set of paths starting at the vertex $i$ and terminating at the vertex $j$.

If $R$ is a ring, then the \emph{path algebra} $RQ$ has as basis the set of
paths and the multiplication $\xi \cdot \zeta$ is given by
the concatenation of paths if $t(\xi)=s(\zeta)$ or zero otherwise. In particular, the $\epsilon_i$ are pairwise orthogonal
idempotents of $RQ$, that is $\epsilon_i \epsilon_j = \delta_{ij} \epsilon_i$. The path algebra is obviously
an associative,
unital algebra with unit $1=\sum_{i\in Q_0} \epsilon_i$.

Let $Q_r$ denote the set of paths of length $r$. This extends the notation of vertices $Q_0$ and arrows $Q_1$. Then
\[ RQ = \bigoplus_{r\ge 0} RQ_r,\]
where $RQ_r$ is the free $R$-module with basis the elements of $Q_r$. By construction,
\[RQ_r RQ_s = RQ_{r+s},\]
thus
$RQ$ is an $\N$-graded $R$-algebra.

From now on let $K$ be a field.
We have the following.
\begin{lmm}
	The $\epsilon_i$ form a complete set of pairwise inequivalent orthogonal primitive idempotents in $KQ$. In
	particular, $KQ_0  = \prod_{i \in Q_0} K\epsilon_i$ is a semisimple algebra and the modules $\epsilon_i KQ$ are
	pairwise non-isomorphic indecomposable projective modules.
\end{lmm}

The \index{root lattice}root lattice $\Z Q_0$ is the free
abelian group on elements $\epsilon_i$ for $i \in Q_0$. We define a partial order on
$\Z Q_0$ by $\dvec{d} = \sum_{i} d_i \epsilon_i \ge 0$ if and only if $d_i \ge 0$ for all $i\in Q_0$.
An element $\dvec{d} \in \Z Q_0$ is called a dimension vector.
We endow $\Z Q_0$ with a bilinear form $\bform{\cdot}{\cdot}_Q$ defined by
\[
\bform{\dvec{d}}{\dvec{e}}_Q := \sum_{i \in Q_0} d_i e_i - \sum_{\mathclap{\alpha \colon i \rightarrow j \in Q_1}} d_i e_j.
\]
This form is generally called the \index{Euler form}\emph{Euler form} or the Ringel form.
We also define its symmetrisation
\[
\sbform{\dvec{d}}{\dvec{e}}_Q :=\bform{\dvec{d}}{\dvec{e}}_Q + \bform{\dvec{e}}{\dvec{d}}_Q.
\]

For each vertex $a \in Q_0$ we have the reflection\index{reflection!on dimension vectors}
\begin{alignat*}{2}
 \sigma_a &\colon\quad & \Z Q_0 & \rightarrow \Z Q_0\\
 & & \dvec{d} &\mapsto \dvec{d} - \sbform{\dvec{d}}{\epsilon_a} \epsilon_a.
\end{alignat*}
If $Q$ has no loop at $a$, one easily checks that $\sigma_a^2 \dvec{d} = \dvec{d}$.
By definition, the \emph{Weyl group}\index{Weyl group} $\mathcal{W}$ is the group generated
by the simple reflections
$\sigma_a$, which are those reflections
corresponding to vertices $a \in Q_0$ without loops. Clearly, $\sbform{\cdot}{\cdot}$ is
$\mathcal{W}$-invariant.

We also define reflections on the quiver\index{reflection!on a quiver} itself. The quiver $\sigma_a Q$ is
obtained from $Q$ by reversing all arrows ending or starting in $a$. If $\alpha\colon a \rightarrow j$
is an arrow in $Q_1$, then we call $\alpha^* : j \rightarrow a$ the arrow in the other direction and
analogously for $\alpha \colon i \rightarrow a$. We have that $(\sigma_a Q)_0 = Q_0$, therefore we can
regard $\sigma_a \dvec{d}$ as a dimension vector on $\sigma_a Q$.
Obviously, $\sigma_a^2 Q =Q$ if we identify $\alpha^{* *}$ with $\alpha$ (and we will do
this in the remainder).

A vertex $i \in Q_0$ is called a \emph{sink}\index{sink} (resp. \emph{source}\index{source}) if
there are no arrows starting (resp. terminating) in
$i$. An ordering $(i_1, \dots, i_n)$ of the vertices of $Q$ is called \emph{admissible}\index{admissible ordering},
if
$i_p$ is a sink in $\sigma_{i_{p-1}} \dotsm \sigma_{i_1} Q$ for each $1 \le p \le n$.
Note that there is an admissible ordering if and only if
$Q$ has no oriented cycles. Such a quiver will be called \emph{acyclic}\index{acyclic}\index{quiver!acyclic}.
\begin{beispiel}
Let
\[ Q =
\begin{array}{@{}c@{}}
\makeatletter%
\let\ASYencoding\f@encoding%
\let\ASYfamily\f@family%
\let\ASYseries\f@series%
\let\ASYshape\f@shape%
\makeatother%
\setlength{\unitlength}{1pt}
\includegraphics{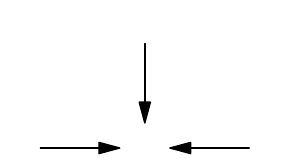}%
\definecolor{ASYcolor}{gray}{0.000000}\color{ASYcolor}
\fontsize{12.000000}{14.400000}\selectfont
\usefont{\ASYencoding}{\ASYfamily}{\ASYseries}{\ASYshape}%
\ASYalign(-43.495099,5.372285)(-0.500000,-0.500000){$1$}
\ASYalign(-80.915173,5.372285)(-0.500000,-0.500000){$2$}
\ASYalign(-43.495099,42.792359)(-0.500000,-0.500000){$3$}
\ASYalign(-6.075025,5.372285)(-0.500000,-0.500000){$4.$}
\ASYalign(-62.205136,8.985785)(-0.500000,0.000000){$\scriptstyle \alpha$}
\ASYalign(-47.108599,24.082322)(-1.000000,-0.281250){$\scriptstyle \beta$}
\ASYalign(-24.785062,8.985785)(-0.500000,0.311112){$\scriptstyle \gamma$}
\end{array}
\]
Then
\[ \sigma_1 Q =
\begin{array}{@{}c@{}}
\makeatletter%
\let\ASYencoding\f@encoding%
\let\ASYfamily\f@family%
\let\ASYseries\f@series%
\let\ASYshape\f@shape%
\makeatother%
\setlength{\unitlength}{1pt}
\includegraphics{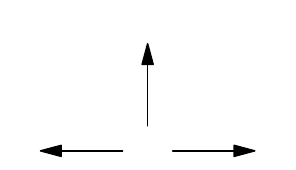}%
\definecolor{ASYcolor}{gray}{0.000000}\color{ASYcolor}
\fontsize{12.000000}{14.400000}\selectfont
\usefont{\ASYencoding}{\ASYfamily}{\ASYseries}{\ASYshape}%
\ASYalign(-42.679134,5.372285)(-0.500000,-0.500000){$1$}
\ASYalign(-80.915173,5.372285)(-0.500000,-0.500000){$2$}
\ASYalign(-42.679134,43.608324)(-0.500000,-0.500000){$3$}
\ASYalign(-4.443095,5.372285)(-0.500000,-0.500000){$4$}
\ASYalign(-61.797153,8.985785)(-0.500000,0.000000){$\scriptstyle \alpha^*$}
\ASYalign(-39.065634,24.490304)(0.000000,-0.283063){$\scriptstyle \beta^*$}
\ASYalign(-23.561114,8.985785)(-0.500000,0.216937){$\scriptstyle \gamma^*$}
\end{array}
\qquad
\text{and}
\qquad
\sigma_2 Q =
\begin{array}{@{}c@{}}
\makeatletter%
\let\ASYencoding\f@encoding%
\let\ASYfamily\f@family%
\let\ASYseries\f@series%
\let\ASYshape\f@shape%
\makeatother%
\setlength{\unitlength}{1pt}
\includegraphics{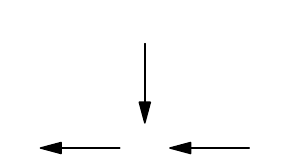}%
\definecolor{ASYcolor}{gray}{0.000000}\color{ASYcolor}
\fontsize{12.000000}{14.400000}\selectfont
\usefont{\ASYencoding}{\ASYfamily}{\ASYseries}{\ASYshape}%
\ASYalign(-43.495099,5.372285)(-0.500000,-0.500000){$1$}
\ASYalign(-80.915173,5.372285)(-0.500000,-0.500000){$2$}
\ASYalign(-43.495099,42.792359)(-0.500000,-0.500000){$3$}
\ASYalign(-6.075025,5.372285)(-0.500000,-0.500000){$4.$}
\ASYalign(-62.205136,8.985785)(-0.500000,0.000000){$\scriptstyle \alpha^*$}
\ASYalign(-47.108599,24.082322)(-1.000000,-0.281250){$\scriptstyle \beta$}
\ASYalign(-24.785062,8.985785)(-0.500000,0.311112){$\scriptstyle \gamma$}
\end{array}\]
Let $\dvec{d} = (2,1,1,1)$.
Then $\sigma_1 \dvec{d} = (2,1,1,1) - (2 + 2 - (1 + 1 + 1)) \epsilon_1 = (1,1,1,1)$ and
$\sigma_2 \dvec{d} = (2,1,1,1) - (1 + 1 - 2) \epsilon_2 = (2,1,1,1)$.
We have that $(1, 2, 3, 4)$ is an admissible ordering of $Q$.
\end{beispiel}

We can define the
set $\Delta \subset \Z Q_0$ of roots of $Q$ combinatorially as follows. We have the
set of \emph{simple roots}\index{simple roots}\index{roots!simple}
\[
\Pi := \Set{ \epsilon_i | i \in Q_0, \text{ no loop at } i}.
\]
The \emph{fundamental region}\index{fundamental region}\index{roots!fundamental region} is
\[
F := \Set{ \dvec{d} > 0 | \sbform{\dvec{d}}{\epsilon_i} \le 0 \text{ for all } \epsilon_i \in \Pi \text{ and } \supp \dvec{d} \text{ connected}},
\]
where $\supp \dvec{d}$ is the \emph{support}\index{support!of a dimension vector}
of $\dvec{d}$, i.e. the full subquiver of $Q$ on the vertices $i$ such that $d_i \neq 0$.
\begin{remark}
  The fundamental region may be empty.
\end{remark}

The sets of \emph{real} and \emph{imaginary roots}\index{roots!real}\index{roots!imaginary} are now
\begin{xalignat*}{3}
    \Delta^{re} & := \mathcal{W} \cdot \Pi & &\text{and} &\Delta^{im} & := \pm \mathcal{W}\cdot  F.
\end{xalignat*}
We define the set of roots $\Delta := \Delta^{re} \cup \Delta^{im}$.
Moreover, each root is either positive or negative and we write $\Delta_+$ for the set of
\emph{positive roots}\index{roots!positive}.
Let $\dvec{d}$ be a root. We note that $\bform{\dvec{d}}{\dvec{d}} = 1$ if $\dvec{d}$ is real and
$\bform{\dvec{d}}{\dvec{d}} \le 0$ if $\dvec{d}$ is imaginary. We call a root $\dvec{d}$ isotropic
if $\bform{\dvec{d}}{\dvec{d}} = 0$.

\section{Representations of Quivers}
For notation and results about representations of quivers and algebras we refer to
\cite{Ringel_integral} and for the geometric aspects to \cite{Kac_infiniterootsystems} and \cite{Schofield_genrep}.
Let $Q$ be a quiver and $K$ be a field. A $K$-\emph{representation}\index{representation}\index{quiver!representation}
$M$ of $Q$ is given by finite dimensional $K$-vector spaces
$M_i$ for each $i\in Q_0$ and $K$-linear maps $M_\alpha \colon M_i \rightarrow M_j$ for each $\alpha \in Q_1$.
If $M$ and $N$ are $K$-representations of $Q$, then a morphism $f \colon M \rightarrow N$ is given by
$K$-linear maps
$f_i \colon M_i \rightarrow N_i$ for each $i \in Q_0$ such that the following diagram commutes
\[
\xymatrix{
M_i \ar[r]^{M_\alpha} \ar[d]_{f_i} & M_j \ar[d]^{f_j}\\
N_i \ar[r]_{N_\alpha} & N_j
}
\]
for all $\alpha \colon i \rightarrow j \in Q_1$.

For each vertex $i \in Q_0$ there is a simple $K$-representation $S_i$ given by setting
$(S_i)_i := K$, $(S_i)_j := 0$ for $i \neq j \in Q_0$ and $(S_i)_\alpha := 0$ for all
$\alpha \in Q_1$.

A $K$-representation $M$ is called nilpotent if it has a filtration by semisimples involving
only the simples $S_i$, $i\in Q_0$. If not otherwise stated, every representation will assumed
to be nilpotent.

The direct sum of two representations is given by
$(X \oplus Y)_i := X_i \oplus Y_i$ and $(X \oplus Y)_\alpha := X_\alpha \oplus Y_\alpha$ and a representation is called
\emph{indecomposable}\index{representation!indecomposable}\index{indecomposable} if it is non-zero and
not isomorphic to a proper direct sum of two representations. We obtain an
additive category $\repK{Q}{K}$ in which the \emph{Krull-Remak-Schmidt theorem}\index{Krull-Remak-Schmidt theorem} holds.
\begin{theom}
	Every representation is isomorphic to a direct sum of indecomposable representations and the isomorphism classes
	and multiplicities are uniquely determined.
	\label{theom:KRS}
\end{theom}

The \emph{dimension vector}\index{representation!dimension vector}\index{dimension vector} of $M$ is defined by
\[
\dvec{\dim}_K M := \sum_{i \in Q_0} \dim_K M_i \epsilon_i \in \Z Q_0.
\]
By fixing bases, we see that
the representations of dimension vector $\dvec{d}$ (probably also non-nilpotent) are parametrised by the affine space
\[
\RepK[Q]{\dvec{d}}{K} = \Rep{\dvec{d}} := \bigoplus_{\alpha \colon i \rightarrow j} \Hom_K(K^{d_i}, K^{d_j}).
\]
In the following we identify each point $x\in\RepK[Q]{\dvec{d}}{K}$ with the corresponding $K$-representation
$M_x \in \repK{Q}{K}$
and write directly $M \in \RepK[Q]{\dvec{d}}{K}$.

The group
\[
\GL_\dvec{d}=\GL_\dvec{d}(K) := \prod_{i \in Q_0} \GL(K^{d_i})
\]
acts on $\RepK[Q]{\dvec{d}}{K}$ by conjugation: given $g = (g_i) \in \GL_\dvec{d}$ and $M \in \RepK[Q]{\dvec{d}}{K}$ we define
$(g \cdot M)_\alpha := g_j M_\alpha g_i^{-1}$ for each $\alpha \colon i \rightarrow j \in Q_1$, i.e. making the following diagram commute
\[
\xymatrix{
K^{d_i} \ar[r]^{M_\alpha} \ar[d]_{g_i} & K^{d_j} \ar[d]^{g_j}\\
K^{d_i} \ar[r]_{(g\cdot M)_\alpha} & K^{d_j}.
}
\]
There is a 1-1 correspondence between the $\GL_\dvec{d}$-orbits in $\RepK[Q]{\dvec{d}}{K}$ and isomorphism classes of
$K$-representations of $Q$ of dimension vector $\dvec{d}$. Denote the orbit of a representation $M$ under this action by
$\Orbit_M$.

A $K$-representation $X$ gives rise to a $KQ$-module $M:=\bigoplus X_i$
where the action of $\alpha\colon i \rightarrow j$ is
given by $\iota_j X_\alpha \pi_i$, $\iota_j \colon X_j \rightarrow \bigoplus X_k = M$ denoting the
canonical inclusion and $\pi_i \colon M= \bigoplus X_k \rightarrow X_i$ the canonical
projection. Vice versa, each $KQ$-module $M$ gives rise to a $K$-representation $X$ with
$X_i := M\epsilon_i$ and $X_\alpha$ is the restriction of the multiplication with $\alpha\colon i \rightarrow j$ to the
domain
$X_i$ and the codomain $X_j$. This yields an equivalence of $\Mod KQ$ with
$\repK{Q}{K}$. Therefore,
$\repK{Q}{K}$ is abelian and we can speak of kernel, cokernel, image and exactness.
Vector space duality $D=\Hom_K(-,K)$ gives a duality $D \colon \Mod KQ \rightarrow \Mod KQ^{op}$.

We will use the following notations for a $K$-algebra $\Lambda$ and two finite dimensional
$\Lambda$-modules $M$ and $N$:
\begin{itemize}
  \item $(M,N)_\Lambda := \Hom_\Lambda(M,N)$,
  \item $(M,N)^i_\Lambda := \Ext^i_\Lambda(M,N)$,
  \item $[M,N]_\Lambda := \dim_K\Hom_\Lambda(M,N)$,
  \item $[M,N]^i_\Lambda := \dim_K\Ext^i_\Lambda(M,N)$.
\end{itemize}
Note that $[M,N]^0 = [M,N]$ and $(M,N)^0 = (M,N)$.
If $Q$ is a quiver, we define $(M,N)_Q := (M,N)_{KQ}$ if $M$ and $N$ are $K$-representations and similarly for
the other notations. For dimension vectors $\dvec{d}$ and $\dvec{e}$ we denote
by $\hom_{KQ}(\dvec{d},\dvec{e})$ the minimal value of $[M,N]_Q$ for $M \in \RepK[Q]{\dvec{d}}{K}$ and
$N \in \RepK[Q]{\dvec{e}}{K}$ and similarly for $\ext^i_{KQ}(\dvec{d},\dvec{e})$. More generally, if
$\mathcal{A}$ and $\mathcal{B}$ are subsets of $\RepK[Q]{\dvec{d}}{K}$ and $\RepK[Q]{\dvec{e}}{K}$ respectively, then
$\hom(\mathcal{A}, \mathcal{B})$ and $\ext^i(\mathcal{A}, \mathcal{B})$ are defined analogously.
Whenever the algebra, the field or the quiver is clear from the context, we omit them from the notation.
If our algebra is hereditary we denote $\Ext^1$ by $\Ext$.

We can generalise the notion of Euler form to $K$-algebras $\Lambda$ of finite global dimension. Namely,
for two $\Lambda$-modules $M$ and $N$ we define
\[ \bform{M}{N}_\Lambda := \sum_{i \ge 0} (-1)^{i} [M,N]^i. \]
By appendix \ref{app:tensalg} we obtain that $KQ$ is hereditary and that for two
$K$-representations $M$ and $N$ of $Q$ we have that
\[ \bform{M}{N}_Q = \bform{M}{N}_{KQ} = [M, N] - [M, N]^1. \]

\section{Reflection Functors}
\label{sec:intro:reflfun}
The main reference for this section is \cite{Ringel_integral}. For a nice
introduction see \cite{HKrause_repsofalgs}.
If $a$ is a sink of $Q$, we define for each $K$-representation $M$ of $Q$ the homomorphism
\[
\phi^M_a \colon \bigoplus_{\alpha: j \rightarrow a} M_j \xrightarrow{\left(
\begin{smallmatrix}
	M_\alpha
\end{smallmatrix}\right)} M_a.
\]
Dually, if $b$ is a source of $Q$, we define
\[
\phi^M_b \colon M_b \xrightarrow{\left(
\begin{smallmatrix}
	    M_\alpha
	\end{smallmatrix}\right)} \bigoplus_{\alpha: b \rightarrow j} M_j .
\]
Note that $D\phi^{DM}_a = \phi^M_a$ and
that $d_a - \rank \phi^X_a = \dim \Hom(X, S_a)$ for $a$ a sink of $Q$ and
a representation $X$ of dimension vector $\dvec{d}$.

We define a pair of \emph{reflection functors}\index{reflection functors}
$S_a^+$ and $S_b^-$. To this end we fix a $K$-representation $M$ of $Q$ of dimension vector $\dvec{d}$.
\begin{enumerate}
	\item[(1)] If the vertex $a$ is a sink of $Q$ we construct
        \[ S_a^+ \colon \repK{Q}{K} \rightarrow \repK{\sigma_a Q}{K} \]
		as follows. We define $S_a^+ M := N$ by letting $N_i := M_i$ for a vertex $i \neq a$ and letting $N_a$ be
		the kernel of the map $\Phi^M_a$. Denote by
		$\iota \colon \Ker \phi^M_a \rightarrow \bigoplus M_j$ the canonical inclusion and by
		$\pi_i \colon \bigoplus M_j \rightarrow M_i$ the canonical projection.
		Then, for each
		$\alpha: i \rightarrow a$ we let $N_{\alpha^*} \colon \Ker \phi^M_a \rightarrow M_i$
		be
		the composition $\pi_i \circ \iota$ making the following diagram commute
		\[
		\begin{CD}
			0 @>>>  \Ker \phi^M_a  @>>\iota> \bigoplus M_j @>>> M_a\\
			@.     @VV{N_{\alpha^*}}V  @VV{\pi_i}V @.\\
			@.    M_i @=       M_i. \\
		\end{CD}
		\]
		We obtain
		a $K$-representation $N$ of $\sigma_a Q$. We call this representation $S^+_a M$.

		Let $s := d_a - \rank \phi^M_a$
		the codimension
		of $\Bild \phi^M_a$ in $M_a$.
		Obviously, $\dim (S^+_a M)_a = e_a$
		where
		\[e_a = \sum_{\alpha: i \rightarrow a} d_i - \rank \Phi^M_a  =
		d_a - (\dvec{d}, \epsilon_a) + s = d_a -s  - (\dvec{d}- s \epsilon_a, \epsilon_a).\]
		Therefore, $\dimv{ S^+_a M } = \sigma_a (\dvec{d} - s \epsilon_a) = \sigma_a (\dvec{d}) + s \epsilon_a$.

		For a morphism $f = (f_i) \colon X \rightarrow Y$ we obtain a morphism
		$S_a^+ f \colon S_a^+ X \rightarrow S_a^+ Y$ by letting
		$(S_a^+ f)_i = f_i$ for $i \neq a$ and letting $(S_a^+ f)_a$ be the map induced on
		the kernels of $\Phi^X_a$ respectively $\Phi^Y_a$.

		If
		\[0 \rightarrow X^1 \rightarrow X^2 \rightarrow X^3 \rightarrow 0\]
		is a short exact sequence in $\repK{Q}{K}$, then we obtain an exact sequence
		\[ 0 \rightarrow S_a^+ X^1 \rightarrow S_a^+ X^2 \rightarrow S_a^+ X^3 \rightarrow S_a^{s} \rightarrow 0\]
        in $\repK{\sigma_a Q}{K}$
		where $s=s_1 - s_2 + s_3$ with $s_i= \dim X^i_a - \rank \Phi^{X^i}_a$.
	\item[(2)] If the vertex $b$ is a source of $Q$, we construct
        \[ S_b^- \colon \repK{Q}{K} \rightarrow \repK{\sigma_b Q}{K} \]
		dually.

		Let $s := d_b - \rank \phi^M_b$
		the dimension
		of $\ker \phi^M_b \subset M_b$.
		We have that $\dim (S_b^-M)_b = e_b$
		where
		\[e_b = \sum_{\alpha: b \rightarrow i} d_i - \rank \phi^M_b  =
		d_b - (\dvec{d}, \epsilon_b) + s = d_b -s  - (\dvec{d}- s \epsilon_b, \epsilon_b).\]
		Therefore, $\dimv{ S^-_b M } = \sigma_b (\dvec{d} - s \epsilon_b) = \sigma_b (\dvec{d}) + s \epsilon_b$.

%
		If
		\[0 \rightarrow X^1 \rightarrow X^2 \rightarrow X^3 \rightarrow 0\]
        is a short exact sequence in $\repK{Q}{K}$, then we obtain an exact sequence
		\[ 0 \rightarrow S_b^s \rightarrow S_b^- X^1 \rightarrow S_b^- X^2 \rightarrow S_b^- X^3  \rightarrow 0\]
        in $\repK{\sigma_b Q}{K}$ where $s=s_1 - s_2 + s_3$ and $s_i= \dim X^i_b - \rank \Phi^{X^i}_b$.
\end{enumerate}
For a sink $a$ of $Q$ we have that $(S_a^-, S_a^+)$ is a pair of adjoint functors and
that $S_a^+$ is left exact and $S_a^-$ is right exact. There is a natural
monomorphism $\iota_{a,M} \colon S_a^- S_a^+ M \rightarrow M$ for $M \in \repK{Q}{K}$ and a natural epimorphism
$\pi_{a,N} \colon N \rightarrow S_a^+ S_a^- N$ for $N \in \repK{\sigma_a Q}{K}$.
We have the following lemma.
\begin{lmm}
 Let $a$ be a sink and $X$ an indecomposable representation of $Q$. Then, the following are equivalent:
 \begin{enumerate}
 \item $X \ncong S_a$.
 \item $S^+_a X$ is indecomposable.
 \item $S^+_a X \neq 0$.
 \item $S_a^- S_a^+ X \cong X$ via the natural inclusion.
 \item The map $\Phi_a^X$ is an epimorphism.
 \item $\sigma_a (\dvec{\dim} X) > 0$.
 \item $\dvec{\dim} S^+_a X = \sigma_a (\dvec{\dim} X)$.
 \end{enumerate}
\end{lmm}

For any admissible sequence of sinks $w=(i_1, \dots, i_r)$ define
\[S_w^+ := S_{i_r}^+ \circ \dots \circ S_{i_1}^+. \]

If $w=(i_1, \dots, i_n)$ is an admissible ordering of $Q$, we define the \emph{Coxeter functors}\index{Coxeter functors}
\begin{align*}
  C^+ &:= S^+_{i_n} \circ \dots \circ S^+_{i_1} & C^- &:= S^-_{i_1} \circ \dots \circ S^-_{i_n}. \\
\end{align*}
Since both reverse each arrow of $Q$ exactly twice, these are endofunctors of $\repK{Q}{K}$. Neither functor
depends on the
choice of the admissible ordering.

A $K$-representation $P$ is projective if and only if $C^+ P =0$. Dually, a $K$-representation $I$
is injective if and only if $C^- I =0$. This motivates the nomenclature in the following definition.
\begin{defn}
  An indecomposable $K$-representation $M$ of $Q$ is called \emph{preprojective}\index{preprojective} if $(C^+)^r M=0$ for some $r\ge0$
  and \emph{preinjective}\index{preinjective} if $(C^-)^r M =0$ for some $r\ge0$. If $M$ is neither preprojective nor preinjective, then
  $M$ is called \emph{regular}\index{regular}.

  An arbitrary $K$-representation is called preprojective (or preinjective or regular) if it is isomorphic to a direct sum of indecomposable
  preprojective (or preinjective or regular) representations.
\end{defn}
Let $Q$ be a connected, acyclic non-Dynkin quiver as introduced in the next section.
Then an arbitrary representation $M$ can be decomposed uniquely into a direct sum $M \cong M_P \oplus M_R \oplus M_I$ such
that $M_P$ is preprojective, $M_R$ is regular and $M_I$ is preinjective. This means that the set of isomorphism classes of
indecomposable $K$-representations $\ind \repK{Q}{K}$ decomposes into a disjoint union $\mathcal{P} \cup \mathcal{R} \cup \mathcal{I}$
where $\mathcal{P}$ denotes the set of indecomposable preprojective, $\mathcal{R}$ denotes the set of indecomposable regular
and $\mathcal{I}$ denotes the set of indecomposable preinjective $K$-representations. Note that
$\Hom(R,P) = \Hom(I,P) = \Hom(I,R) = 0$ and $\Ext(P,R)=\Ext(P,I)=\Ext(R,I)=0$ for all representations
$P \in \mathcal{P}$, $R \in \mathcal{R}$ and $I \in \mathcal{I}$.
There is a partial order $\preceq$ on $\mathcal{P} \cup \mathcal{I}$ given by
$M \preceq N$ for $M, N \in \mathcal{P} \cup \mathcal{I}$ if and only if
there is a sequence of non-zero morphisms $M \rightarrow M_1 \rightarrow M_2 \rightarrow \dots \rightarrow N$
for some indecomposable representations $M_i \in \mathcal{P} \cup \mathcal{I}$.
Note that, if $\Ext^1(M, N) \neq 0$ for two indecomposables $M, N \in \mathcal{P} \cup \mathcal{I}$, then
$N \prec M$ .

Fix an admissible ordering $(a_1, \dots, a_n)$ of the vertices of $Q$.
  For each indecomposable preprojective representation $M$ there is a natural number
  $r=kn+s$ for some $k\ge0$ and $0 \le s < n$ such that
  \[ S^+_{a_{s}} \dotsm S^+_{a_1}(C^+)^k M =0. \]
  Let $\sigma(M)$ be the minimal such number. Note that $\sigma(M)$ depends on the
  choice of the admissible ordering.
We have the following easy
\begin{lmm}
	\label{lmm:intro_refl:order}
	The map
    \begin{alignat*}{2}
		\sigma &\colon\quad & \mathcal{P} & \rightarrow \N\\
		&& M & \mapsto \sigma(M)
    \end{alignat*}
	is an injection respecting the partial order $\preceq$ on $\mathcal{P}$,
	i.e. for all $M, N \in \mathcal{P}$ with $M \preceq N$ we have that
	$\sigma(M) \le \sigma(N)$ for the natural ordering on $\N$.
\end{lmm}
\begin{proof}
    Choosing an admissible ordering is the same as refining the partial order $\preceq$ on the projectives
    in $\repK{Q}{K}$ to a total order. Since the Auslander-Reiten quiver is just
    a number of copies of the quiver given by the projectives with morphisms only going from
    left to right, the lemma follows.
\end{proof}

\section{Dynkin and Extended Dynkin Quivers}
The references for this chapter are \cite{Ringel_integral} and \cite{Gabriel_unzdst}.
The representation type of a quiver $Q$ is governed by its underlying graph $\Gamma$. Note
that the symmetric bilinear form $\sbform{\cdot}{\cdot}_Q$ depends only on $\Gamma$.

\begin{theom}
    \label{dynkineuclwild} Suppose $\Gamma$ is connected.
    \begin{enumerate}
        \item $\Gamma$ is \emph{Dynkin}\index{Dynkin!diagram}
            if and only if $\sbform{\cdot}{\cdot}_Q$ is positive definite. By
            definition the (simply-laced) Dynkin diagrams are:
            \shorthandoff{"}
            \begin{align*}
                A_n: &
                {\xygraph{ !{0;(.7,0):0} []\bullet -[r]\bullet - [r] \bullet -[r] -@{ *{.\ } }
                [r] - [r]\bullet } }  \quad (n\ge1)&
                E_6: & {\xygraph{!{0;(.7,0):0} []\bullet -[r]\bullet -[r] \bullet="B" -[r] \bullet -[r] \bullet
                "B" -[u] \bullet } } \\
                D_n: &
                {\xygraph{ !{0;(.7,0):0} [] \bullet="B" -[r] \bullet -[r] -@{ *{.\:} }
                [r] - [r]\bullet "B" -[dl(0.77)] \bullet "B" - [ul(0.77)] \bullet } }  \quad\ \; (n\ge4)&
                E_7: & {\xygraph{!{0;(.7,0):0} []\bullet -[r]\bullet -[r] \bullet="B" -[r] \bullet -[r] \bullet
                -[r] \bullet
                "B" -[u] \bullet } } \\
                &(n=\text{number of vertices}) &
                E_8: & {\xygraph{!{0;(.7,0):0} []\bullet -[r]\bullet -[r] \bullet="B" -[r] \bullet -[r] \bullet
                -[r] \bullet -[r] \bullet
                "B" -[u] \bullet } }
            \end{align*}
        \item $\Gamma$ is \emph{extended Dynkin}\index{Dynkin!extended} if and only if $\sbform{\cdot}{\cdot}_Q$ is positive semi-definite and
            not positive definite. We have that
            $\rad\sbform{\cdot}{\cdot}_Q=\Z \dvec{\delta}$ for some dimension vector $\dvec{\delta}$.
            By definition
            the extended Dynkin diagrams are as below. We have marked each vertex $i$ with the value of $\delta_i$.
            Note that $\dvec{\delta}$ is sincere and $\dvec{\delta} \ge 0$.
            \begin{align*}
                \widetilde{A}_n:&
                {\quad\xygraph{ !{0;(.7,0):0} [] 1="A" -[ur(.77)] 1 -[r]  -@{ *{.\,} }[r]
                -[r] 1 -[dr(.77)] 1 -[dl(.77)] 1 -[l]  -@{ *{.\,} }[l] -[l] 1 -"A"} } & (n\ge 0)\\
                \widetilde{D}_n: &
                { \quad\xygraph{ !{0;(.7,0):0} [] 2="A" -[r] 2 -@{ *{.\,} }[r] -[r] 2="B" -[ur(.77)] 1
                "B" - [dr(.77)] 1 "A" -[ul(.77)] 1 "A" -[dl(.77)] 1} } & (n \ge 4)\\
                \widetilde{E}_6: &
                {\quad\xygraph{!{0;(.7,0):0} []1 -[r]2 -[r] 3="B" -[r] 2 -[r] 1
                "B" -[u] 2 -[u] 1 } }
                 &\quad (n+1 = \text{number of vertices})\\
                \widetilde{E}_7: &
                {\quad\xygraph{!{0;(.7,0):0} []1 -[r]2 -[r] 3 -[r] 4="B" -[r] 3 -[r] 2 -[r] 1
                "B" -[u] 2 } } \\
                \widetilde{E}_8: &
                {\quad\xygraph{!{0;(.7,0):0} []2 -[r]4 -[r] 6="B" -[r] 5 -[r] 4 -[r] 3 -[r] 2 -[r] 1
                "B" -[u] 3 } }\\
            \end{align*}
            \shorthandon{"}
            Note that $\widetilde{A}_0$ has one vertex and one loop and $\widetilde{A}_1$ has two
            vertices joined by two edges.
          \item Otherwise, there is a dimension vector $\dvec{d} \ge 0$ with $\sbform{\dvec{d}}{\dvec{d}} < 0$ and
            $\sbform{\dvec{d}}{\epsilon_i} \le 0$
            for all $i$.
    \end{enumerate}
\end{theom}
A vertex $i$ of an extended Dynkin graph $\Gamma$ with $\delta_i=1$ is called an extending vertex. Removing
$i$ from $\Gamma$ gives the corresponding Dynkin diagram. We call a quiver of Dynkin type if it
is a disjoint union of quivers with underlying graphs being Dynkin. Similarly for extended Dynkin.

Gabriel's theorem states the following.
\begin{theom}
    If $Q$ is a connected, acyclic quiver with underlying graph $\Gamma$, then there are only finitely many
    isomorphism classes of indecomposable representations of $Q$ if and only if $\Gamma$ is Dynkin.
    In this case the assignment
    $X \mapsto \dimve X$ induces a bijection between the isomorphism classes of indecomposable
    representations and the positive roots $\Delta_+$.
\end{theom}
More generally, Kac's theorem \cite{Kac_infiniterootsystems} yields that
    the map $X \mapsto \dimve X$ is a surjection from isomorphism classes
    of indecomposable representations to the set of positive roots $\Delta_+$.

    Now let us review the representation theory of connected, acyclic, extended Dynkin quivers, the
    main references being \cite{DlabRingel_repsgraphsalgs} and $\cite{Ringel_species}$.
    Let $Q$ be an acyclic extended Dynkin quiver and $K$ a field.
    Let
    $\Gamma_{AR}$ be the Auslander-Reiten quiver of $Q$. Vertices of $\Gamma_{AR}$ correspond
    to indecomposable representations and arrows to irreducible morphisms. We have
    a decomposition of $\Gamma_{AR}$
    into the preprojective, the preinjective and the regular part.
    The set of regular representations
    $\add \mathcal{R}$ is an abelian subcategory of the category of all
    representations. We say that
    $M$ is regular simple, of regular length $k$, \ldots if $M$ is simple,
    of length $k$, \ldots in $\add \mathcal{R}$. The regular part is
    the disjoint union of pairwise orthogonal tubes\index{tube} $\mathcal{T}_x$ for
    each scheme theoretic, closed point $x \in \Proj^1_K = \ProjS K[X,Y]$ in
    such a way that each regular
    simple representation $R$ in the tube labelled by $x$ satisfies $\End(R) \cong \kappa(x)$,
    $\kappa(x)$ denoting the residue field at the point $x$. The
    degree of $x$ is defined to be $[\kappa(x):K]$ as degree of field extensions. We have
    that each regular simple $R$ in the tube $\mathcal{T}_x$ is of dimension vector
    $\dimve R = (\deg x) \dvec{\delta}$. Note that
    if $K$ is algebraically closed, then all tubes have degree one. We define $\rank \mathcal{T}_x$
    to be the number of regular simples
    in the tube $\mathcal{T}_x$. A tube $\mathcal{T}_x$ is called homogeneous if $\rank \mathcal{T}_x =1$ and
    inhomogeneous otherwise.
    Each tube $\mathcal{T}_x$
    is equivalent to the category of (nilpotent) representations of an
    oriented cycle with $\rank \mathcal{T}_x$ vertices.
    For a $K$-representation $M$ we define $M_x$ to be the summand
    of $M$ living in the tube $x \in \Proj^1_K$.
    Finally, we may also assume that the non-homogeneous tubes
    are labelled by some subset of $\{0, 1, \infty\}$, whereas the homogeneous tubes are labelled
    by the points of the scheme $\mathbb{H}_K = \mathbb{H}_\Z \tensor K$ for some open integral subscheme
    $\mathbb{H}_\Z \subset \Proj^1_\Z$.

    For the representation theory of $Q$ the \emph{defect}\index{defect} $\partial$ plays a major role.
    For a dimension vector $\dvec{d}$ we define
    \[
    \partial \dvec{d} := \bform{\dvec{\delta}}{\dvec{d}}_Q.
    \]
    If $M$ is an indecomposable representation, then
    \begin{itemize}
        \item $M$ is preprojective if and only if $\partial \dvec{d} < 0$;
        \item $M$ is preinjective if and only if $\partial \dvec{d} > 0$;
        \item $M$ is regular if and only if $\partial \dvec{d} = 0$.
    \end{itemize}

\section{Canonical Decomposition}
In this section we recall the canonical decomposition of a dimension
vector $\dvec{d}$ introduced by V. Kac \cite{Kac_infiniterootsystems,Kac_infiniterootsystems2} and examined
by A. Schofield \cite{Schofield_genrep}. Let $K$ be an algebraically closed field and let $Q$ be
a quiver. A dimension vector $\dvec{d}$ is called a \emph{Schur root}\index{Schur root}\index{roots!Schur}
if the general representation of dimension vector $\dvec{d}$ has endomorphism ring $K$, i.e.
there is an open non-empty subset $U \subset \Rep[Q]{\dvec{d}}$ such that for all
$M \in U$ we have that $\End(U) \cong K$.
We have that all preprojective and preinjective roots are real Schur roots.

A decomposition
$\dvec{d} = \sum \dvec{f}^i$ is called the \emph{canonical decomposition}\index{canonical decomposition}
if and only if
a general representation of dimension vector $\dvec{d}$ is isomorphic to the direct sum
of indecomposable representations of dimension vectors $\dvec{f}^1, \dots, \dvec{f}^r$.

This is closely related to the question if all representations of a dimension vector
$\dvec{d} + \dvec{e}$ have a subrepresentation of dimension vector $\dvec{d}$.
A. Schofield \cite{Schofield_genrep} proved the following theorem for $\chark K = 0$ and
W. Crawley-Boevey \cite{Bill_subreps} for
$\chark K$ arbitrary.
\begin{theom}
    Every representation of dimension vector
    $\dvec{d} + \dvec{e}$ has a subrepresentation of dimension vector $\dvec{d}$
    if and only if $\ext(\dvec{d}, \dvec{e}) = 0$.
    Moreover, the numbers $\ext(\dvec{d}, \dvec{e})$ are given
    combinatorially and therefore are independent of the algebraically closed base field.
    \label{intro:theom:extallsubreps}
\end{theom}
\vspace{-0.5cm}

\begin{remark}
    Note that, if $K$ is not algebraically closed, then the theorem fails. Let
    $Q$ be the Kronecker quiver. Then there is a regular simple $K$-representation $M$ of
    dimension vector $(2,2)$. The representation $M$ does not have a subrepresentation of dimension
    vector $(1,1)$ even though $\ext( (1,1), (1,1)  ) = 0$ as one easily sees by
    taking two regular simple representations of dimension vector $(1,1)$ living
    in different tubes.
\end{remark}

By using this theorem, A. Schofield showed the following.
\begin{theom}
    A decomposition $\dvec{d} = \sum \dvec{f}^i$ is the canonical decomposition if and only if each $\dvec{f}^i$
    is a Schur root and
    $\ext(\dvec{f}^i, \dvec{f}^j) = 0$ for all $i \neq j$. Moreover, the decomposition
    is independent of the field.
\end{theom}

We often write the canonical decomposition as $\dvec{d} = \sum r_i\dvec{f}^i$ such
that $r_i > 0$ and the $\dvec{f}^i$ are pairwise different Schur roots. Note that
for every Schur root $\dvec{f}^i$ appearing in this sum with multiplicity $r_i >1$
we have that $\ext(\dvec{f}^i, \dvec{f}^i) = 0$. Moreover, for each $i\neq j$ we have that
$\ext(\dvec{f}^i, \dvec{f}^j)=0$.

We have the following.
\begin{lmm}
    Let $\dvec{d}^1, \dots, \dvec{d}^k$ be dimension vectors such that
    $\ext(\dvec{d}^i, \dvec{d}^j) = 0$ for all $i \neq j$. Then the
    canonical decomposition of $\sum \dvec{d}^i$ is a refinement of
    the decomposition given by $\sum \dvec{d}^i$, i.e. there are
    Schur roots $\dvec{f}^i_j$ such that $\sum_j \dvec{f}^i_j$ is the
    canonical decomposition of $\dvec{d}^i$ and
    $\sum_{i,j} \dvec{f}^i_j$ is the canonical decomposition of $\sum \dvec{d}^i$.
    \label{lmm:candecrefine}
\end{lmm}
\begin{proof}
    Let $\dvec{f}^i_j$ be Schur roots
    such that for each $i$ we have that $\sum_j \dvec{f}^i_j$ is the canonical decomposition
    of $\dvec{d}^i$. A general representation of dimension vector $\dvec{d}^i$ is therefore a direct sum of
    indecomposable representations of dimension vectors $\dvec{f}^i_j$.
    By A. Schofield \cite[Theorem 3.4]{Schofield_genrep} a general representation 
    of dimension vector $\sum \dvec{d}^i$ is a direct sum of representations of
    dimension vector $\dvec{d}^i$. Therefore,
    a general representation of dimension vector $\sum \dvec{d}^i$ is a direct sum
    of indecomposable representations of dimension vectors $\dvec{f}^i_j$. Hence,
    $\sum \dvec{f}^i_j$ is the canonical decomposition of $\sum \dvec{d}^i$.
\end{proof}

\section{Degenerations}
The basic reference for the following is \cite{Bongartz_degens}.
Let $K$ be an algebraically closed field. Let $M$ and $N$ be $K$-representations of
a quiver $Q$. We say that
$M \ledeg N$ if $\Orbit_N \subseteq \overline{\Orbit_M}$. The degeneration
order $\ledeg$ on the representation variety has been investigated by various authors.

C. Riedtmann \cite{Riedtmann_degens} and G. Zwara \cite{Zwara_degensext}
were able to describe $\ledeg$ in purely representation theoretic terms.
\begin{theom}
    Let $M$ and $N$ be $K$-representations of $Q$. The following are equivalent:
    \begin{itemize}
        \item $M \ledeg N$.
        \item There is representation $Y$ and a short exact sequence
            \[ 0 \rightarrow Y \rightarrow Y \oplus M \rightarrow N \rightarrow 0. \]
        \item There is representation $Z$ and a short exact sequence
            \[ 0 \rightarrow N \rightarrow M \oplus Z \rightarrow Z \rightarrow 0. \]
    \end{itemize}
\end{theom}

There are two other orders on the isomorphism classes of representations. Firstly, the
$\Hom$-order $\le$ is given by
$M \le N$ if $[M,X] \le [N,X]$ for all representations $X$.
A result of Auslander yields that $\le$ is a partial order. Moreover,
$M \le N$ if and only if $[X,M] \le [X,N]$ for all representations $X$. Therefore,
the definition of $\le$ is symmetric.

Secondly, we have the $\Ext$-order $\le_{\ext}$. We define $M \le_{\ext} N$ if
there are representations $M_i$, $U_i$ and $Q_i$ and short exact
sequences
\[
0 \rightarrow U_i \rightarrow M_i \rightarrow Q_i \rightarrow 0
\]
for $0 \le i < r$ such that $M = M_0$, $M_{i+1} = U_i \oplus Q_i$ and $N = M_r$.
This also yields a partial order and we have the following.
\begin{theom}[\cite{Bongartz_degens}]
  Let $M$ and $N$ be $K$-representations. Then we have the following implications:
  \[ M \le_{\ext} N \Rightarrow M \ledeg N \Rightarrow M \le N. \]
\end{theom}

It is interesting to investigate when the orders agree. We have the following theorem.
\begin{theom}
    If $Q$ is Dynkin or extended Dynkin, then
    \[ M \le_{\ext} N \Leftrightarrow M \ledeg N \Leftrightarrow M \le N. \]
    \label{theom:extishomorder}
\end{theom}
\vspace{-0.7cm}
This theorem is due to K. Bongartz \cite{Bongartz_degens} for Dynkin quivers
and to G. Zwara \cite{Zwara_degensbiserial, Zwara_degensextdyn} for extended Dynkin quivers.

If $K$ is not algebraically closed we have the obvious definitions of $\le$ and
$\le_{\ext}$. We define $M \ledeg N$ if $M\tensor \overline{K} \ledeg N \tensor \overline{K}$,
$\overline{K}$ being the algebraic closure of $K$.

\section{Segre Classes}
\label{sec:intro_segre}
For this chapter the main references are \cite{Bongartz_decomp}, \cite{Bongartz_schichten}
and \cite{Hubery_hallpolys}.
Let
\[
Q=
\begin{array}{@{}c@{}}
\makeatletter%
\let\ASYencoding\f@encoding%
\let\ASYfamily\f@family%
\let\ASYseries\f@series%
\let\ASYshape\f@shape%
\makeatother%
\setlength{\unitlength}{1pt}
\includegraphics{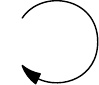}%
\definecolor{ASYcolor}{gray}{0.000000}\color{ASYcolor}
\fontsize{12.000000}{14.400000}\selectfont
\usefont{\ASYencoding}{\ASYfamily}{\ASYseries}{\ASYshape}%
\ASYalign(-23.947121,12.331180)(-0.500000,-0.500000){$\bullet$}
\end{array}
\]
be the Jordan quiver, having one vertex and one loop, and let $K$ be a field.

In order to describe the isomorphism classes of $Q$ we use \emph{partitions}\index{partition}. A partition
$\lambda=(\lambda_1 \ge \lambda_2 \ge \dots\ge \lambda_r)$ of $n$ is a sequence
of decreasing positive integers $\lambda_i \in \N$ such that $n = \sum_{i=1}^r \lambda_i$.
Each $\lambda_i$ is called a part of size $\lambda_i$ of $\lambda$. The length $l(\lambda)$ of
a partition $\lambda$ is the number of parts.
An alternative way to write a partition is in exponential form $\lambda=(1^{l_1}2^{l_2}\dotsm m^{l_m})$ for
an integer $m\in \N$ and non-negative integers $l_1, \dots, l_m \in \N$. This means that
$\lambda$ has exactly $l_i$ parts of size $i$. Therefore, $\lambda$ is a partition of $\sum_i i\cdot l_i$.
\begin{beispiel}
    For example $\lambda=(3,3,2)$ is a partition of $8$. $\lambda$ can be written in exponential form
    as $(1^0 2^1 3^2)$. The partition $\lambda$ has two parts of size $3$ and one part of size $2$.
\end{beispiel}

Now we come back to the Jordan quiver $Q$.
Note
that $\RepK[Q]{d}{K} = \End(K^d)$ for a dimension $d$. We know that $KQ = K[t]$ is a
principal ideal domain, so finite dimensional modules are described by their elementary divisors.
In particular,
we can associate to a finite dimensional module $M$ the data $\{(\lambda_1, p_1), \dots, (\lambda_r, p_r)\}$
consisting of partitions $\lambda_i$ and distinct monic irreducible polynomials $p_i \in K[t]$ such
that
\[ M \cong \bigoplus M(\lambda_i, p_i) \]
where, for a partition $\lambda = (1^{l_1}\dotsm n^{l_n})$ and a monic irreducible polynomial $p$, we write
\[ M(\lambda, p) = \bigoplus_i \left(k[t]/(p^i)\right)^{l_i}.\]

Clearly, the primes $p_i$ depend on the field, but we can partition the set of isomorphism
classes by considering only their degrees. We call this the Segre decomposition.
More precisely, a
\emph{Segre symbol}\index{Segre symbol} is a multiset $\sigma=\{(\lambda_1, d_1), \dots, (\lambda_r, d_r)\}$ of
pairs $(\lambda, d)$ consisting of a partition $\lambda$ and a positive integer $d$. The corresponding
\emph{Segre class}\index{Segre class} $\SegreK(\sigma, K)$ consists of all modules of isomorphism type
$\{(\lambda_1, p_1), \dots, (\lambda_r, p_r)\}$ where the $p_i \in K[t]$ are distinct monic irreducible
polynomials with $\deg p_i = d_i$.
\begin{theom}[\cite{Bongartz_decomp}]
    Let $K$ be an algebraically closed field. Then the Segre classes stratify the variety
    $\RepK[Q]{d}{K}=\End(K^d)$ into
	smooth, irreducible, $\GL_d(K)$-stable subvarieties, each admitting a smooth and rational geometric quotient. Moreover,
	the stabilisers of any two matrices in the same Segre class are conjugate inside $\GL_d(K)$.
\end{theom}

Let $Q$ be a connected, acyclic, extended Dynkin quiver and $K$ a field.
The indecomposable preprojective and preinjective representations are all exceptional, as are the regular simple
representations in the non-homogeneous tubes. Hence, the isomorphism class of a representation without homogeneous regular
summands can be described combinatorially, whereas homogeneous regular representations are determined by pairs
consisting of a partition together with a point of the scheme $\mathbb{H}_K$.

Now we can define the decomposition of K. Bongartz and D. Dudek \cite{Bongartz_decomp} or more
precisely the generalisation given by A. Hubery \cite{Hubery_hallpolys}. A
\emph{decomposition symbol}\index{decomposition symbol} is a
pair $\alpha = (\mu, \sigma)$ such that $\mu$ specifies a representation without
homogeneous regular summands and $\sigma = \{(\lambda_1, d_1), \dots, (\lambda_r, d_r)\}$ is
a Segre symbol. Given a decomposition symbol $\alpha = (\mu, \sigma)$ and a field $K$, we define
the \emph{decomposition class}\index{decomposition class}
$\SegreK_Q(\alpha, K) = \SegreK(\alpha, K)$ to be the set
of representations $X$ of the form $X \cong M(\mu, K) \oplus R$
where $M(\mu, K)$ is the $K$-representation determined by $\mu$ and
\[ R = R(\lambda_1, x_1) \oplus \dots \oplus R(\lambda_r, x_r) \]
for some distinct points $x_1, \dots, x_r \in \mathbb{H}_K$ such that $\deg x_i = d_i$ and $R(\lambda, x)$ is
the representation associated to the partition $\lambda$ living in the tube $\mathcal{T}_x$ of rank one. We call
$\mu$ the discrete part and $\sigma$ the continuous part of $\alpha$. If $\sigma = \emptyset$, we
say that $\alpha$ is discrete.

Let $a$ be a sink of $Q$ and let $\alpha$ be a decomposition symbol. Let $M \in \SegreK_Q(\alpha, K)$ and
let $\beta$ be the decomposition symbol of $S_a^+ M$. Then $S_a^+$ gives a
bijection between isomorphism classes in $\SegreK_Q(\alpha, K)$ and isomorphism classes in $\SegreK_{\sigma_a Q}(\beta, K)$.
This is obvious since
$S_a^+$ is additive and gives a bijection from $\mathcal{R}_Q$ to $\mathcal{R}_{\sigma_a Q}$ with inverse
$S_a^-$ which deals with the continuous part and for the discrete part there is only one choice.

\section{Ringel-Hall Algebra}
For this section the main reference is \cite{Ringel_hallalgsandquantumgroups} or the
lecture notes \cite{Schiffmann_lectureshallalgs} and \cite{Hubery_ringelhall}.
Let $\mathcal{A}$ be a skeletally small abelian category such that for two objects $M, N \in \mathcal{A}$ the sets
$\Ext^i_{\mathcal{A}}(M,N)$ are finite for all $i \ge 0$. Such a category is called \emph{finitary}\index{finitary}.
For three objects $M,N,X \in \mathcal{A}$ define
\[
F^X_{M N} := \#\Set{ U \le X | U \cong N, X/U \cong M }.
\]
Let $\mathcal{H}(\mathcal{A})$ be the $\Q$-vector space with basis $u_{[X]}$ where $[X]$ is
the isomorphism class of $X$. For convenience we write $u_X$ instead of $u_{[X]}$. Define
\[
u_M \diamond u_N := \sum_{X} F^X_{M N} u_X.
\]
Then $(\mathcal{H}(\mathcal{A}), +, \diamond)$ is an associative
$\Q$-algebra with unit $1=u_0$, the \emph{Ringel-Hall algebra}
\index{Ringel-Hall algebra} or just Hall algebra. The \emph{composition algebra}\index{composition algebra}
is the subalgebra $\mathcal{C}(\mathcal{A})$ of $\mathcal{H}(\mathcal{A})$ generated by
the simple objects without self-extensions. Note that $\mathcal{H}(\mathcal{A})$ and
$\mathcal{C}(\mathcal{A})$ are naturally graded by the Grothendieck group of $\mathcal{A}$.

Let $Q$ be a quiver and $q$ a prime power. Then the finite dimensional $\F_q$-representations of $Q$ form a finitary
abelian category. Define
\[ \mathcal{H}_{\F_q}(Q) := \mathcal{H}(\repK{Q}{\F_q})\]
and
\[ \mathcal{C}_{\F_q}(Q) := \mathcal{C}(\repK{Q}{\F_q}).\]
Note that $\mathcal{H}_{\F_q}(Q)$ and $\mathcal{C}_{\F_q}(Q)$ are naturally graded by dimension vector.
We set $u_i := u_{S_i}$ for each $i \in Q_0$. If $w = (i_1, \dots, i_r)$ is
a word in vertices of $Q$, we define
\[ u_w := u_{i_1} \diamond \dots \diamond u_{i_r}.\]
By definition, there are numbers $F_w^X$ for each $\F_q$-representation $X$ of $Q$ such that
\[ u_w = \sum_{X} F_w^X u_X. \]

\section{Generic Composition Algebra}
For the following let $Q$ be a finite quiver with vertex set $Q_0$ and arrow set $Q_1$. We consider only finite
dimensional and nilpotent representations and modules.
We define the \emph{generic composition algebra}\index{generic composition algebra}
\index{composition algebra!generic} via \emph{Hall polynomials}\index{Hall polynomials}.

The main references for this section are C. M. Ringel
\cite{Ringel_hallalgsandquantumgroups, Ringel_hallpolysforrepfiniteheralgs}
for the representation finite case, J. Guo \cite{Guo_hallpolyscyclic} and
C. M. Ringel \cite{Ringel_compalgofcyclic} for the oriented cycle and
A. Hubery \cite{Hubery_hallpolys} for the acyclic, extended Dynkin case.

Now let $Q$ be Dynkin or extended Dynkin.
Then, there is a partition of the isomorphism classes of representations of each dimension vector given
by some combinatorial set $\Segre = \bigcup_\dvec{d} \Segre_\dvec{d}$ such that $\Segre_\dvec{d}$ is finite.
This means, for each
dimension vector $\dvec{d}$ and each field $K$ we have subsets
$\SegreK(\alpha, K) \subset \RepK[Q]{\dvec{d}}{K}$ for each $\alpha \in \Segre_\dvec{d}$ such that
\[ \RepK[Q]{\dvec{d}}{K} = \bigcup_{\alpha \in \Segre_\dvec{d}} \SegreK(\alpha, K) \]
and that there are polynomials $a_\alpha, n_\alpha \in \Q[q]$ such that
for each finite field $K$ we have that $a_\alpha(|K|) = \# \Aut(M)$ for all
$M \in \SegreK(\alpha, K)$ and that $\# [\SegreK(\alpha, K)] = n_\alpha(|K|)$.

Following \cite{Hubery_hallpolys}, we say that Hall polynomials exist with respect to this decomposition
if there are polynomials $f^\gamma_{\alpha \beta} \in \Q[q]$ such that for each finite field $K$
we have that
\[ \sum_{\substack{[A] \in [\SegreK(\alpha, K)]\\ [B] \in [\SegreK(\beta, K)]}} F^C_{A B} = f^{\gamma}_{\alpha \beta}(|K|)
\text{ for all } C \in \SegreK(\gamma, K) \]
and further that
\begin{align*}
    n_\gamma(|K|) f^\gamma_{\alpha \beta}(q) &=
    n_\alpha(|K|) \sum_{\substack{[B] \in [\SegreK(\beta, K)]\\ [C] \in [\SegreK(\gamma, K)]}} F^C_{A B} \quad \text{for all }
    A \in \SegreK(\alpha, K)\\
    &=n_\beta(|K|) \sum_{\substack{[A] \in [\SegreK(\alpha, K)]\\ [C] \in [\SegreK(\gamma, K)]}} F^C_{A B} \quad \text{for all }
    B \in \SegreK(\beta, K).
\end{align*}

For the Dynkin case the isomorphism classes of indecomposable representations are in bijection with the positive roots
$\Delta_+$ of the corresponding Lie algebra, which are independent of the field $K$. Therefore, we can take
$\Segre$ to be the set of all
$\alpha \colon \Delta_+ \rightarrow \N$ with finite support. For each such $\alpha$ and any field $K$ there is a unique
isomorphism class such that the indecomposable representation corresponding to the root $\rho$ appears
$\alpha(\rho)$ times. Choose an element $M(\alpha, K)$ of this isomorphism class.
C. M. Ringel showed that there are polynomials $f^{\gamma}_{\alpha \beta} (q) \in \Z[q]$ such that for each finite
field $K$ we have
\[ F^{M(\gamma,K)}_{M(\alpha,K) M(\beta,K)} = f^{\gamma}_{\alpha \beta} (|K|).\]

If $Q$ is an oriented cycle, we have that Hall polynomials exist if we take $\Segre_\dvec{d}$
to be the set of isomorphism classes of representations of dimension vector $\dvec{d}$ and
$\SegreK(\alpha, K)$ then to be the set of representations $M$ lying in this isomorphism class. See
chapter \ref{chap:cyclic} for more on this.

For $Q$ an acyclic, extended Dynkin quiver we have to take decomposition symbols as defined in section \ref{sec:intro_segre}.
See \cite{Hubery_hallpolys}.

Note that
for all three choices of $\Segre$ we have that each simple representation $S_i$ gives a class $\SegreK(\alpha, K)$ for
some $\alpha \in \Segre$ and that the classes $[\SegreK(\alpha, K)]$ are stable under $S_a^+$ and $S_a^-$, meaning
that for an $M \in \SegreK_Q(\alpha,K)$ such that $S_a^+ M \in \SegreK_{\sigma_a Q}(\beta, K)$ we have
that $S_a^+$ induces a bijection from $[\SegreK_Q(\alpha, K)]$ to $[\SegreK_{\sigma_a Q}(\beta, K)]$.

We define the generic Hall algebra $\mathcal{H}_q (Q)$ to be the free $\Q[q]$-module with basis
\[\Set {u_{\alpha} | \alpha \in \Segre}\] and multiplication given by
\[u_{\alpha} \diamond u_{\beta} = \sum_{\gamma} f^{\gamma}_{\alpha \beta} (q) u_{\gamma}. \]
The generic composition algebra $\mathcal{C}_q(Q)$ is then the subalgebra of $\mathcal{H}_q(Q)$ generated
by the simple representations without self-extensions, i.e. the elements of $\Segre$ corresponding
to those. If the quiver is fixed, then we often write $\mathcal{H}_q$ and
$\mathcal{C}_q$ instead of $\mathcal{H}_q(Q)$
and $\mathcal{C}_q(Q)$. Note that the definition of the generic Hall algebra seems to be non-standard.
Again, $\mathcal{H}_q(Q)$ and $\mathcal{C}_q(Q)$ are graded by dimension vector.

We can then specialise $\mathcal{C}_q(Q)$ to any $n \in \Q$ by evaluating the structure constants
given by the polynomials $f^\gamma_{\alpha \beta}$ at $n$. We call this algebra
$\mathcal{C}_n(Q)$. By definition, we have that
\[ \mathcal{C}_{q'} (Q) \cong \mathcal{C}_{\F_{q'}}(Q)\]
for any prime power $q' \in \N$.
In the following, we identify these two algebras.

For a word $w$ in vertices of the quiver $Q$ we define $u_w \in \mathcal{C}_q(Q)$ in the
obvious way. There are polynomials $f^\alpha_w$ for each class
$\alpha \in \Segre$ such that
\[
u_w = \sum_{\alpha} f^\alpha_w u_\alpha.
\]
Note that we have for each finite field $K$ and each $X \in \SegreK(\alpha, K)$ that
$F^X_w = f^\alpha_w(|K|)$.

This thesis mainly deals with the algebras
$\mathcal{C}_0(Q)$ for $Q$ Dynkin or extended Dynkin.

Let $\Q(v)$ be the function field in one variable. Consider it as a $\Q[q]$-algebra
via $v^2 = q$. Denote the twisted composition algebra by
\[ \widetilde{\mathcal{C}}_q(Q) := \mathcal{C}_q(Q) \tensor_{\Q[q]} \Q(v) \]
with the multiplication given by
\[ u_\alpha * u_\beta := v^{\bform{\alpha}{\beta}} u_\alpha \diamond u_\beta. \]
Let $\mathfrak{g}$ be the Kac-Moody Lie algebra associated to the Cartan datum given by
$Q$ (or by $(.,.)_Q$).
C. M. Ringel \cite{Ringel_greenstheom}, J. Green \cite{Green_hallalgs} and G. Lusztig \cite{Lusztig_quantum} showed that
$\widetilde{\mathcal{C}}_q(Q) \cong U^+_q(\mathfrak{g})$
where $U^+_q(\mathfrak{g})$ is the positive part of the quantised enveloping algebra of $\mathfrak{g}$.

G. Lusztig \cite{Lusztig_quantum} showed that
the $\Q(v)$-dimension of the $\dvec{d}$-th graded part $U^+_q(\mathfrak{g})_\dvec{d}$ is
equal to the $\Q$-dimension of the $\dvec{d}$-th graded part $U^+(Q)_\dvec{d}$ of the positive part
of the universal enveloping algebra of $\mathfrak{g}$. Therefore, $\dim \widetilde{\mathcal{C}}_q(Q)_\dvec{d}$
is equal to $\dim U^+(\mathfrak{g})_\dvec{d}$.
The twist does not change the dimension of the graded parts and $\mathcal{C}_q(Q)$ is
free as a $\Q[q]$-module, as a submodule of a free module. Therefore, we obtain that the $\Q[q]$-rank of
the $\dvec{d}$-th graded part
$\mathcal{C}_q(Q)_\dvec{d}$ of the generic composition algebra is equal to $\dim U^+(\mathfrak{g})_\dvec{d}$.
Finally, specialisation does not change the rank of a free module, therefore
$\dim \mathcal{C}_0(Q)_\dvec{d} =\dim U^+(\mathfrak{g})_\dvec{d}$.

Note that $\dim U^+(\mathfrak{g})_\dvec{d}$ is the number of ways of
writing $\dvec{d}$ as a sum of positive roots with multiplicities, in the sense
that we count each positive root $\dvec{f}$ with multiplicity
$\dim \mathfrak{g}_\dvec{f}$ of the corresponding root space.

\section{Generic Extension Monoid}
In \cite{Reineke_monoid}, M. Reineke introduced the generic extension monoid and the composition monoid. We recall briefly how
this is done. Fix
an algebraically closed field $K$.

For two arbitrary sets $U \subseteq \Rep{\dvec{d}}, V \subseteq \Rep{\dvec{e}}$ we define
\begin{align*}
    \mathcal{E}(U,V) := \{ M \in \Rep{\dvec{d} + \dvec{e}} \: |& \: \exists \:A \in U, B \in V \
    \text{and a short exact sequence } \\
     & 0 \rightarrow B \rightarrow M \rightarrow A \rightarrow 0 \}.
\end{align*}

The multiplication on closed irreducible $\GL_\dvec{d}$-stable respectively $\GL_\dvec{e}$-stable
subvarieties $\mathcal{A} \subseteq \Rep{\dvec{d}},
\mathcal{B} \subseteq \Rep{\dvec{e}}$ is defined as:
\[ \mathcal{A} * \mathcal{B} := \mathcal{E}(\mathcal{A}, \mathcal{B}).\]
M. Reineke \cite{Reineke_monoid} showed
that then $\mathcal{A} * \mathcal{B}$ is again closed, irreducible and $\GL_{\dvec{d}+\dvec{e}}$-stable.
Moreover, he showed that $*$ is associative and has a unit: $\Rep{\dvec{0}}$.
The set
\[
\mathcal{M}(Q):= \coprod_{\dvec{d}} \Set{ \mathcal{A} \subset \Rep[Q]{\dvec{d}} | \mathcal{A}
\text{ is a closed, irreducible and } \GL_\dvec{d}\text{-stable subvariety}} 
\]
with this multiplication is therefore a monoid, the
\index{generic extension monoid}\emph{generic extension monoid}
$\mathcal{M}(Q)$.
The \index{composition monoid}\emph{composition monoid}
$\mathcal{CM}(Q)$ is the submonoid
generated by the orbits of simple representations without self-extensions. Note that $\mathcal{M}(Q)$ and
$\mathcal{CM}(Q)$ are graded by dimension vector.

For any word $w = (N_1, \dots, N_r)$ in semisimples we define
$\mathcal{A}_w := \Orbit_{N_1} * \dots * \Orbit_{N_r}$. This is an element of $\mathcal{M}(Q)$. If $N_1, \dots, N_r$
are simple, then $\mathcal{A}_w \in \mathcal{CM}(Q)$. Since there is a simple for each vertex of a quiver, we can
similarly define $\mathcal{A}_w$ for $w$ a word in the vertices of $Q$. Note that the definition of
$\mathcal{A}_w$ makes sense even
if $K$ is not algebraically closed, at least as a set. Therefore, we obtain a set $\mathcal{A}_w$
consisting of all $K$-representations of $Q$ having a filtration of type $w$.

If $\mathcal{A}$ is a closed, irreducible, $\GL_\dvec{d}$-stable subvariety of some $\Rep{\dvec{d}}$, define
$[\mathcal{A}]:= \mathcal{A}/\GL_{\dvec{d}}$ as the set of orbits.
Hence elements of
$[\mathcal{A}]$ correspond to isomorphism classes in $\mathcal{A}$.

More generally, if $\mathcal{A}$ is any subset of $\RepK{\dvec{d}}{K}$ for any field $K$ denote
by $[\mathcal{A}]$ the set of isomorphism classes in $\mathcal{A}$.

If $\mathcal{M}$ is any monoid and $R$ any ring we denote by $R\mathcal{M}$ the monoid algebra
of $\mathcal{M}$ given formal $R$-linear combinations of elements of $\mathcal{M}$ and
the obvious multiplication induced by the multiplications in $\mathcal{M}$ and $R$.

\section{Quiver Flags and Quiver Grassmannians}
Let $M$ be a finite dimensional $K$-vector space and let
\[ \seqv{d} = (d^0, d^1, \dots, d^\nu) \]
be a $(\nu+1)$-tuple of integers. Then one usually defines
$\Fl[K]{\seqv{d}}{M}$ to be the set of sequences of subspaces
$0=U^0  \subset U^1 \subset \dots \subset U^\nu = M$ such that
$\dim_K (U^i) = d^i$. If this set is non-empty, we have that $d^0 = 0$, $d^i \le d^{i+1}$ and
$d^\nu = \dim M$. We call such a sequence a \emph{filtration} of $M$ or of dimension $\dim M$. If we have the
first two conditions we call $\seqv{d}$ a filtration.

Let $\Lambda$ be a $K$-algebra and $M$ a finite dimensional $\Lambda$-module.
Then it is natural to define
\[\Fl[\Lambda]{\seqv{d}}{M} := \Set{ U^0  \subset U^1 \subset \dots \subset U^\nu \in \Fl[K]{\seqv{d}}{M} |
U^i \text{ is a } \Lambda \text{-submodule of } M}.\]

Let $Q$ be a quiver, $K$ a field and $\Lambda := KQ$.
Each module $U$ decomposes, as a $K$-vector space, into the
direct sum of $U \epsilon_i$, $i \in Q_0$.
The flag of submodules $\Fl[\Lambda]{\seqv{d}}{M}$ is therefore the disjoint union into open subvarieties
\[\Fl[Q]{\seqv{\dvec{d}}}{M} := \Set{ U^0  \subset U^1 \subset \dots \subset U^\nu \in \Fl[\Lambda]{\seqv{d}}{M} |
\dvec{\dim}(U^i) = \dvec{d}^i}\]
such that $\sum_{j \in Q_0} d^i_j = d^i$.
We call a sequence of dimension vectors $\flvec{d}=(\dvec{d}^0, \dots, \dvec{d}^\nu)$
a filtration of a $K$-representation $M$, or of a dimension vector $\dimve M$, if
$\dvec{d}^0 = 0$, $\dvec{d}^i \le \dvec{d}^{i+1}$ and $\dvec{d}^\nu = \dimve M$. If
we only have the first two conditions we call $\flvec{d}$ a filtration.

If $\dvec{d}$ is a dimension vector and $M$ a $K$-representation of dimension
vector $\dvec{d}+\dvec{e}$, then we define the quiver Grassmannian as
\[ \Gr[Q]{\dvec{d}}{M}:=\Fl[Q]{(0,\dvec{d}, \dvec{d}+\dvec{e})}{M}. \]

For an algebraically closed
field $K$ and a filtration $\flvec{d}$
we denote by $\mathcal{A}_\flvec{d}(K)$ the closed, irreducible and $\GL_{\dvec{d}^\nu}$-stable
subvariety of $\RepK{\dvec{d}^\nu}{K}$ consisting of all $K$-representations $M$ such that
$\Fl[Q]{\seqv{\dvec{d}}}{M}$ is non-empty. For the geometric statements see chapter \ref{chap:geometry}.

\section{Main Theorems}
We want to investigate the relationship between the composition algebra at $q=0$ and the
composition monoid. For $Q$ a Dynkin or extended Dynkin quiver we obtain a complete answer.
\begin{mtheom}
  Let $Q$ be a Dynkin quiver or an oriented cycle. Then the map
  \begin{alignat*}{2}
    \Psi &\colon\quad & \Q\mathcal{M}(Q) & \rightarrow \mathcal{H}_0(Q)\\
    && \mathcal{A} & \mapsto \sum\limits_{M \in [\mathcal{A}]} u_M \\
  \end{alignat*}
  is an isomorphism of graded $\Q$-algebras.
\end{mtheom}
\begin{proof}
    This is the combination of theorems \ref{cyclic:theom:iso} and \ref{reflflags:theom:dynkin}.
\end{proof}
\begin{mtheom}
  Let $Q$ be an acyclic, extended Dynkin quiver. Then the map
  \[ \Phi \colon\quad \mathcal{C}_0(Q) \rightarrow \Q\mathcal{CM}(Q) \]
  sending $u_{S_i}$ to $\Orbit_{S_i}$ is a graded $\Q$-algebra
  homomorphism with kernel generated by the relations
  \[ (u_{\dvec{\delta}})^r = u_{r\dvec{\delta}} \quad \forall \; r \in \N.\]
\end{mtheom}
\begin{proof}
    This is corollary \ref{extdyn:cory:morph}.
\end{proof}
Moreover, we obtain a basis of PBW-type for the generic composition algebra.

On a more geometric side, we can prove that the quiver flag variety is
irreducible under certain circumstances.
\begin{mtheom}
    Let $K$ be an algebraically closed field.
   Assume that there is an $M \in \mathcal{A}_\flvec{d}(K)$ such that $\dim \Ext^1_{Q}(M,M) =
   \codim \mathcal{A}_\flvec{d}(K)$. Then $\Fl[Q]{\flvec{d}}{M}$ is smooth and
   irreducible.
\end{mtheom}
\begin{proof}
    This is theorem \ref{theom:geomflags:irred}.
\end{proof}
\begin{remark}
    We prove this statement for arbitrary fields $K$.
\end{remark}

\chapter{Cyclic Quiver Case}
\label{chap:cyclic}
\dictum[Archimedes]{Don't disturb my circles.}
\vspace{1cm}
Consider the cyclic quiver $Q=C_n$ of type $\widetilde{A}_n$ (all arrows in one direction):
\[
\makeatletter%
\let\ASYencoding\f@encoding%
\let\ASYfamily\f@family%
\let\ASYseries\f@series%
\let\ASYshape\f@shape%
\makeatother%
\setlength{\unitlength}{1pt}
\includegraphics{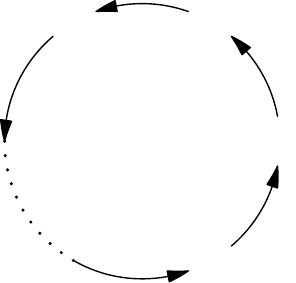}%
\definecolor{ASYcolor}{gray}{0.000000}\color{ASYcolor}
\fontsize{12.000000}{14.400000}\selectfont
\usefont{\ASYencoding}{\ASYfamily}{\ASYseries}{\ASYshape}%
\ASYalign(-4.443095,40.829370)(-0.500000,-0.500000){$0$}
\ASYalign(-24.333386,75.280364)(-0.500000,-0.500000){$1$}
\ASYalign(-64.113968,75.280364)(-0.500000,-0.500000){$2$}
\ASYalign(-24.333386,6.378375)(-0.500000,-0.500000){$n$}
\]
Let $K$ be a field and $\Lambda:=K C_n$.
We want to examine $\repK{C_n}{K}$,
the category of nilpotent $K$-representations over $C_n$. Note that
$\repK{C_n}{K}$ is equivalent to $\Mod \Lambda$, the category of nilpotent $\Lambda$-modules.
Let $S_0, S_1, \dots, S_{n}$ denote
the simple representations in $\repK{C_n}{K}$ corresponding to the vertices of $Q$. We have that
$\Ext(S_i, S_{i+1})\neq 0$ (now and in the
remainder of this section always count modulo $n+1$). For generalities on the cyclic quiver see
\cite{dengdu_monbasesforsln}.

Up to isomorphism there is exactly
one indecomposable $K$-representation $S_i[l]$ of length $l$ with socle
$S_i$. For a partition $\lambda=(\lambda_1 \ge \lambda_2 \ge \dots \ge \lambda_r)$ we set
\[
S_i[\lambda]:= \bigoplus_{k = 1}^r S_i[\lambda_k].
\]
The set of isomorphism classes of representations in $\repK{C_n}{K}$ is therefore in bijection with
\[
\Pi := \Set{ (\pi^{(0)},\dots,\pi^{(n)}) | \pi^{(i)}
\text{ is a partition } \forall i},
\]
where each partition $\pi^{(i)}$
describes the indecomposable summands with socle $S_i$. Hall polynomials exist
with respect to $\Pi$, as shown in \cite{Guo_hallpolyscyclic}.

Now let $M$ be an arbitrary representation of isomorphism class $\pi \in \Pi$. Then
we denote by $u_M$ or $u_\pi$ its symbol in the generic Hall algebra $\mathcal{H}_q(C_n)$.
Most times we will write a partition $\lambda$ of $m'$ in exponential form, i.e.
$\lambda = (1^{s_1} 2^{s_2} \dotsm m^{s_m})$ such that $m'= \sum_{i=1}^{m} s_i i$.

When calculating Hall polynomials, certain quantum numbers appear. Let $R$ be some commutative ring
and let $q \in R$. Usually $R$ will be $\Q[q]$, the polynomial ring in one variable. We define for $r,n \in \N$,
$0 \le r \le n$:
\begin{align*}
    [n]_q &:=  1 + q + \dots + q^{n-1} \\
    [n]_q ! &:= \prod_{i=1}^{n} [i]_q \\
    \qbinom{n}{r}_q &:= \frac{[n]_q !}{[r]_q ! [n-r]_q!}.
\end{align*}
Obviously, $\qbinom{n}{r}_0 = 1$.

Now we want to calculate the Hall polynomial between direct sums of a simple and two arbitrary representations.
We use the following well-known lemma which describes the structure of the Grassmannian of subspaces.
\begin{lmm}
    Let $V, W$ be two $K$-vector spaces and denote by $\pi_W \colon V \oplus W \rightarrow W$
    the second projection.
    Then the map
    \begin{alignat*}{2}
        \phi_{V,W}& \colon\quad&
\Gr{r}{V \oplus W} & \rightarrow  \coprod_{s+t=r} \Gr{s}{V} \times \Gr{t}{W} \\
        \\
               &&  U          & \mapsto   (U \cap V, \pi_W (U))
    \end{alignat*}
    is surjective and the fibre over $(A,B)$ is isomorphic to $\Hom_K (B, W/A)$.
\end{lmm}
\begin{proof}
    See for example \cite{CalderoChapoton_clusterhall} for a generalisation to quiver representations.
\end{proof}
Let $\lambda = (1^{s_1} 2^{s_2} \dotsm t^{s_t})$ be a partition and let $S$ be a simple representation.
Then $\soc S[\lambda] \cong S^{\sum s_i}$ and
every inclusion of $S$ into $S[\lambda]$ corresponds to a one dimensional subspace of
$\soc S[\lambda]$. More generally, every inclusion of $S^k$, $k\in \N$, into $S[\lambda]$ corresponds to a
$k$ dimensional subspace of $\soc S[\lambda]$.
\begin{lmm}
    Let $\lambda = (1^{s_1} 2^{s_2} \dotsm m^{s_m})$ be a partition, $S$ a simple representation,
    $M := S[m]^{s_m}$ and $N := \bigoplus_{i=1}^{m-1} S[i]^{s_i}$. So
    $X:=S[\lambda]= M \oplus N$. Fix $r \in \N$. Let $U$ be a subspace of $\soc X$ of dimension $r$.
    Then
    $X/U \cong M/A \oplus N/B$ where $(A,B) = \phi_{\soc M,\soc N}(U)$ and the fibre of
    $\phi_{\soc M, \soc N}$ over $(A,B)$ has dimension $(s_m - \dim A) \dim B$.
    \label{quotconstoncell}
\end{lmm}
\begin{proof}
    For any $i \le j$ there is an
    inclusion $S[i] \hookrightarrow S[j]$ which restricts to an isomorphism of the one
    dimensional socles.
    Therefore, each vector space homomorphism $\soc N \rightarrow \soc M$ extends to a 
    homomorphism $N \rightarrow M$ of $K$-representations. Note that this extension is not unique.

    We have $U \in \Gr{r}{\soc M \oplus \soc N}$. Let $(A,B) :=  \phi_{\soc M,\soc N}(U)$.
    We can choose complements $A'$ of $A$ in $\soc M$ and $B'$ of $B$ in $\soc N$. Then, by the
    preceding lemma, $U$ is given by a homomorphism $g \in \Hom_K(B, A')$ and we denote its extension by $0$ on
	$A$ and $B'$ by
	$g'\colon \soc N \to \soc M$. Let $\tilde g \colon N \to M$ be any
    homomorphism of $K$-representations induced by $g'$. Now the
    following diagram commutes.
    \[
    \xymatrix{
    A \oplus B \: \ar@{^(->}[rr]^(0.4){\left( \begin{smallmatrix}
        \id_A & 0\\
          0   & g\\
          0   & \id_B \\
          0   &  0
      \end{smallmatrix}\right)}
      \ar@{=}[dd]
      & & A \oplus A' \oplus B \oplus B' \: \ar@{^(->}[rr] & &
      M \oplus N \ar@{<-}[dd]^{\left(
          \begin{smallmatrix}
              \id_M & \tilde g\\
              0 & \id_N
          \end{smallmatrix}
          \right)}_{\begin{sideways}$\sim$ \end{sideways}}\\
    \\
    A \oplus B \: \ar@{^(->}[rr]^(0.4){\left( \begin{smallmatrix}
        \id_A & 0\\
          0   & 0\\
          0   & \id_B \\
          0   &  0
      \end{smallmatrix}\right)} & & A \oplus A' \oplus B \oplus B' \: \ar@{^(->}[rr] \ar[uu]_{\left( \begin{smallmatrix}
          \id_A & 0 & 0 & 0\\
          0   & \id_{A'} & g & 0\\
          0   &    & \id_{B} & 0\\
          0 & 0 & 0 & \id_{B'}
      \end{smallmatrix}\right)}^{\begin{sideways}$\sim$ \end{sideways}} & & M \oplus N
    }
    \]
    Hence the
    cokernel of the map in the top row is the same as the cokernel of the map in the bottom row,
    therefore $(M \oplus N) /U \cong M/A \oplus N/B$.
\end{proof}
\begin{propn}
    Let $s=(s_1, \dots, s_m)$ and $t=(t_1, \dots, t_m) \in \N^m$ such that $t_i \le s_i$.
    Let $S$ and $T$ be simples such that
    $\Ext(T,S) \neq 0$. Set
    \begin{align*}
        \lambda &:= (1^{s_1}2^{s_2}\dotsm m^{s_m})&\mu&:=(1^{s_1 - t_1} 2^{s_2-t_2} \dotsm m^{s_m - t_m})\\
        \nu&:=(1^{t_2}2^{t_3}\dotsm (m-1)^{t_m})& r&:= \sum t_i\\
        X&:=S[\lambda] & Y(s,t)&:=S[\mu]\oplus T[\nu].
    \end{align*}
    Then
    \[f^X_{Y(s,t) \; S^r}(q) = \prod_i \qbinom{s_i}{t_i}_q q^{\sum_{j<i} t_j (s_i - t_i)}. \]
    Moreover, all quotients of $X$ by $S^r$ are of the form of $Y(s,t)$ for some choice of $t$.
    \label{hallalgs_ssarb}
\end{propn}
\begin{proof}
    We proof this by induction on $m$. For $m=0$ the claims are true. Now let $m > 0$.
    Then $X \cong M \oplus N$, where $M := S[m]^{s_m}$ and $N := \bigoplus_{i=1}^{m-1} S[i]^{s_i}$,
    hence, without loss of generality, we can assume $X = M \oplus N$.
    Now, a subrepresentation of $X$ isomorphic to $S^r$ is given by an element
    $U \in \Gr{r}{\soc M \oplus \soc N}$. Let $(A,B) := \phi_{\soc M, \soc N}(U)$ as in the
    previous lemma. Then $X /U \cong M/A \oplus N/B$. This yields by induction that
    all quotients are of the desired form. One has that
    $M/A \cong S[m]^{s_m-\dim A} \oplus T[m-1]^{\dim A}$ and no
    quotient of $N$ by $B$ has a summand isomorphic to $S[m]$. Therefore, if $X/U \cong Y(s,t)$,
    $\dim A$ has to be equal to $t_m$.

    Set $s':= (s_1, s_2, \dots, s_{m-1})$, $t':= (t_1, t_2, \dots, t_{m-1})$ and analogous
    as before
    $\lambda'$, $\mu'$, $\nu'$, $X'$, $Y(s',t')$ and $r'$.
    The number of subspaces $U$ such that $\phi_{\soc M, \soc N}(U) = (A,B) $ is given
    by $q^{(s_m - \dim A) \dim B } = q^{(s_m - t_m) \sum_{j<m} t_j}$. Therefore, we have:
    \[ f^X_{Y(s,t) \; S^r}(q) = q^{(s_m - t_m) \sum_{j<m} t_j}
    \qbinom{s_m}{t_m}_q f^{X'}_{Y(s',t') \; S^{r'}}(q).\]
    By induction this is equal to
    \[ = q^{(s_m - t_m) \sum_{j<m} t_j}
    \qbinom{s_m}{t_m}_q \prod_{i=1}^{m-1} \qbinom{s_i}{t_i}_q q^{\sum_{j<i} t_j (s_i - t_i)} =
    \prod_{i=1}^{m} \qbinom{s_i}{t_i}_q q^{\sum_{j<i} t_j (s_i - t_i)}. \]
\end{proof}
\begin{defn}
Let $S$ and $T$ be simples such that $\Ext(T,S) \neq 0$.
Let $X = S[\lambda]$ for a partition $\lambda=(1^{s_1}2^{s_2}\dotsm m^{s_m})$ and $k \in \N$
such that $k \le \sum_i s_i = l(\lambda)$.
Define $t=(t_1, \dots, t_m)$ recursively by
\begin{align*}
    t_m &:= \min\left\{s_m, k\right\},\\
    t_i &:= \min\left\{s_i, k - \sum_{j>i} t_j\right\} \quad \text{for all } 0 \le i \le m-1.
\end{align*}
We define
\[ Q(X, S^k) := Y(s,t) = \bigoplus_{i=1}^{m} S[i]^{s_i - t_i} \oplus \bigoplus_{i=1}^{m-1} T[i]^{t_{m+1}}.\]

More generally, if $X = \bigoplus_{i=0}^{n} S_i[\pi^{(i)}]$ for a $\pi \in \Pi$ and $N = \bigoplus S_i^{k_i}$ for
some $k_i \in \N$ such that $k_i \le l(\pi^{(i)})$, define
\[ Q(X,N) := \bigoplus_{i=0}^{n} Q(S_i[\pi^{(i)}], S_i^{k_i}). \]
\end{defn}
We will prove in corollary \ref{degen-minimal-elts} that $Q(X,N)$ is maximal with respect to the degeneration order,
hence the
quotient of $X$ by $N$ with the smallest orbit dimension.

We obtain the following.
\begin{cory}
    Let $X = S[\lambda]$ for some $\lambda=(1^{s_1}2^{s_2}\dotsm m^{s_m})$. Let $M$ be a quotient of $X$ by
    $S^k$, $k \in \N$ with $k \le l(\lambda)$. Then
    \[ f^X_{M S^k}(0) =
    \begin{cases}
        1 & \text{if } M \cong Q(X,S^k)\\
        0 & \text{otherwise}.
    \end{cases}\]
    \label{cyclic-one-semisimple-factor}
\end{cory}
\begin{proof}
	By proposition \ref{hallalgs_ssarb}, $M = Y(s,t)$ for some $t=(t_1, \dots, t_m) \in \N^m$ with $t_i \le s_i$
    and $\sum t_i = k$.
    We have that
    \[
    \qbinom{s_i}{t_i}_0 = 1.
    \]
    Therefore,  $f^X_{M S^k}(0) \neq 0$ if and only if
    $\sum_{j<i} t_j (s_i - t_i) =0$ for all $1 \le i \le m$. This is the case if either
    $s_i = t_i$ or that for all $j<i$ we have that $t_j = 0$.
    This is exactly the way we chose $Q(X, S^k)$ in the definition
    and in this case we have
    $f^X_{M S^k}(0) = 1$.
\end{proof}
Now we are able to describe the coefficients modulo $q$ for an extension with a semisimple.
\begin{lmm}
    Let $N=\bigoplus_{i=0}^{n} S_i^{k_i}$, $k_i \in \N$, be a semisimple representation. Let
    $X=\bigoplus_{i=0}^{n} S_i[\pi^{(i)}]$, $\pi \in \Pi$, be arbitrary and let
    $M\in \repK{C_n}{K}$ be a quotient of $X$ by $N$.
    Then
    \[ f_{M N} ^X (0) = \begin{cases}1 & \text{if }
        M \cong Q(X,N) \\
         0 & \text{otherwise}.
     \end{cases}
    \]
    \label{hallnumber-semisimple-arbitrary}
\end{lmm}
\begin{proof}
    Since $\Hom(S_i^{k_i}, S_j[\pi^{(j)}])=0$ for $i \neq j$ every short exact sequence
    $0 \rightarrow N \rightarrow X \rightarrow M \rightarrow 0$ is the direct sum of
    short exact sequences of the form
    \[
    0 \rightarrow S_i^{k_i} \rightarrow X_i \rightarrow M_i \rightarrow 0,
    \]
    where $X_i = S_i[\pi^{(i)}]$ and for some representations $M_i$ such that $\bigoplus_{i=0}^{n} M_i \cong M$. So we have
    \[f_{M N}^X = \sum_{\substack{(M_0,\dots,M_{n}) :\\ \bigoplus M_i \cong M}}
    \prod_{i=0}^{n} f_{M_i S_i^{k_i}}^{X_i} \]
    where $X_i := S_i[\pi^{(i)}]$. But now $f_{M_i S_i^{k_i}}^{X_i}(0)$ is non-zero if and only if
    \[
    M_i \cong Q(X_i,S_i^{k_i})\]
    by lemma \ref{hallnumber-semisimple-arbitrary}. Moreover, the same lemma yields
    $f_{Q(X_i,S_i^{k_i}) \; S_i^{k_i}}^{X_i} (0)= 1$. Hence we are done.
\end{proof}
\begin{lmm}
    Let $w = (N_1, N_2, \dots, N_r)$ be a word in semisimples. Let $M \in \repK{Q}{K}$ and
    $M' \in \repK{Q}{L}$ for two arbitrary fields $K$ and $L$ such that $M$ and $M'$ are
    of isomorphism type $\pi \in \Pi$. Then $M$ has a filtration of type $w$ if and only
    $M'$ has. In other words, the sets $[\mathcal{A}_w]$ can be considered as subsets of $\Pi$
    and do not depend on the field we are working over.
    \label{cyclic:lmm:filtrindept}
\end{lmm}
\begin{proof}
    We prove the claim by induction on $r$. If $r=1$, then the claim is trivial.
    Now let $r>1$ and set $w':=(N_1, N_2, \dots, N_{r-1})$. A representation $M$ has a
    filtration of type $w$ if and only if it has a subrepresentation $U$ isomorphic to $N_r$
    such that the quotient $M/U \in \mathcal{A}_{w'}$. Since the polynomials
    of proposition \ref{hallalgs_ssarb} have positive coefficients and by
    a similar argumentation as in lemma \ref{hallnumber-semisimple-arbitrary} we have that
    the possible isomorphism classes of the quotients only depend on the isomorphism classes of
    $M$ and $N_r$. By induction we have that the isomorphism classes in $\mathcal{A}_{w'}$ do not
    depend on the field. This finishes the proof.
\end{proof}
Now we show that $Q(X,N)$ is maximal with respect to the degeneration order among the quotients of $X$ by $N$.
\begin{lmm}
    Let $N$, $X_1, X_2$ be arbitrary representations, $f_1 \colon N \rightarrow X_1$
    and $f_2 \colon N \rightarrow X_2$ two injections. Moreover, let $g \colon X_1 \rightarrow X_2$
    be a morphism such that $g  f_1 = f_2$.

    Then we have a short exact sequence
    \[ \xymatrix{ 0 \ar[r] &
    X_1 \ar[r] &
    X_2 \oplus X_1/f_1(N) \ar[r]
    & X_2/f_2(N) \ar[r]& 0} \]
    and therefore
    $X_2  \oplus (X_1 /f_1(N)) \ledeg X_1  \oplus (X_2 /f_2(N))$.
    \label{inclusion_hence_degenerates}
\end{lmm}
\begin{proof}
    We construct an extension degeneration. Let
    $\pi_1 \colon X_1 \rightarrow X_1 /f_1(N)$
    and $\pi_2 \colon X_2 \rightarrow X_2/f_2(N)$ be the canonical projections and let
    $\overline{g} \colon X_1/f_1(N) \rightarrow X_2/f_2(N)$ be the map induced by $\pi_2 g$, which exists since
    $f_1(N) \subseteq \ker (\pi_2 g)$.
    By construction we have the following commutative diagram:
    \[\begin{CD}
        0 @>>> N @>{f_1}>> X_1 @>{\pi_1}>> X_1/f_1(N) @>>> 0\\
        @. @| @VV{g}V @VV{\bar g}V\\
        0 @>>> N @>{f_2}>> X_2 @>{\pi_2}>>  X_2/f_2(N) @>>> 0.
    \end{CD}\]
    This is a pullback, therefore
    \[ \xymatrix{ 0 \ar[r] &
    X_1 \ar[r]^(0.3){\left(\begin{smallmatrix}
        g\\
        \pi_1
    \end{smallmatrix}\right)} &
    X_2 \oplus X_1/f_1(N) \ar[r]^(0.6){\left(\begin{smallmatrix}
        \pi_2 & -\overline{g}
    \end{smallmatrix}\right)}
    & X_2/f_2(N) \ar[r]& 0} \]
    is the desired short exact sequence.
\end{proof}
\begin{cory}
    Let $M$, $N$ be arbitrary representations and let $h \colon N \rightarrow M$
    be an injection. Let $g \in \End(M)$ be any endomorphism such that $g h(N) \cong N$.
    Then
    \[ M/h(N) \ledeg M/gh(N). \]
    \label{deg_one_zero}
\end{cory}
\begin{proof}
From lemma \ref{inclusion_hence_degenerates} we have an extension of the form
\[
    0 \rightarrow M \rightarrow M \oplus M/h(N)
    \rightarrow M/gh(N) \rightarrow 0.
    \]
By a result of C. Riedtmann \cite[prop. 4.3]{Riedtmann_degens} this yields that
\[M/h(N) \ledeg M/gh(N).\]
\end{proof}

\begin{cory}
    Let $X, M \in \repK{C_n}{K}$ and $N$ a semisimple representation such that there is a short exact sequence
    $0 \rightarrow N \rightarrow X \rightarrow M \rightarrow 0$. Then $M$ degenerates to
    $Q(X,N)$.
    \label{degen-minimal-elts}
\end{cory}
\begin{proof}
    As before, the short exact sequence
	is a direct sum of short exact sequences of the form 
	\[
	0 \rightarrow S_i^{k_i} \rightarrow X_i \rightarrow M_i \rightarrow 0,
	\]
	where $X \cong \bigoplus X_i$ and $M \cong \bigoplus M_i$ for some
    $X_i$ with socle in $\add S_i$.
	
	It is obviously enough to show the claim for every short exact
	sequence of this form. We can
    therefore assume that $X=S[\lambda]$, $M=S[\mu] \oplus T[\nu]$ and $N=S^k$ for
	$k \in \N$, $S$ and $T$ simples such that $\Ext(T,S) \neq 0$ and
	partitions $\lambda=(\lambda_1,\cdots,\lambda_m)$, $\mu$ and $\nu$.
    Let $U$ be the $k$-dimensional subspace of $\soc X=\bigoplus \soc S[\lambda_i]$ induced
	by the inclusion of $S^k$ into $X$. We have that $X/U \cong M$.
    We want to apply the previous
    corollary to $X$, so we have to find an endomorphism $f$ which maps $U$ to $U'$, where $U'$ is a subspace
    of $\soc X$ such that $X/U' \cong Q(X,S^k)$.

    By repeated application of the construction used
    in the proof of lemma \ref{quotconstoncell} we obtain an automorphism $\phi$ of
	$X$ such that 
	\[
	\phi(U) = \bigoplus_{j=1}^k \soc S[\lambda_{i_j}] \subset X,
	\]
	for a sequence of integers $1 \le i_1 < i_2 < \dots < i_k \le m$,
	since $\soc S[\lambda_{i_j}] \cong k$. Note that $X/\Phi(U) \cong M$, we can therefore
	assume that $U= \phi(U)$.
	Let 
	\[
	U':= \bigoplus_{j=1}^k \soc S[\lambda_j] \subset X.
	\]
	By definition we have that $X/U' = Q(X, S^k)$.

	For $i \le j$ there is an injection $\varphi_{i,j}$ from $S[i]$ to $S[j]$ which restricts
    to an isomorphism on the one dimensional socles. Therefore, the map
	\begin{alignat*}{2}
	  g &\colon\quad& \bigoplus_{j=1}^k S[\lambda_{i_j}] &\rightarrow
	  \bigoplus_{j=1}^k S[\lambda_j]\\
	  && (x_1, \dots, x_k) &\mapsto  (\varphi_{\lambda_{i_1}, \lambda_1}(x_1), \dots,
	  \varphi_{\lambda_{i_k},\lambda_k}(x_k))\\
	\end{alignat*}
	induces an isomorphism from $U$ to $U'$ and can be extended by $0$ on the remaining
	summands to an endomorphism of $X$.
	By applying the previous lemma we obtain that
	\[
	M= X/U \ledeg X/U' = Q(X,S^k),
	\]
	and this yields the desired result.
\end{proof}
Now we are able to describe the monomial elements of $\mathcal{H}_0(C_n)$, the generic Hall
algebra of the cyclic quiver specialised at $q=0$.
\begin{theom}
    Let $w=(N_1,\dots,N_r)$ be a sequence of semisimple representations.
    Then
    \[u_w := u_{[N_1]} \diamond u_{[N_2]} \diamond \dots \diamond u_{[N_r]}
    = \sum_{[M] \in [\mathcal{A}_w]} u_{[M]} \in \mathcal{H}_0(C_n).\]
    \label{coeff-compmonoid-areone}
\end{theom}
\begin{proof}
    By lemma \ref{cyclic:lmm:filtrindept} the expression on the right hand side is well-defined.
    We prove the theorem by induction on $r$. For $r=1$ the statement is trivial.
    Now let $w'=(N_1, \dots, N_{r-1})$. We need to show that
    $\sum_{[M]\in [\mathcal{A}_{w}]} u_{[M]} = u_{w} \in \mathcal{H}_0(C_n)$.
    First, it is obvious that
    every isomorphism class appearing in $u_w = u_{w'} \diamond u_{[N_r]}$
    with non-zero coefficient is in $[\mathcal{A}_w]$.

    Now let $X$ be in $\mathcal{A}_w$.
    We have to show that the coefficient of $u_{[X]}$ in $u_w = u_{w'} \diamond u_{[N_r]} \in \mathcal{H}_0(C_n)$
    is $1$.
    This coefficient is
    \[
    \sum_{[M] \in [\mathcal{A}_{w'}]} f_{M N_r}^X(0).
    \]
    By lemma \ref{hallnumber-semisimple-arbitrary}
    this sum is $1$ if and only if $[Q(X,N_r)]$ is in $[\mathcal{A}_{w'}]$ and $0$ otherwise.
    So it remains to show that $[Q(X,N_r)] \in [\mathcal{A}_{w'}]$. By lemma \ref{cyclic:lmm:filtrindept}
    we can assume that we work over an algebraically closed field $K$ to check this.
    Now, since $X \in \mathcal{A}_{w}$,
    we know that there is at least one $M \in \mathcal{A}_{w'}$ such that there is a short exact
    sequence $0 \rightarrow N_r \rightarrow X \rightarrow M \rightarrow 0$. Moreover, $\mathcal{A}_{w'}$ is closed
    as a variety over $K$.
    Via lemma \ref{degen-minimal-elts}
    we know that $M$ degenerates
    to $Q(X,N_r)$ and therefore $Q(X,N_r) \in \mathcal{A}_{w'}$.
\end{proof}
%
By using this we obtain the following result.
\begin{theom}
    The map
    \begin{alignat*}{2}
    \Psi &\colon\quad& \Q\mathcal{M}(C_n) &\rightarrow \mathcal{H}_0(C_n) \\
    &&\mathcal{A} &\mapsto \sum\limits_{[M] \in [\mathcal{A}]} u_{[M]}
    \end{alignat*}
    is an isomorphism  of graded rings.
    \label{cyclic:theom:iso}
\end{theom}
\begin{proof}
    The sets
    $\Orbit_{N}$ for $N$ semisimple generate $\mathcal{M}(C_n)$ and therefore
    we can apply lemma \ref{cyclic:lmm:filtrindept} to show that the map is well-defined.
    The map $\Psi$ obviously maps $\Orbit_N$ to $u_{[N]}$ for $N$ a semisimple representation.
    Applying theorem \ref{coeff-compmonoid-areone} yields that $\Psi$ is
    a homomorphism of rings. More precisely, if $w$ and $v$ are words in semisimples, we have
    that
    \begin{multline*}
        \Psi(\mathcal{A}_w * \mathcal{A}_v)
        = \Psi(\mathcal{A}_{wv})
        = \sum_{[M] \in [\mathcal{A}_{wv}]} u_{[M]} = u_{wv} \\
        =u_w \diamond u_v =
        \sum_{[M] \in [\mathcal{A}_w]} u_{[M]} \diamond
        \sum_{[M] \in [\mathcal{A}_v]} u_{[M]}
        = \Psi(\mathcal{A}_w) \diamond \Psi(\mathcal{A}_v).
    \end{multline*}
    Moreover, $\Psi$ is obviously a graded homomorphism.

    Now we show that $\Psi$ is an isomorphism by
    showing that it is an isomorphism of $\Q$-vector spaces
    for every graded component. Note that every element $\mathcal{A} \in \mathcal{M}(C_n)$
    is given by the orbit closure of some representation since we are only considering nilpotent representations.
    Therefore, the $\Q$-dimension of the graded components of $\mathcal{M}(C_n)$ and $\mathcal{H}_0(C_n)$ agree,
    since both are equal to the number of isomorphism classes of the given dimension vector.

    The degeneration order is a partial order on the isomorphism classes of representations and we have that
    $\Psi(\overline{\Orbit_{X}}) = u_X + \sum_{X <_{\text{deg}} Y} u_Y$.
    Therefore, on the graded components, $\Psi$ is given by a unipotent matrix and so is an isomorphism.
\end{proof}
\begin{cory}
    The isomorphism $\Psi$ restricts to an isomorphism
    \[\Q \mathcal{CM}(C_n) \cong \mathcal{C}_0(C_n).\]
\end{cory}
\begin{proof}
    Everything follows from the theorem since $\Psi$ maps $\Orbit_{S_i}$ to $u_{[S_i]}$ and these are the generators
    of $\Q \mathcal{CM}(C_n)$ respectively $\mathcal{C}_0(C_n)$.
\end{proof}
\begin{remark}
    We call an isomorphism class $\pi=(\pi^{(0)}, \dots, \pi^{(n)}) \in \Pi$ separated
    if for each $k\ge 1$ there is some $i_k \in \{0,\dots,n\}$ such
    that $\pi^{(i_k)}_j \neq k$ for all $j\ge 1$. In other words, $\pi^{(i_k)}$ has
    no part of size $k$. We denote by $\Pi^s$ the set of all
    separated isomorphism classes. We call a representation separated if its isomorphism class
    is separated. B. Deng and J. Du \cite[Theorem 4.1]{dengdu_monbasesforsln} show that
    \[
    \mathcal{CM}(C_n) = \Set{ \overline{\Orbit_M} | M \text{ is separated}}.
    \]
    Therefore, an element of the generic extension monoid is in the composition monoid if and
    only if it is the orbit closure of a separated representation.
\end{remark}


\chapter{Composition Monoid}
\dictum[Oscar Wilde]{I love hearing my relations abused. It is the only thing that makes me put up with them at all.}
\label{chap:compmon}
\section{Generalities}
Let $Q$ be an arbitrary quiver and $K$ an algebraically closed field.
M. Reineke proved the following.
\begin{theom}[\cite{Reineke_monoid}]
	\label{compmon:codimform}
    Let $\chark K = 0$.
    Let $\mathcal{A} \subseteq \RepK{\dvec{d}}{K}$, $\mathcal{B} \subseteq \RepK{\dvec{e}}{K}$ be closed irreducible
	subvarieties. Then we have:
	\[\codim (\mathcal{A} * \mathcal{B}) \le \codim \mathcal{A} + \codim \mathcal {B} +
	\ext(\mathcal{B}, \mathcal{A}). \]
	If $\ext(\mathcal{A}, \mathcal{B}) = 0$, or $\mathcal{A} = \rep{\dvec{d}}$ and $\mathcal{B} = \rep{\dvec{e}}$,
	equality holds.
\end{theom}
\begin{remark}
  A similar theorem holds for elements of the composition monoid in arbitrary characteristic. We
  prove this in chapter \ref{chap:geometry}.
\end{remark}
From this we obtain the following.
\begin{cory}
	\label{compmon:noext->dsum}
	Let $M$ and $N$ be representations such that $[M,N]^1=0$. Then
	\[ \overline{\Orbit_M} * \overline{\Orbit_N} = \overline{\Orbit_{M \oplus N}}.\]
\end{cory}
We can prove this corollary, without using the theorem, in arbitrary characteristic.
For this, we need the following.
\begin{lmm}
	Let $\mathcal{A}, \mathcal{B}$ be two sets in the generic extension monoid and let
	$U \subseteq \mathcal{A}, V \subseteq \mathcal{B}$ be two open subsets contained in
	those. Then $\mathcal{A} * \mathcal{B} = \overline{\mathcal{E}(U, V)}$.
	\label{compmon:openopenisopen}
\end{lmm}
\begin{proof}
	To prove this we use the setup in \cite{Reineke_monoid}. The subvariety $U \times V$ is open, thus dense
	in $\mathcal{A} \times \mathcal{B}$. Let $\mathcal{Z}(U,V)$ be the set of
	$x \in \rep{\dvec{d}+\dvec{e}}$ such that
	\[x_{\alpha} =\begin{pmatrix}
		u_{\alpha} & \zeta_{\alpha} \\
		0        & v_{\alpha}
	\end{pmatrix}\]
	for every arrow $\alpha \in Q_1$, where $u \in U, v \in V$ and $\zeta$ arbitrary. Now,
	$\mathcal{Z}(U,V)$ is an open subvariety of the irreducible variety
    $\mathcal{Z}(\mathcal{A},\mathcal{B})$, thus dense.
    Therefore, the $\GL_{\dvec{d}+\dvec{e}}$-saturation of $\mathcal{Z}(U,V)$ is dense
    in the $\GL_{\dvec{d}+\dvec{e}}$-saturation of $\mathcal{Z}(\mathcal{A}, \mathcal{B})$.
	But the former is equal to $\mathcal{E}(U,V)$ whereas the latter is equal to
	$\mathcal{A} * \mathcal{B}$, hence the claim follows.
\end{proof}
\begin{proof}[Proof of corollary \ref{compmon:noext->dsum}]
	We know that $\Orbit_M$ is open in $\overline{\Orbit_M}$ and
	$\Orbit_N$ is open in $\overline{\Orbit_N}$. By using lemma \ref{compmon:openopenisopen}
	we have $\overline{\Orbit_M} * \overline{\Orbit_N} = \overline{\mathcal{E}(\Orbit_M, \Orbit_N)}$.
	But $\mathcal{E}(\Orbit_M , \Orbit_N) = \Orbit_{M \oplus N}$ since every 
	extension is split and the corollary follows.
\end{proof}

We have the following fundamental relation in $\mathcal{CM}(Q)$.
\begin{propn}
    Let $\dvec{d}$ and $\dvec{e}$ be two dimension vectors such that $\ext(\dvec{e}, \dvec{d}) = 0$. Then
    \[ \Rep{\dvec{d}} * \Rep{\dvec{e}} = \Rep{\dvec{d} + \dvec{e}}. \]
    Note that $\ext(\dvec{e}, \dvec{d})$ only depends on the dimension vectors and the quiver and not on the field we
    are working over.
    \label{compmon:propn:dimvecmult}
\end{propn}
\begin{proof}
    This is an immediate consequence of theorem
    \ref{intro:theom:extallsubreps}.
\end{proof}
We will deduce from this relation some relations on multiples of Schur roots.
\begin{lmm}
    Let $\dvec{d}$, $\dvec{e}$ be two dimension vectors with $\ext(\dvec{e}, \dvec{d}) = 0$.
    Then
    \[\Rep{\dvec{d}} * \Rep{\dvec{e}}=
    \Rep{r_1\dvec{f}^1} * \dots * \Rep{r_l\dvec{f}^l}\] where $\sum r_i \dvec{f}^i$
    is the
    canonical decomposition of $\dvec{d} + \dvec{e}$.
    \label{compmon:schurrootsmult} Moreover,
    if $\dvec{a}$ is any other dimension vector, then
    \[\ext(\dvec{a}, \dvec{d}) + \ext(\dvec{a}, \dvec{e}) \ge \sum \ext(\dvec{a}, r_i\dvec{f}^i)\] and
    \[\ext(\dvec{d}, \dvec{a}) + \ext(\dvec{e}, \dvec{a}) \ge \sum \ext(r_i\dvec{f}^i, \dvec{a}).\]
\end{lmm}
\begin{proof}
    Using the previous proposition we obtain $\rep{\dvec{d}} * \rep{\dvec{e}} = \rep{\dvec{d} + \dvec{e}}$.

    Let $\dvec{a}$, $\dvec{b}$ and $\dvec{c}$ be three arbitrary dimension vectors.
    We have 
    \[
    \ext(\dvec{a}, \dvec{b}) + \ext(\dvec{a}, \dvec{c}) \ge
    \ext(\dvec{a}, \dvec{b}+\dvec{c}).
    \]
    To see
    this take representations $A \in \rep{\dvec{a}}$, $B \in \rep{\dvec{b}}$ and
    $C \in \rep{\dvec{c}}$ such that $[A,B]^1 = \ext(\dvec{a}, \dvec{b})$ and
    $[A,C]^1= \ext(\dvec{a},\dvec{c})$. This is possible since the set of representations taking minimal $\Ext$
    values is open in $\rep{\dvec{a}}$. Now we have
    \[
        \ext(\dvec{a}, \dvec{b}) + \ext(\dvec{a}, \dvec{c}) = [A,B]^1 + [A,C]^1\\
        = [A, B \oplus C]^1 \ge \ext(\dvec{a}, \dvec{b} + \dvec{c}).
    \]

    Therefore, if
    $\ext(\dvec{a}, \dvec{b}) =0 = \ext(\dvec{a}, \dvec{c})$ for three arbitrary dimension vectors
    $\dvec{a}$, $\dvec{b}$ and $\dvec{c}$, then
    $\ext(\dvec{a}, \dvec{b} + \dvec{c}) = 0$.
    If $\sum r_i \dvec{f}^i$ is the canonical decomposition of $\dvec{d}+\dvec{e}$, then, by definition,
    $\ext(\dvec{f}^i, \dvec{f}^j) = 0$ for all $1 \le i\neq j \le l$.
    Therefore, we can iteratively apply the previous proposition
    to obtain that
    \[
    \Rep{r_1\dvec{f}^1} * \dots * \Rep{r_l\dvec{f}^l}=
    \Rep{\dvec{d}+\dvec{e}}
    \]
    and this proves the first claim.

    Now we prove the second claim.
    We can choose representations $C' \in \rep{\dvec{a}}$ and $A_i \in \rep{r_i\dvec{f}^i}$ such that
    \[
	[C', A_1 \oplus \dots \oplus A_l]^1 = \ext(\dvec{a}, r_1\dvec{f}^1 + \dots + r_l\dvec{f}^l) =
    \ext(\dvec{a}, \dvec{d}+\dvec{e}),
	\]
	since a general representation in $\rep{\dvec{d} + \dvec{e}}$ will have a
    decomposition into summands like this. So we have
    \begin{multline*}
        \ext(\dvec{a}, \dvec{d}) + \ext(\dvec{a}, \dvec{e}) \ge \ext(\dvec{a}, \dvec{d} + \dvec{e})\\
        = \ext(\dvec{a}, r_1\dvec{f}^1 + \dots + r_l\dvec{f}^l) = [C', A_1 \oplus \dots \oplus A_l]^1\\
        = \sum_{i=1}^l [C',A_i]^1 \ge \sum_{i=1}^l \ext(\dvec{a}, r_i\dvec{f}^i).
    \end{multline*}
    The other statement is proved dually.
\end{proof}

\section{Partial Normal Form in Terms of Schur Roots}
For the following let $Q$ be a connected acyclic quiver.
We say that a root $\dvec{d}$ is preprojective/regular/preinjective if one
(and therefore all) indecomposable representation of dimension vector $\dvec{d}$ is. If $Q$ is
Dynkin, then all roots are preprojective and preinjective, but for convenience we set them
to be preprojective and not preinjective. We also choose a total order $\prec_t$ on the preprojective
and preinjective Schur roots refining the order $\prec$ on $\mathcal{P} \cup \mathcal{I}$.

We define a new monoid $\RMon(Q)$, the Schur root monoid. We will use it to obtain a partial normal
form in the composition monoid. For $Q$ Dynkin we will show that it is isomorphic
to $\mathcal{CM}(Q)$. Moreover we will later show that the relations of $\RMon(Q)$ also hold in $\mathcal{C}_0(Q)$
for $Q$ an extended Dynkin quiver.
\begin{defn}
    The monoid $\RMon(Q)$ is given by the generators
    \[\Set{ \SRoot{r\dvec{d}} | \dvec{d} \text{ is a Schur root}, r\in\N}\] and the relations
    \begin{align}
        \SRoot{s\dvec{d}} * \SRoot{t\dvec{e}} &= \label{eq:rootrel1}
        \SRoot{r_1\dvec{f}^1} * \dots * \SRoot{r_l\dvec{f}^l}& &\text{if }\ext(\dvec{e}, \dvec{d})=0,\ 
        \dvec{d}\text{ or }\dvec{e}\text{ real},\\
        \intertext{where $\sum r_i \dvec{f}^i$ is the canonical decomposition of $s\dvec{d} + t\dvec{e}$ and}
        \SRoot{s\dvec{d}} * \SRoot{t\dvec{e}} &= \SRoot{t\dvec{e}} * \SRoot{s\dvec{d}} && \text{if }
        \ext(\dvec{d}, \dvec{e}) = \ext(\dvec{e}, \dvec{d}) = 0.\label{eq:rootrel2}
    \end{align}
    If $w=(i_1, i_2, \dots, i_r)$ is a word in vertices of $Q$ we write
    \[ \SRoot{w} := \SRoot{\epsilon_{i_1}} * \dots * \SRoot{\epsilon_{i_r}} \in \RMon(Q). \]
\end{defn}
\begin{remark}
    Note that for real Schur roots $\dvec{d}$ we have that
    \[ \SRoot{s\dvec{d}} * \SRoot{t\dvec{d}} = \SRoot{(s+t)\dvec{d}}\]
    by relation (\ref{eq:rootrel1}).
\end{remark}
We have the following observation.
\begin{propn}
    The map
    \begin{alignat*}{2}
        \Theta &\colon\quad& \RMon(Q) &\rightarrow \mathcal{CM}(Q)\\
                          &&  \SRoot{s\dvec{d}} &\mapsto \Rep{s\dvec{d}}
    \end{alignat*}
    is an epimorphism of monoids.
\end{propn}
\begin{proof}
    Since $Q$ is acyclic we have that $\Rep{\dvec{d}} \in \mathcal{CM}(Q)$ for
    each dimension vector $\dvec{d}$. The map is well-defined since the defining relations
    of $\RMon(Q)$ hold in $\mathcal{CM}(Q)$ by proposition \ref{compmon:propn:dimvecmult} and
    lemma \ref{compmon:schurrootsmult}.
\end{proof}
\begin{defn}
    We say that an element of $\RMon(Q)$ is in partial normal form if it is equal
    to $\mathcal{P}*\mathcal{R}*\mathcal{I}$ where
    \begin{align*}
        \mathcal{P}&=\SRoot{r_1\dvec{d}^1} * \dots * \SRoot{r_l\dvec{d}^l} && \text{with } r_i > 0,\ \dvec{d}^i
        \text{ preprojective and }
        \dvec{d}^i \prec_t \dvec{d}^j \text{ for all } i < j;\\
        \mathcal{R}&=\SRoot{s_1\dvec{e}^1} * \dots * \SRoot{s_m\dvec{e}^m} && \text{with } s_i > 0,\ \dvec{e}^i
        \text{ regular};\\
\mathcal{I}&=\SRoot{t_1\dvec{f}^1} * \dots * \SRoot{t_n\dvec{f}^n} && \text{with } t_i > 0,\ \dvec{f}^i
        \text{ preinjective and }
        \dvec{f}^i \prec_t \dvec{f}^j \text{ for all } i < j
    \end{align*}
    where $\dvec{d}^i, \dvec{e}^i$ and $\dvec{f}^i$ are Schur roots.
    \label{compmon:normform}
\end{defn}

\begin{remark}
    Note that if $M \prec_t N$, for two indecomposable representations
    $M, N \in \mathcal{P} \cup \mathcal{I}$, then $[M,N]^1=0$.
    Therefore, this is a similar partial normal form as in theorem 5.8 of \cite{Reineke_monoid}.
\end{remark}

\begin{lmm}
    Let $\dvec{d}^1, \dots, \dvec{d}^r$ be Schur roots such that
    \begin{align*}
        &\dvec{d}^1\text{ or }\dvec{d}^r\text{ is not regular},\\
        &\ext(\dvec{d}^1, \dvec{d}^r)\neq 0\text{ and }\\
        &\ext(\dvec{d}^i, \dvec{d}^j)=0 &&\text{ for all }
        1 \le i < j \le r\text{ with }(i,j) \neq (1,r).
    \end{align*}
    Let $s_1, \dots, s_r \in \N$ be positive integers.
    Then there is a
    permutation $\pi$ of $\{1,\dots,r\}$ such that
    such that
    \begin{align*}
    \SRoot{s_1\dvec{d}^1} * \dots * \SRoot{s_r\dvec{d}^r} &=
    \SRoot{s_{\pi(1)}\dvec{d}^{\pi(1)}} * \dots * \SRoot{s_{\pi(r)}\dvec{d}^{\pi(r)}} \in \RMon(Q),\\
    \ext(\dvec{d}^{\pi(i)}, \dvec{d}^{\pi(j)})&=0\quad\text{for all }
    1 \le i < j \le r\text{ with }(i,j) \neq (\pi^{-1}(1),\pi^{-1}(r)) \text{ and }\\
    \pi^{-1}(r)&=\pi^{-1}(1)+1.
    \end{align*}
    In other words, we can interchange the
    roots with each other such that $\SRoot{s_1 \dvec{d}^1}$ is next to $\SRoot{s_r \dvec{d}^r}$.
    \label{compmon:lmm:reorderext}
\end{lmm}
\begin{proof}
    We prove the claim by induction on $r$. If $r=2$ we are done. Now assume $r>2$.
    If there exists a $1<l<r$ such that either $\ext(\dvec{d}^l, \dvec{d}^i) = 0$ for
    all $i<l$ or $\ext(\dvec{d}^j, \dvec{d}^l) = 0$ for all $j>l$, then we
    can use relation (\ref{eq:rootrel2}) to move $\SRoot{s_l\dvec{d}^l}$ either to the left of
    $\SRoot{s_1\dvec{d}^1}$ or to the right of $\SRoot{s_r\dvec{d}^r}$ and we
    are done by induction.

    Assume now that for each $1<l<r$ there is an $i<l$ such that $\ext(\dvec{d}^l, \dvec{d}^i) \neq 0$ and
    a $j>l$ such that $\ext(\dvec{d}^j, \dvec{d}^l) \neq 0$. We therefore have a sequence
    $1=i_0 < i_1 < \dots < i_s =r$ such that $\ext(\dvec{d}^{i_{j+1}}, \dvec{d}^{i_{j}}) \neq 0$ for
    all $1 \le j < s$.

    This yields a path in the Auslander-Reiten quiver from
    $\dvec{d}^1$ to $\dvec{d}^r$.\footnote{More precisely, from the indecomposables corresponding to each.} Moreover,
    since $\ext(\dvec{d}^1, \dvec{d}^r) \neq 0$, we have a path from $\dvec{d}^r$ to
    $\dvec{d}^1$. This yields a cycle in the Auslander-Reiten quiver which involves
    at least one non-regular Schur root, a contradiction.
%
%
%
\end{proof}
\begin{theom}
    Now let $\mathcal{X}$ be any element of $\RMon(Q)$. Then
    $\mathcal{X}$ can be written in partial normal form, i.e.
    \[\mathcal{X} = \SRoot{s_1\dvec{d}^1} * \dots * \SRoot{s_r\dvec{d}^r} \in \RMon(Q), \]
    where the right hand side is in partial normal form.
    \label{compmon:hasnormalschur}
\end{theom}
\begin{proof}
  By definition, $\mathcal{X} = \SRoot{s_1\dvec{d}^1} * \dots * \SRoot{s_r\dvec{d}^r}$ for some Schur roots
  $\dvec{d}^i$.

  Let $n:= \sum_{i<j} \ext(s_i\dvec{d}^i, s_j\dvec{d}^j)$. We proof the claim by induction on $n$.
  For $n=0$ we can reorder the roots in the desired way by using relation
  (\ref{eq:rootrel2}), because there are no extensions
  between regular and preinjective, preprojective and preinjective and preprojective and regular
  Schur roots. After reordering we probably have to use relation (\ref{eq:rootrel1}) to
  obtain $\SRoot{r\dvec{d}} * \SRoot{s\dvec{d}} = \SRoot{(r+s)\dvec{d}}$ for every
  preprojective or preinjective Schur root $\dvec{d}$ to end up with an expression in partial normal form.

  Let $n \ge 1$. If
  \[
  \mathcal{X} = \SRoot{s_1\dvec{d}^1} * \dots * \SRoot{s_r\dvec{d}^r}
  \]
  cannot be reordered in the desired way by relation (\ref{eq:rootrel2}), then there exist $1\le i_0< j_0 \le r$ such that
  $\dvec{d}^{i_0}$ or $\dvec{d}^{j_0}$ is not regular,
  $\ext(\dvec{d}^{i_0}, \dvec{d}^{j_0}) \neq 0$ and $\ext(\dvec{d}^i, \dvec{d}^j) = 0$ for all
  $i_0 \le i<j \le j_0$ with $(i,j) \neq (i_0, j_0)$.  By lemma \ref{compmon:lmm:reorderext} we can
  assume that $j_0=i_0+1$.

  Now we can apply
  relation (\ref{eq:rootrel1}) to obtain that
  \[\SRoot{s_{i_0}\dvec{d}^{i_0}}*\SRoot{s_{j_0}\dvec{d}^{j_0}}=
  \SRoot{t_1\dvec{e}^1} * \dots *\SRoot{t_r\dvec{e}^r},
  \]
  where $\sum t_j \dvec{e}^j$ is the canonical decomposition of $s_{i_0} \dvec{d}^{i_0} + s_{j_0} \dvec{d}^{j_0}$.
  Replacing $\SRoot{s_{i_0}\dvec{d}^{i_0}}*\SRoot{s_{j_0}\dvec{d}^{j_0}}$ by
  $\SRoot{t_1\dvec{e}^1} * \dots *\SRoot{t_r\dvec{e}^r}$ makes $n$ smaller, as the following calculation shows, and
  we are done by induction.
  \begin{multline*}
      \sum_{i<j, \{i,j\} \cap \{i_0,j_0\} = \varnothing  } \ext(s_i\dvec{d}^i, s_j\dvec{d}^j) +
      \sum_{i < i_0, j} \ext(s_i\dvec{d}^i,t_j\dvec{e}^j) + \sum_{j > j_0, i} \ext(t_i\dvec{e}^i, s_j\dvec{d}^j) \\
      \stackrel{\ref{compmon:lmm:reorderext}}{\le}
      \sum_{i<j, \{i,j\} \cap \{i_0,j_0\} = \varnothing  } \ext(s_i\dvec{d}^i, s_j\dvec{d}^j) +
      \sum_{i < i_0} \left( \ext(s_i\dvec{d}^i, s_{i_0}\dvec{d}^{i_0})+\ext(s_i\dvec{d}^i, s_{j_0}\dvec{d}^{j_0}\right) \\
      + \sum_{j> j_0} \left(\ext(s_{i_0}\dvec{d}^{i_0}, s_j\dvec{d}^j) + \ext(s_{j_0}\dvec{d}^{j_0}, s_j\dvec{d}^j) \right)\\
      =\sum_{i<j} \ext(s_i\dvec{d}^i, s_j\dvec{d}^j) - \ext(s_{i_0}\dvec{d}^{i_0}, s_{j_0}\dvec{d}^{j_0}) < n.
  \end{multline*}
\end{proof}
As a corollary we obtain a partial normal form in $\mathcal{CM}(Q)$.
\begin{cory}
    Let $\mathcal{A}_w \in \mathcal{CM}(Q)$ for a word $w$ in vertices of $Q$. Then $\mathcal{A}_w$ can be written as
    $\mathcal{P}*\mathcal{R}*\mathcal{I}$ where
    \begin{align*}
        \mathcal{P}&=\Rep{r_1\dvec{d}^1} * \dots * \Rep{r_l\dvec{d}^l} && \text{with } r_i > 0,\ \dvec{d}^i
        \text{ preprojective and }
        \dvec{d}^i \prec_t \dvec{d}^j \ \forall\ i < j;\\
        \mathcal{R}&=\Rep{s_1\dvec{e}^1} * \dots * \Rep{s_m\dvec{e}^m} && \text{with } s_i > 0,\ \dvec{e}^i
        \text{ regular};\\
        \mathcal{I}&=\Rep{t_1\dvec{f}^1} * \dots * \Rep{t_n\dvec{f}^n} && \text{with } t_i > 0,\ \dvec{f}^i
        \text{ preinjective and }
        \dvec{f}^i \prec_t \dvec{f}^j \ \forall\ i < j.
    \end{align*}
\end{cory}
\begin{proof}
    By definition of the morphism $\Theta \colon \RMon(Q) \rightarrow \mathcal{CM}(Q)$ we have that
    $\Theta(\{w\}) = \mathcal{A}_w$. Since $\SRoot{w}$ can be written in partial normal form and
    $\Theta$ is a homomorphism sending $\SRoot{r\dvec{d}}$ to $\Rep{r\dvec{d}}$ for each Schur root
    $\dvec{d}$ and each $r\in\N$ we obtain the result.
\end{proof}
\begin{cory}
  If $Q$ is Dynkin, then
    \begin{alignat*}{2}
        \Theta &\colon\quad& \RMon(Q) &\rightarrow \mathcal{CM}(Q)\\
                          &&  \SRoot{s\dvec{d}} &\mapsto \Rep{s\dvec{d}}
    \end{alignat*}
    is an isomorphism and the partial normal form is a normal form. In particular
    $\mathcal{CM}(Q)$ is independent of the base field.
\end{cory}
\begin{proof}
    By definition, there are no preinjective or regular Schur roots, therefore each element
    of $\RMon(Q)$ can be written as
        $\SRoot{r_1\dvec{d}^1} * \dots * \SRoot{r_l\dvec{d}^l}$ with $r_i > 0$, $\dvec{d}^i$
        preprojective and
        $\dvec{d}^i \prec_t \dvec{d}^j$ for all $i < j$. For each $\dvec{d}^i$ there is a unique
        indecomposable representation $M_i$ without self-extensions such that $\overline{\Orbit_{M_i}} = \Rep{\dvec{d}^i}$.
        Therefore, $\Rep{r_i\dvec{d}^i} = \overline{\Orbit_{M_i^{r_i}}}$. By the condition on the partial normal
        form we have that $[M_i, M_j]^1 = 0$ for all $i < j$. We can therefore apply lemma \ref{compmon:noext->dsum}
        and obtain
        \[
        \Theta(\SRoot{r_1\dvec{d}^1} * \dots * \SRoot{r_l\dvec{d}^l}) = 
        \Rep{r_1\dvec{d}^1} * \dots * \Rep{r_l\dvec{d}^l} =
        \overline{\Orbit_{M_1^{r_1}}} * \dots * \overline{\Orbit_{M_l^{r_l}}} =
        \overline{\Orbit_{\bigoplus M_i^{r_i}}}.
        \]
        Since every representation in $\mathcal{CM}(Q)$ is uniquely given by the orbit closure of some representation $M$,
        and the isomorphism class of $M$ is uniquely given by some Schur roots with multiplicities
        we have that $\Theta$ is bijective.
\end{proof}

\section{Extended Dynkin Case}
%
Let $K$ be an algebraically closed field and let $Q$
be a connected, acyclic, extended Dynkin quiver. We show that $\mathcal{CM}_K(Q)$ is a quotient
of $\RMon(Q)$ by two more classes of relations which are independent of $K$. As a consequence we will obtain that
$\mathcal{CM}(Q)$ does not depend on the base field $K$.

For a connected, acyclic, extended Dynkin quiver the tubes are indexed by $\mathbb{P}_K^1$. There is
exactly one isotropic Schur root, $\dvec{\delta}$. Let
$\mathbb{H}_K\subseteq\mathbb{P}_K^1$ be the set of indices of the homogeneous tubes, which is an
open subset. For each $x \in \mathbb{H}_K$ there is a unique (up to isomorphism) regular
simple representation in the corresponding tube $\mathcal{T}_x$. Let us call this representation $R_x$.

Let $x \in \Proj^1_K$ such that $\mathcal{T}_x$ is inhomogeneous. Denote by $\Pi_x$ the
set of isomorphism classes in this tube. Since $\mathcal{T}_x$ is equivalent to
$\repK{C_{k-1}}{K}$ where $k = \rank \mathcal{T}_x$ we have that
\[\Pi_x \cong \Set{ \pi = (\pi^{(0)}, \dots, \pi^{(k-1)}) | \pi^{(i)} \text{ is a partition}}.\]
In the following we identify $\Pi_x$ with this set.

K. Bongartz and D. Dudek show in \cite{Bongartz_decomp} the following.
\begin{theom}
    Let $M$ be a representation without homogeneous direct summands and $\dvec{d}=\dimve M + r \dvec{\delta}$ for
    $\dvec{\delta}$ being the isotropic Schur root and some $r \in \N$. Then the set
    \[ D(M,r):=\Set{ N \in \rep{\dvec{d}} |
	\begin{array}{l}
		N \cong M \oplus \bigoplus_{j=1}^r R_{x_j} \ \text{with}\
    x_j \\ \text{pairwise different elements of}\  \mathbb{H}_K
    \end{array}
    }
\]
    is a smooth, locally closed subset of $\rep{\dvec{d}}$.
    \label{compmon:tamesmoothopensubset}
\end{theom}
Note that $D(M,r)$ is the decomposition class corresponding to $(\mu, \sigma)$, where $\mu$ is
the isomorphism class of $M$ and $\sigma = \{ ((1),1), ((1),1), \dots, ((1),1) \}$.

First, we look at an inhomogeneous tube of rank $k$
and use our result about the composition monoid of a cyclic quiver.
So let $T_1, \dots, T_k$ be the regular simples of an inhomogeneous tube. We know that
$\dimve T_j$ is a Schur root and $\Orbit_{T_j}$ is open in $\rep{\dimve T_j}$. Moreover, the tube
is closed under extensions. By work of K. Bongartz \cite{Bongartz_degenstame} and G. Zwara \cite{Zwara_degensextdyn}
we also know
that the $\Hom$-order and the $\Ext$-order agree on $\repK{Q}{K}$.
We have the following.
\begin{lmm}
    Let $M$ and $N$ be $K$-representations of $Q$
    in the same inhomogeneous tube $\mathcal{T}_x$ of rank $k$. Then $M \le_{\text{deg}} N$ if and only
    if $M'\ledeg N'$ in
    $\repK{C_{k-1}}{K}$, where $M'$ and $N'$ are the images of $M$ and $N$ under
    the equivalence $\mathcal{T}_x \cong \repK{C_{k-1}}{K}$.
    \label{compmon:cyclicistameinhom}
\end{lmm}
\begin{proof}
    If $M'\le_{\text{deg}} N'$, then this degeneration is given by successive extensions. These can be
    transformed to successive extensions in $\repK{Q}{K}$ to obtain a degeneration from $M$ to $N$.

    Now, if $M \le_{\text{deg}} N$, then $[M,B] \le [N, B]$ for all representations $B \in \repK{Q}{K}$.
    But then, by the equivalence, we also have $[M',B'] \le [N',B']$ for all representations
    $B' \in \repK{C_{k-1}}{K}$. Hence $M' \le_{\text{deg}} N'$ since in $\repK{C_{k-1}}{K}$ the degeneration order
    agrees with the $\Hom$-order.
\end{proof}
\begin{theom}
    Let $M$ and $N$ be $K$-representations in an inhomogeneous tube $\mathcal{T}_x$ of rank $k$.
    Then there is a representation $E \in \mathcal{T}_x$ such that
    \[ \overline{\Orbit_M} * \overline{\Orbit_N} = \overline{\Orbit_E},\]
    i.e. for an inhomogeneous tube we have generic extensions. Moreover, the
    isomorphism class $\pi \in \Pi_x$ of $E$ only depends on the isomorphism classes $\pi', \pi'' \in \Pi_x$ of
    $M$ and $N$ and not on the field $K$.
    \label{compmon:genericextsregular}
\end{theom}
\begin{proof}
    Let $M', N' \in \repK{C_{k-1}}{K}$ be the images of $M$ and $N$ under the equivalence
    $\mathcal{T}_x \cong \repK{C_{k-1}}{K}$.
    In $\repK{C_{k-1}}{K}$ we have generic extensions. Therefore, there is an $E' \in \repK{C_{k-1}}{K}$ such
    that $E' \in \mathcal{E}(\Orbit_{M'}, \Orbit_{N'})$ and $E' \ledeg X'$ for all
    $X' \in \mathcal{E}(\Orbit_{M'}, \Orbit_{N'})$.
    Let $E \in \repK{Q}{K}$ be the representation corresponding to $E'$ under the equivalence.
    By lemma \ref{compmon:cyclicistameinhom} $E \ledeg X$ for all $X \in \mathcal{E}(\Orbit_M,\Orbit_N)$
    since $\mathcal{E}(\Orbit_M,\Orbit_N) \subset \mathcal{T}_x$.
    Therefore, $\Orbit_E$ is dense in $\mathcal{E}(\Orbit_M,\Orbit_N)$.
    Now we can apply lemma \ref{compmon:openopenisopen}:
    \[ \overline{\Orbit_M} * \overline{\Orbit_N} = \overline{\mathcal{E}(\Orbit_M,\Orbit_N)} = \overline{\Orbit_E}.\]

    Moreover, the set of isomorphism classes
    $[\overline{\mathcal{E}(\Orbit_{M'},\Orbit_{N'})}] \subset \Pi_x$ only depends on the isomorphism classes
    $\pi'$ of $M'$ and $\pi''$ of $N'$
    by lemma \ref{cyclic:lmm:filtrindept} since $\overline{\Orbit_{M'}}$ and
    $\overline{\Orbit_{N'}}$ are
    elements of the generic extension monoid $\mathcal{M}(C_{k-1})$.
    Therefore, the
    isomorphism class $\pi \in \Pi_x$ of $E'$ only depends on the isomorphism classes $\pi', \pi'' \in \Pi_x$
    since it is the isomorphism class in $[\overline{\mathcal{E}(\Orbit_{M'},\Orbit_{N'})}]$ with
    the smallest endomorphism ring dimension. Now $\pi$ is by definition the isomorphism class of $E$
    and this proves the claim.
\end{proof}
\begin{cory}
    Let $\dvec{e}_1,\dots, \dvec{e}_l$ be regular Schur roots living in one inhomogeneous tube $\mathcal{T}_x$ and
    $r_1, \dots, r_l \in \N$.
    Then there is a regular
    representation $M \in \mathcal{T}_x$ whose isomorphism class $\pi \in \Pi_x$ only depends
    on $r_1\dvec{e}_1,\dots, r_l\dvec{e}_l$ and not on $K$ such that
    \[ \rep{r_1\dvec{e}_1} * \dots * \rep{r_l\dvec{e}_l} = \overline{\Orbit_M}. \]
    \label{compmon:cory:rootmultonetube}
\end{cory}
\begin{proof}
    Each $\Rep{r_i\dvec{e}_i}$ is given by $\overline{\Orbit_{M_i^{r_i}}}$ for some indecomposable
    $K$-representation $M_i \in \mathcal{T}_x$ whose
    isomorphism class only depends on $\dvec{e}_i$. Iteratively applying theorem $\ref{compmon:genericextsregular}$
    yields the result.
\end{proof}
We therefore have additional relations in $\mathcal{CM}(Q)$ which seem not to be a consequence
of the relations of proposition \ref{compmon:propn:dimvecmult}. It would be interesting to
make this precise. We add them to
the defining relations of the Schur root monoid and call the resulting monoid $\EMon(Q)$.
\begin{defn}
  We define $\EMon(Q)$ as the quotient of $\RMon(Q)$ by the following type of relation. For each $x \in \Proj^1_K$
  such that $\mathcal{T}_x$ is inhomogeneous we identify
  \begin{align}
    \SRoot{r_1\dvec{d}_1} * \dots * \SRoot{r_l\dvec{d_l}} &= \SRoot{s_1 \dvec{e}_1} * \dots *
    \SRoot{s_m\dvec{e}_m}
    \label{eq:rootrel3}
  \end{align}
  if $\dvec{d}_1, \dots, \dvec{d}_l$ and $\dvec{e}_1, \dots, \dvec{e}_m$ live in the inhomogeneous tube
  $\mathcal{T}_x$ and there is a representation $M \in \mathcal{T}_x$ such that
  \[\SRoot{r_1\dvec{d}_1} * \dots * \SRoot{r_l\dvec{d_l}} = \overline{\Orbit_M} = \SRoot{s_1 \dvec{e}_1} * \dots *
    \SRoot{s_m\dvec{e}_m}.
    \]
    In this case we denote $\SRoot{r_1\dvec{d}_1} * \dots * \SRoot{r_l\dvec{d_l}}$ by $\SRoot{(\pi,x)}$,
    $\pi \in \Pi_x^s$ being the isomorphism class of $M$.
    Note that for every $\pi \in \Pi_x^s$ we have that $\SRoot{(\pi,x)} \in \EMon(Q)$ and that
    this relation does not depend on $K$.
\end{defn}
\begin{remark}
  We will later show that, for $Q$ acyclic, extended Dynkin, $\Q\EMon(Q) \cong \mathcal{C}_0(Q)$.
\end{remark}
\begin{propn}
    The morphism $\Theta \colon \RMon(Q) \rightarrow \mathcal{CM}(Q)$ factors via $\EMon(Q)$.
    We denote the induced morphism $\EMon(Q) \rightarrow \mathcal{CM}(Q)$ again by
    $\Theta$.
\end{propn}
\begin{proof}
    By corollary \ref{compmon:cory:rootmultonetube} the additional relations of $\EMon(Q)$
    hold in $\mathcal{CM}(Q)$ and the claim follows.
\end{proof}
\begin{remark}
    Note
    that, for an inhomogeneous tube $\mathcal{T}_x$,
    $\Theta(\SRoot{(\pi,x)}) = \overline{\Orbit_M}$ for $M \in \mathcal{T}_x$ a $K$-representation of
    isomorphism class $\pi \in \Pi^s_x$.
\end{remark}
\begin{lmm}
    Let $x_1, \dots, x_r \in \Proj^1_K$ be the indices of the inhomogeneous tubes. Then
    every element $\mathcal{X} \in \EMon(Q)$ can be written in the form
    $\mathcal{X} = \mathcal{P}*\mathcal{C}_1*\dots*\mathcal{C}_r * \mathcal{R} *\mathcal{I}$ where
    \begin{align*}
        \mathcal{P}&=\SRoot{s_1\dvec{p}^1} * \dots * \SRoot{s_k\dvec{p}^k} && \text{with } s_i > 0,\ \dvec{p}^i
        \text{ preprojective and }
        \dvec{p}^i \prec_t \dvec{p}^j \text{ for all } i < j;\\
        \mathcal{C}_i&=\SRoot{(\pi,x_i)} && \text{for a } \pi \in \Pi^s_{x_i};\\
        \mathcal{R}&=\SRoot{\lambda_1\dvec{\delta}} * \dots * \SRoot{\lambda_l\dvec{\delta}} && \text{with }
        \lambda=(\lambda_1\ge \dots \ge \lambda_l)
        \text{ a partition};\\
        \mathcal{I}&=\SRoot{t_1\dvec{q}^1} * \dots * \SRoot{t_m\dvec{q}^m} && \text{with } t_i > 0,\ \dvec{q}^i
        \text{ preinjective and }
        \dvec{q}^i \prec_t \dvec{q}^j \text{ for all } i < j.
    \end{align*}
    \label{compmon:lmm:normformemon}
\end{lmm}
\begin{proof}
    By lemma \ref{compmon:hasnormalschur} $\mathcal{X}$ can be written as $\mathcal{P}*\mathcal{R}'*\mathcal{I}$ where
    $\mathcal{P}$ and $\mathcal{I}$ are already in the desired form and
    $\mathcal{R}' = \SRoot{i_1 \dvec{e}^1} *\dots* \SRoot{i_n \dvec{e}^n}$ for some regular Schur roots
    $\dvec{e}^1, \dots, \dvec{e}^n$. For each $j$ we have that either $\dvec{e}^j \in \mathcal{T}_x$
    for some $x \in \{x_1,\dots, x_r\}$ or that $\dvec{e}^j = \dvec{\delta}$. Since there
    are no extensions between different tubes we can use relation (\ref{eq:rootrel2}) to
    reorder them in $\EMon(Q)$ in such a way that roots from the same inhomogeneous tube are next to each other as
    are multiples of $\dvec{\delta}$.
    Using relation (\ref{eq:rootrel3}) we obtain the desired form.
\end{proof}
\begin{remark}
    We will show in chapter \ref{chap:extdynkin} that this is a normal form.
\end{remark}
We finally define the last monoid, $\TEMon(Q)$, which will be isomorphic to $\mathcal{CM}(Q)$.
\begin{defn}
    We define $\TEMon(Q)$ as the quotient of $\EMon(Q)$ by the relations
    \begin{align}
        \SRoot{s\dvec{\delta}}*\SRoot{t\dvec{\delta}} &= \SRoot{(s+t)\dvec{\delta}}&& \text{ for all } s, t > 0.
        \label{eq:rootrel4}
    \end{align}
\end{defn}
\begin{propn}
    The morphism $\Theta \colon \EMon(Q) \rightarrow \mathcal{CM}(Q)$ induces
    a morphism
    \[\overline{\Theta} \colon \TEMon(Q) \rightarrow \mathcal{CM}(Q).
    \]
\end{propn}
\begin{proof}
    Relation (\ref{eq:rootrel4}) holds in $\mathcal{CM}(Q)$ by proposition \ref{compmon:propn:dimvecmult} and
    the claim follows.
\end{proof}
Now we obtain a normal form in $\TEMon(Q)$.
\begin{theom}
    Every element $\mathcal{X} \in \TEMon(Q)$ can be uniquely written as
    $\mathcal{P}*\mathcal{C}_1 * \dots \mathcal{C}_r * \SRoot{l\dvec{\delta}} * \mathcal{I}$ where
    \label{compmon:firstnormal}
    \begin{align*}
        \mathcal{P}&=\SRoot{s_1\dvec{p}^1} * \dots * \SRoot{s_k\dvec{p}^k} && \text{with } s_i > 0,\ \dvec{p}^i
        \text{ preprojective and }
        \dvec{p}^i \prec_t \dvec{p}^j \text{ for all } i < j;\\
        \mathcal{C}_i&=\SRoot{(\pi_i,x_i)} && \text{for a } \pi_i \in \Pi^s_{x_i};\\
         l &\ge 0; \\
        \mathcal{I}&=\SRoot{t_1\dvec{q}^1} * \dots * \SRoot{t_m\dvec{q}^m} && \text{with } t_i > 0,\ \dvec{q}^i
        \text{ preinjective and }
        \dvec{q}^i \prec_t \dvec{q}^j \text{ for all } i < j.
    \end{align*}
    Moreover, $\overline{\Theta} \colon \TEMon(Q) \rightarrow \mathcal{CM}(Q)$ is an isomorphism.
    Since the relations of $\TEMon(Q)$ do not depend on $K$ we have that $\mathcal{CM}(Q)$
    is independent of $K$.
\end{theom}
\begin{proof}
    By using lemma \ref{compmon:lmm:normformemon} and relation (\ref{eq:rootrel4}) we have that
    every element can be written in this way. We have to show that this form is unique.
    Let $\mathcal{X} = \mathcal{P} * \mathcal{C}_1 * \dots * \mathcal{C}_r * \SRoot{l\dvec{\delta}} * \mathcal{I}$
    as above. Then for each $i$ there is a unique (up to isomorphism) indecomposable
    preprojective $K$-representation $P_i$
    of dimension vector $\dvec{p}^i$. We have that $\overline{\Orbit_{P_i^{s_i}}} = \Rep{s_i \dvec{p}^i}$.
    Dually, for each $i$ there is unique indecomposable preinjective $K$-representation $I_i$ of
    dimension vector $\dvec{q}^i$ and $\overline{\Orbit_{I_i^{t_i}}} = \Rep{t_i\dvec{q}^i}$.
    Finally, for each $1 \le i \le r$ there is a unique $K$-representation $M_i \in \mathcal{T}_{x_i}$
    of isomorphism class $\pi_i$. Let $P:=\bigoplus P_i^{s_i}$, $I:= \bigoplus I_i^{t_i}$
    and $M := \bigoplus M_i$.
    By the Krull-Remak-Schmidt theorem $\mathcal{P}$, $\mathcal{C}_1$, \ldots, $\mathcal{C}_r$ and $\mathcal{I}$ are
    uniquely determined by $[P]$, $[M]$ and $[I]$.
 
    By applying corollary \ref{compmon:noext->dsum} and using that $\overline{\Theta}$ is a homomorphism we have that
    \begin{multline*}
        \overline{\Theta}(\mathcal{X}) =
        \overline{\Theta}(\mathcal{P}) * \overline{\Theta}(\mathcal{C}_1)* \dots *\overline{\Theta}(\mathcal{C}_r)*
        \Rep{l\dvec{\delta}} * \overline{\Theta}(\mathcal{I})\\
        =\overline{\Orbit_{P}} * \overline{\Orbit_{M}} * \Rep{l \dvec{\delta}} * \overline{\Orbit_{I}}=
        \overline{\Orbit_{P \oplus M}} * \Rep{l \dvec{\delta}} * \overline{\Orbit_{I}}.
    \end{multline*}
    If we show that $[P]$, $[M]$, $[I]$ and $l$ are uniquely determined by the element
    $\mathcal{A} := \overline{\Theta}(\mathcal{X})$, then we have that $\mathcal{X}$ can
    be uniquely written as $\mathcal{P}*\mathcal{C}_1*\dots*\mathcal{C}_r * \SRoot{l\dvec{\delta}} * \mathcal{I}$
    and that $\overline{\Theta}$ is injective.

    Assume therefore that there is a preprojective representation $P'$, a regular inhomogeneous representation $M'$,
    a preinjective representation $I'$ and an integer $l'$ such that
    \[
        \overline{\Orbit_{P \oplus M}} * \Rep{l \dvec{\delta}} * \overline{\Orbit_{I}} =
        \overline{\Orbit_{P' \oplus M'}} * \Rep{l' \dvec{\delta}} * \overline{\Orbit_{I'}}.
    \]
 
    By
    theorem \ref{compmon:tamesmoothopensubset} we have that the set $D(0,l)$ is smooth and locally closed in
    $\Rep{\dvec{\delta}}^l$. Moreover, it is open
    in $\Rep{\dvec{\delta}}^l$. Now we can apply lemma \ref{compmon:openopenisopen}
    to obtain that $D(P \oplus M \oplus I, l)$ is dense in $\mathcal{A}$. Since it is also
    locally closed by theorem \ref{compmon:tamesmoothopensubset} we have that it is open
    in $\mathcal{A}$. With the same arguments we have that
    $D(P' \oplus M' \oplus I', l')$ is open in $\mathcal{A}$.
    Since $\mathcal{A}$ is irreducible these two open sets have to intersect. Hence there are
    $x_1,\dots, x_l \in \mathbb{H}_K$ and $y_1,\dots, y_{l'} \in \mathbb{H}_K$ such that
    \[P \oplus M \oplus I \oplus \bigoplus_{i=1}^l R_{x_i} \cong P' \oplus M' \oplus I'
    \oplus \bigoplus_{i=1}^{l'} R_{y_i}.\]
    By the Krull-Remak-Schmidt theorem this yields $P \cong P'$, $R \cong R'$, $I \cong I'$ and $l = l'$.

    Therefore, we have uniqueness of the expression and that $\overline{\Theta}$ is injective. Since
    it is also surjective by definition we have that it is an isomorphism.
\end{proof}
\begin{cory}
    Every element $\mathcal{A}$ of the composition monoid $\mathcal{CM}(Q)$ can be written in the form
    \[ \overline{\Orbit_{P\oplus M}}  * \rep{\dvec{\delta}}^l * \overline{\Orbit_I}, \]
    where $P$ is preprojective, $M$ is regular having no homogeneous summand, $I$ is
    preinjective and $l \in \N$. Moreover, $P$, $R$, $I$ and $l$ are, up to
    isomorphism, uniquely determined by $\mathcal{A}$.
\end{cory}

\chapter{Geometry of Quiver Flag Varieties}
\label{chap:geometry}
\dictum[H.L. Mecken]{For every complex problem, there is a solution that is simple, neat, and wrong.}
\vspace{1cm}
In this chapter we will talk about the geometry of quiver Grassmannians and
quiver flag varieties. Most times when we use algebraic geometry we take the
functorial viewpoint, i.e. we identify a scheme $X$ with its functor of points
$R \mapsto \Hom(\Spec R, X)$ from the category of commutative rings to sets. A good reference for
this is \cite{DG}.
\section{Basic Constructions}
We want to define some functors from commutative rings to sets
which are schemes, ending up with the quiver Grassmannian. These should turn out to be
useful to define morphisms between schemes naturally coming up in representation theory of quivers.
All of this section should be standard. We say that a scheme $X$ over some field $K$ is a variety
if it is separated, noetherian and of finite type over $K$.
\begin{defn}
    Let $d, e \in \N$. We define the scheme $\Sch{Hom}(d,e)$ via its functor of points
    \[\Sch{Hom}(d,e) (R) := \Hom_R (R^d, R^e).\]
    For any $I \subset \{1,\dots,d+e\}$ with $|I| = c \le d+e$ we have a morphism
    \[\rho_I \colon \Sch{Hom}(d,d+e) \rightarrow \Sch{Hom}(d,c)\] given by
    removing the rows not indexed by $I$ and a morphism
    \[\delta_I \colon \Sch{Hom}(d+e, d) \rightarrow \Sch{Hom}(c,d)\] given by
    removing columns not indexed by $I$. For $f \in \Sch{Hom}(d,d+e)(R)$ we
    write $f_I$ for $\rho_I (f)$ and for $g \in \Sch{Hom}(d+e, d) (R)$ we write
    $g_I$ for $\delta_I (g)$.
\end{defn}
The scheme $\Sch{Hom}(d,e)$ is an affine space since it is obviously isomorphic to
\[
\Spec \Z[X_{ij}]_{\substack{1\le i \le e\\ 1 \le j \le d}}.
\]
The maps $\rho_I$ and $\delta_I$ are natural transformations,
hence morphism of schemes.

Now we want to define open subschemes resembling the injective and surjective, or more general $\rank r$
elements of the homomorphism spaces.
\begin{defn}
    Let $d, e, r \in \N$. We define the schemes $\Sch{Hom}(d,e)_r$, $\Sch{Inj}(d,d+e)$ and
    $\Sch{Surj}(d+e, d)$ via their functors of points
\begin{align*}
    \Sch{Hom}(d,e)_r(R)\; & := \Set{ f\in \Sch{Hom}(d,e) (R) |
    \Bild f \text{ is a direct summand of } R^e \text{ of } \rank r},\\
    \Sch{Inj}(d,d+e)(R)\; & := \Set{ f\in \Sch{Hom}(d,d+e) (R) | f \text{ is a split injection}}\\
    & \; = \Sch{Hom}(d,d+e)_d(R)
    \intertext{and}
    \Sch{Surj}(d+e,e)(R) \; & := \Set{ g\in \Sch{Hom}(d+e,e) (R)| g \text{ is surjective} }\\
    & \; = \Sch{Hom}(d+e,e)_e(R).
\end{align*}
\end{defn}
\begin{remark}
    All elements of $\Sch{Surj}(d+e,e)(R)$ are automatically split since $R^e$ is a projective module.
\end{remark}
\begin{propn}
    Let $d, e, r \in \N$.
    We have that $\Sch{Hom}(d,e)_r$ is a locally closed subscheme of $\Sch{Hom}(d,d+e)$.
    Moreover, $\Sch{Inj}(d,d+e)$ is an open subscheme of $\Sch{Hom}(d,d+e)$ and $\Sch{Surj}(d+e,e)$ is an
    open subscheme of $\Sch{Hom}(d+e,e)$.
\end{propn}
\begin{proof}
    The functor $\Sch{Hom}(d,e)_r$ is a locally closed subscheme of $\Sch{Hom}(d,e)$ if for each
    $f \in \Sch{Hom}(d,e) (R)$ there are two ideals $N_1$ and $N_2$ such that
    for each $R$-algebra $S$, $f \tensor S \in \Sch{Hom}(d,e)_r (S)$ if and only if
    $(N_1 \tensor S) S =S$ and $N_2 \tensor S = (0)$ (\cite[I, \S 1, 3.6]{DG} and
    \cite[I, \S 2, 4.1]{DG}).

    Let $f \in \Sch{Hom}(d,e)(R)$. Let $N_1:= \mathcal{F}_{e-r} (R^e/ \Bild(f))$ and
    $N_2 := \mathcal{F}_{e-r+1}(R^e/ \Bild(f))$, the Fitting
    ideals of the finitely
    presented module $R^e/ \Bild(f)$.\footnote{For more on Fitting ideals see appendix \ref{app:fitting}.} Then
    $(N_1 \tensor S) S =S$ and $N_2 \tensor S = (0)$ if and only if $S^e/ \Bild(f \tensor S)$ is projective
    of $\rank e-r$ (\cite{Campillo_projfittideal}), if and only if $\Bild(f\tensor S)$
    is a direct summand of $S^e$ of rank $r$.

    The subschemes $\Sch{Inj}(d,d+e)$ and $\Sch{Surj}(d+e,e)$ are open since $N_2$ will be always $(0)$
    and only the open condition remains.
\end{proof}

\begin{defn} Let $d, e \in \N$. We define the scheme $\Sch{Ex}(d,e)$ via its functor of points
    \[ \Sch{Ex}(d,e)(R) := \Set{ (f,g) \in (\Sch{Inj}(d,d+e) \times \Sch{Surj}(d+e,e)) (R) | g \circ f = 0}. \]
\end{defn}
\begin{propn}
    Let $d, e \in \N$.
    $\Sch{Ex}(d,e)$ is a closed subscheme of $\Sch{Inj}(d,d+e) \times \Sch{Surj}(d+e,e)$ and therefore
    a locally closed subscheme of $\Sch{Hom}(d,d+e) \times \Sch{Hom}(d+e,e)$.
\end{propn}
\begin{proof}
    Let $R$ be a commutative ring and
    $(f,g) \in \Sch{Inj}(d,d+e)(R) \times \Sch{Surj}(d+e,e)(R)$. Let $I$
    be the ideal generated by the entries of the matrix of $g\circ f$. Then
    $(f\tensor S, g\tensor S) \in \Sch{Ex}(d,e)(S)$ if and only if $I \tensor S = (0)$ for
    all $R$-algebras $S$.
\end{proof}
\begin{propn}
    Let $d, e \in \N$ and let $R$ be a commutative ring.
    Then the elements $(f,g) \in \Sch{Ex}(d,e)(R)$ are exactly the short exact sequences
    $0 \rightarrow R^d \overset{f}{\rightarrow} R^{d+e} \overset{g}{\rightarrow} R^e \rightarrow 0$.
\end{propn}
\begin{proof}
    Obviously, every such short exact sequence given by $(f,g)$ is in $\Sch{Ex}(d,e)(R)$. So let
    $(f,g) \in \Sch{Ex}(d,e)(R)$. We have to show that the sequence
    $0 \rightarrow R^d \overset{f}{\rightarrow} R^{d+e} \overset{g}{\rightarrow} R^e \rightarrow 0$
    is exact. The $R$-homomorphisms $f$ and $g$ are split by definition, so we can choose sections:
    \[
    \xymatrix{ 0 \ar[r]& R^d \ar[r]^f & R^{d+e} \ar[r]^g \ar@/^/@{.>}[l]^r& R^e \ar[r] \ar@/^/@{.>}[l]^s & 0}.
    \]
    Now the map $\left(\begin{smallmatrix}
        f & s
    \end{smallmatrix}\right) \colon R^{d+e} \rightarrow R^{d+e}$ is a split monomorphism since it has
    $\left(\begin{smallmatrix}
        r-rsg \\ g
    \end{smallmatrix}\right)$ as a section. Therefore, the cokernel of $\left(\begin{smallmatrix}
                f & s
    \end{smallmatrix}\right)$ is projective of rank $0$, hence equal to $0$
    and the sequence is a short exact sequence.

\end{proof}
\begin{defn}
    Let $n \in \N$. We define the group scheme $\GL_n$ via its functor of points
    \[ \GL_n (R) := \GL(R^n). \]
    $\GL_n$ is a smooth group scheme (\cite[II, \S 1, 2.4]{DG}) over $\Z$.
\end{defn}
\begin{propn}
    Let $d, e \in \N$. Then
    $\GL_d \times \GL_e$ acts on $\Sch{Ex}(d,e)$ via the action on $R$-valued points
    \[ (\sigma, \tau) \cdot (f,g) := (f \circ \sigma^{-1}, \tau \circ g) \]
    for any $\sigma \in \GL_d(R)$, $\tau \in \GL_e(R)$ and $(f,g) \in \Sch{Ex}(d,f)(R)$.
\end{propn}
\begin{proof}
    obvious.
\end{proof}
We recall the definition of a principal $G$-bundle in our setting, which is just a bundle being
trivial for the Zariski topology.
\begin{defn}
    Let $X$ be a scheme, $P$ a scheme with the action $\mu$ of a group scheme $G$ and a projection
    $\pi \colon P \rightarrow X$ which is $G$-invariant. Then $(P,X,\mu,\pi)$, or simply $\pi$,
    is a principal $G$-bundle if for each ring $R$ and for every morphism $\Spec R \rightarrow X$
    there are elements $s_1, \dots, s_k \in R$ which generate $R_R$ such that
    $\Spec R_{s_i} \times_X P$ is trivial, i.e. isomorphic to $\Spec R_{s_i} \times G$, for each $i$.
\end{defn}
\begin{theom}
    Let $d, e \in \N$.
    The natural projection $\pi \colon \Sch{Ex}(d,e) \rightarrow \Sch{Inj}(d,d+e)$ is a principal
    $\GL_e$-bundle.
\end{theom}
\begin{proof}
    Let $R$ be a ring and $f \in \Sch{Inj}(d,d+e)(R)$. For every
    $I \subset \{1,\dots,d+e\}$ with $|I| = d$ let $s_I:= \Delta_I(f)$,
    where $\Delta_I(f) := \det f_I$.
    The elements $\Delta_I := \Delta_I(f)$
    generate $R_R$. We have to show that
    \[
    D(\Delta_I) \times_\Sch{Inj} \Sch{Ex}(d,e)
    \]
    is trivial, where $D(\Delta_I) = \Spec R_{\Delta_i}$. We construct
    an isomorphism to $D(\Delta_I) \times \GL_e$. For $I \subset \{1,\dots,d+e\}$ denote by
    $\theta_I \colon \Sch{Hom}(d,e) \rightarrow \Sch{Hom}(d+e,e)$ the section
    of $\delta_I$ replacing the columns not indexed by $I$ with $-\id_{R^e}$.
    Let $f \in D(\Delta_I)$. Then
    $f_I$ is invertible. Let $\phi \colon D(\Delta_I) \rightarrow \Sch{Surj}(d+e,e)$
    be the morphism sending $f$ to $\theta_I(f_{I^C} f_I^{-1})$. Then
    $\phi(f) f = f_{I^C} f_I^{-1} f_I - f_{I^C} = 0$. Here, $I^C$ denotes the
    complement of $I$ in $\{1,\dots, d+e\}$.
    Moreover, $\phi(f)$ is surjective since $\Delta_{I^C}(\phi(f))= \det (-\id_{R^e}) = (-1)^{e}$. Let
    $\psi : D(\Delta_I) \times \GL_e \rightarrow D(\Delta_I) \times_\Sch{Inj} \Sch{Ex}(d,e)$ be the
    morphism given by $\psi(f,\sigma) = (f, \sigma \phi(f))$. This is an isomorphism. Therefore, the
    claim follows.
\end{proof}
\begin{defn}
    Let $d, e \in \N$. We define the scheme $\Gr{d}{d+e}$ via its functor of points
    \[ \Gr{d}{d+e} (R) := \Set{ P \subset R^{d+e} | P \text{ is a direct summand of }
    \rank d}. \]
    This is a projective scheme, the Grassmannian, smooth over $\Z$.
\end{defn}
\begin{proof}
    See \cite{DG}.
\end{proof}
More generally, we can define the flag scheme.
\begin{defn}
    Let $\seqv{d} = (d^0,\dots,d^\nu)$
    be a sequence of integers.
    Define
    $\OFl(\seqv{d}) \subset \prod \Gr{d^i}{d^\nu}$ as
    \[ \OFl(\seqv{d}) (R) := \Set{ (U^0, \dots, U^\nu) \in \prod_{i=0}^\nu \Gr{d^i}{d^\nu}(R) |
    U^i \subset U^{i+1} \text{ for all } 0 \le i \le \nu-1}. \]
    This is a closed subscheme of the product of Grassmannians and therefore projective.
\end{defn}
\begin{theom}
    Let $d, e \in \N$.
    The projection $\Sch{Inj}(d,d+e) \rightarrow \Gr{d}{d+e}$ sending each
    $f$ to $\Bild(f)$ is a principal $\GL_d$-bundle.
\end{theom}
\begin{proof}
    See, for example, \cite{bill_homsgenrep}.
\end{proof}
Let $Q$ be a quiver, $\dvec{d}$ and $\dvec{e}$ dimension vectors and
$\flvec{d} = (\dvec{d}^0, \dots, \dvec{d}^\nu)$ a filtration. By taking
fibre products we can generalise the previous constructions and define
$\Sch{Hom}(\dvec{d}, \dvec{e})$,
$\Sch{Inj}(\dvec{d},\dvec{d}+\dvec{e})$,
$\Sch{Surj}(\dvec{d}+\dvec{e},\dvec{e})$,
$\Sch{Ex}(\dvec{d},\dvec{e})$,
$\Gr{\dvec{d}}{\dvec{d}+\dvec{e}}$
and $\OFl(\flvec{d})$ pointwise.  For all above schemes we can do base change
to any commutative ring $K$, and we will denote these schemes
by e.g. $\Sch{Hom}(\dvec{d}, \dvec{e})_K$.
\begin{defn}
    Define the representation scheme
    \[
    \Rep[Q]{\dvec{d}} := \prod_{\alpha:i \rightarrow j} \Sch{Hom}(d_i,d_j).
    \]
    This is isomorphic to an affine space.
\end{defn}
\begin{remark}
   For each scheme defined for a quiver $Q$ we often omit the index if $Q$ is obvious, e.g.
   $\Rep{\dvec{d}} = \Rep[Q]{\dvec{d}}$.
\end{remark}
More generally we have the module scheme.
\begin{defn}
    Let $K$ be a field, $\Lambda$ a finitely generated $K$-algebra and
    \[\rho \colon K\langle x_1, \dots, x_m\rangle \rightarrow \Lambda\]
    a surjective map from the free associative algebra to $\Lambda$. The affine $K$-scheme $\MMod_\Lambda(d)$ is defined by
    \[ \MMod_\Lambda(d) (R) =
    \Set{ (M^1, \dots, M^m) \in (\End(R^d))^m | f(M^1, \dots, M^m) = 0 \  \forall\  f \in \Ker \rho}.\]
    There is a natural $\GL_d$-action on $\MMod_\Lambda(d)$ given by conjugation.
\end{defn}
\begin{remark}
    We have that $\MMod_\Lambda(d)$ is isomorphic to the functor
    \[
    R \mapsto \Hom_{R-alg}(\Lambda \tensor R, \End(R^d)).
    \]
    Therefore, another choice of $\rho$ gives an isomorphic scheme.
\end{remark}

Assume that $\Lambda$ is finite dimensional and $\Lambda = S \oplus J$
as $K$-module with $J$ being the Jacobson radical and
$S \cong K^n$ being a semisimple subalgebra. Let $e_i$, $1 \le i \le n$, be the canonical basis
of $S$ and $f_1, \dots, f_m$ elements in the union of the $e_i J e_j$, $1 \le i,j \le n$, such that
the residue classes form a basis of the direct sum of the $e_i J/J^2 e_j$. Then the $e_i$ together
with the $f_i$ generate $\Lambda$. We can choose a presentation
$\rho \colon K\langle x_1, \dots, x_n, y_1, \dots, y_m\rangle \rightarrow \Lambda$ sending
$x_i$ to $e_i$ and $y_i$ to $f_i$. Then a point $M$ in $\MMod_\Lambda(d)(R)$ starts with some
matrices $M^1, \dots, M^n$ corresponding to the $e_i$ and these define a point
in $\MMod_S(d)(R)$ so that we have a functorial $\GL_d$-equivariant morphism $p$ from
$\MMod_{\Lambda}(d)$ to $\MMod_S(d)$. In this case, one has the following lemma,
see \cite[lemma 1]{Bongartz_morita} and \cite[1.4]{Gabriel_finiteisopen}. Note
that both authors do the proof only for $K$ algebraically closed, but if one has a
decomposition of the algebra as above, the same proof works.
\begin{lmm}
  For any vector $\dvec{d}\in \N^n$ such that $\sum_{i=1}^n d_i = d$, there is a connected
  component $\MMod_S(\dvec{d})$ of $\MMod_S(d)$ characterised by the fact that
  $\MMod_S(\dvec{d})(K)$ consists of the $\GL_d(K)$-orbit of the semisimple module containing
  $Se_i$ with multiplicity $d_i$. Every connected component of $\MMod_S(d)$ is of this type and smooth. The
  connected components $\MMod_\Lambda(\dvec{d})$ of $\MMod_\Lambda(d)$ are the inverse
  images of the $\MMod_S(\dvec{d})$ under $p$.
\end{lmm}

Let $Q$ be a quiver, $K$ a field, $d$ an integer and $\dvec{d}$ a dimension vector such
that $\sum d_i = d$. Then, after choosing an isomorphism
$K^d \rightarrow \bigoplus K^{d_i}$, there is a natural immersion
\[
\Rep[Q]{\dvec{d}} \rightarrow \MMod_{KQ}(\dvec{d}) \subset \MMod_{KQ}(d).
\]

We define now quiver flags and quiver Grassmannians as varieties. Choose a $K$-representation
$M \in \Rep[Q]{\dvec{d} + \dvec{e}}(K)$.
We set
\[ \Gr[Q]{\dvec{d}}{M}(R) := \Set{ U \in \Gr{\dvec{d}}{\dvec{d} + \dvec{e}}(R) | U \text{ is a subrepesentation of }
M \tensor_K R}
\]
for each $K$-algebra $R$. $\Gr[Q]{\dvec{d}}{M}$ is a closed subscheme of $\Gr{\dvec{d}}{\dvec{d}+\dvec{e}}_K$ and
therefore projective. In a similar way we obtain a closed subscheme
$\Fl[Q]{\flvec{d}}{M}$ of $\OFl(\flvec{d})_K$ for a filtration $\flvec{d}$ of $\dimve M$.

Let $\Lambda$ be a $K$-algebra. Let $d$ and $e$ be two integers.
For a $\Lambda$-module $M \in \MMod_\Lambda(d+e)$
we define
\[ \Gr[\Lambda]{d}{M}(R) := \Set{ U \in \Gr{d}{d + e}(R) | U \text{ is a submodule of }
M \tensor_K R}.
\]
If $\Lambda \cong S \oplus J$ as before, then we have that
\[ \Gr[\Lambda]{d}{M} = \coprod \Gr[\Lambda]{\dvec{d}}{M},\]
each $\Gr[\Lambda]{\dvec{d}}{M}$ being open in \Gr[\Lambda]{d}{M}.
Moreover, for $\Lambda = KQ$ we obtain that
\[ \Gr[Q]{\dvec{d}}{M} \cong \Gr[KQ]{\dvec{d}}{M} \]
via the immersion $\Rep[Q]{\dvec{d}} \rightarrow \MMod_{KQ}(\dvec{d})$.

We have the following well-known result.
\begin{propn}
    Let $X, Y$ be two schemes, $Y$ being irreducible.
    Let $f \colon X \rightarrow Y$ be a morphism of schemes such that $f$ is open and for each $z \in Y$ the
    fibre $f^{-1}(z)$ is irreducible. Then $X$ is irreducible.
    \label{isch:propn:openirred}
\end{propn}
\begin{proof}
    Take $\emptyset \neq U, V \subset X$ and $U, V$ open. We need to show that $U \cap V$ is non-empty. We know that
    $f(U)$ and $f(V)$ are non-empty and open in $Y$. Therefore, they intersect non-trivially. Let
    $y \in f(U) \cap f(V)$.
    By definition, we have
    that $U \cap f^{-1}(y) \neq \emptyset$ and $V \cap f^{-1}(y) \neq \emptyset$ and both are open in $f^{-1}(y)$. By assumption,
    the fibre $f^{-1}(y)$ is irreducible and therefore $U \cap V \cap f^{-1}(y) \neq \emptyset$ and we are done.
\end{proof}
\begin{remark}
    The same is true for arbitrary topological spaces, since the proof relies purely on topology.
\end{remark}

\section{Grassmannians and Tangent Spaces}
Let $K$ be a field and $\Lambda$ a finitely generated $K$-algebra. We want to calculate the tangent
space at a point of $\Gr[\Lambda]{d}{M}$.
\begin{lmm}
    Let $d, e \in \N$, $K$ a field, $\Lambda$ a finitely generated $K$-algebra and $M \in \MMod_\Lambda(d+e)(K)$.
  Let $U \in \Gr[\Lambda]{d}{M}(F)$, for a field extension $F$ of $K$. Then
  \[ T_U \Gr[\Lambda]{d}{M} \cong \Hom_{\Lambda\tensor_K F} (U, (M\tensor_K F)/U). \]
\end{lmm}
\begin{proof}
  In this proof we use left modules since the notation is more convenient.
  For the tangent space we use the definition of \cite[I, \S 4, no 4]{DG}.
  By base change, we can assume that $F=K$. We prove the claim by doing a $K[\varepsilon]$ valued point calculation ($\varepsilon^2 = 0$).
  The short exact sequence
    \[
    \xymatrix{ 0 \ar[r]& U \ar[r]^\iota & M  \ar[r]^\pi \ar@/^/@{.>}[l]^p
    & M/U \ar[r] \ar@/^/@{.>}[l]^j & 0}
    \]
  is split in the category of $K$-vector spaces, therefore there is a retraction $p$ of $\iota$ and a
  section $j$ of $\pi$. We consider elements of $U$ as elements of $M$ via the inclusion $\iota$.

  The map $K[\varepsilon] \rightarrow K$ given by $s + r \epsilon \mapsto s$ induces
  a map $\theta\colon V \tensor K[\varepsilon] \rightarrow V$ for each $K$-vector space $V$.

  For each homomorphism $f \in \Hom_\Lambda(U, M/U)$ we define
  \[
  S_f := \Set{ u + v \varepsilon | u \in U,\ \pi(v)=f(u) } \subset M \tensor K[\varepsilon].
  \]
  Note that $\theta(S_f) = U$.
  We need to show that $S_f \in \Gr[\Lambda \tensor {K[\varepsilon]}]{d}{M \tensor K[\varepsilon]}$
  and that every element $S$ of the Grassmannian with $\theta(S) = U$ arises in this way.

  First, we show that $S_f$ is a $\Lambda \tensor K[\varepsilon]$-submodule. Let $u + v\varepsilon \in S_f$
  and $r +s\varepsilon \in \Lambda \tensor K[\varepsilon]$. Then we have
  \[ (r+s\varepsilon)(u+v\varepsilon)= ru + (rv + su)\varepsilon. \]
  Since $\pi(v) = f(u)$, $u \in U$ and $f$ is a $\Lambda$-homomorphism we obtain that
  \[ \pi(rv + su) = r\pi(v) = r f(u) = f(ru). \]
  Therefore, $S_f$ is a $\Lambda \tensor K[\varepsilon]$ submodule.

  Now we show that $S_f$ is, as a $K[\varepsilon]$-module, a summand of $M \tensor K[\varepsilon]$ of
  rank $d$. Let $\tilde f\in \Hom_K(U, M):=j \circ f$
  be a $K$-linear lift of $f$.Let
  \begin{alignat*}{2}
      \phi &\colon\quad& U \tensor K[\varepsilon] & \rightarrow M \tensor K[\varepsilon]\\
      && u + v\varepsilon & \mapsto u + (\tilde f (u) +v)\varepsilon.
  \intertext{Obviously, $\phi$ is $K[\varepsilon]$-linear and $\Bild \phi = S_f$. Moreover,
  $\phi$ is split with retraction}
      \psi &\colon\quad& M \tensor K[\varepsilon] & \rightarrow U \tensor K[\varepsilon]\\
      && x + y\varepsilon & \mapsto p(x) + p(y-\tilde f\circ p(x))\varepsilon.
  \end{alignat*}
  Therefore, $S_f$ is a summand of $M \tensor K[\varepsilon]$ of rank $d$.
  
  On the other hand, let $S \in \Gr[\Lambda \tensor {K[\varepsilon]}]{d}{M \tensor K[\varepsilon]}$
  such that $\theta(S)=U$. Then, $U\varepsilon$ is a $K$-subspace of $S$ with
  $\dim_K U\varepsilon = d$, therefore $\dim_K S/(U\varepsilon) = d$. The map sending
  $u + v\varepsilon \in S$ to $u \in U$ is surjective, therefore the induced map from
  $S/(U\varepsilon)$ is an isomorphism.
  Hence, for each $u \in U$ there is a $v \in M$,
  such that $u + v\varepsilon \in S$ and $v$ is unique up to a element in $U\varepsilon$.
  Set $f(u) := \pi(v)$ and, by the discussion before, this does not depend on the choice of $v$
  and we have that $S = S_f$. Moreover,
  $f \in \Hom_{\Lambda}(U,M/U)$: Let $r \in \Lambda$, $u \in U$ and $v \in M$, such that $u+v\varepsilon \in S$.
  Then, $ru + rv\varepsilon \in S$. By definition, $f(ru) = \pi(rv) = r\pi(v) = rf(u)$.
\end{proof}

%
%
\section{Geometry of Quiver Flags}
Now we come back to quiver flags. Let $K$ be a field.
We can consider $\Rep[Q]{\dvec{d}}$ as an affine scheme over $K$ with the obvious functor of points. More
generally, we work in the category of schemes over $K$.
Fix a filtration
$\seqv{\dvec{d}}=(\dvec{d}^0 = 0\le \dvec{d}^1 \le \dots \le \dvec{d}^\nu)$.

Let $\Lambda := (KQ) A_{\nu+1}$. Then $\Mod \Lambda$ is the category of sequences of $K$-representations of
$Q$ of length $\nu+1$ and
chain maps between them, i.e. a morphism between two modules
\[
\boldsymbol{M}= M^0 \rightarrow M^1 \rightarrow \dots \rightarrow M^\nu
\]
and
\[
\boldsymbol{N}= N^0 \rightarrow N^1 \rightarrow \dots \rightarrow N^\nu
\]
is given by a commuting diagram
\[
\begin{CD}
    M^0 @>>> M^1 @>>> \cdots @>>> M^\nu\\
     @VVV   @VVV    @.          @VVV\\
    N^0 @>>> N^1 @>>> \cdots @>>> N^\nu.
\end{CD}
\]
For $\Lambda$ it is easy to calculate the Euler form of two modules $\seqv{M}, \seqv{N} \in \Mod \Lambda$ by using
theorem $\ref{euf_quiverquiver}$:
\[ \euf{\seqv{M},\seqv{N}}_\Lambda = \sum_{i=0}^r \euf{M^i,N^i}_{KQ} - \sum_{i=0}^{r-1} \euf{M^i, N^{i+1}}_{KQ}. \]

We show that $\Gr[\Lambda]{\seqv{\dvec{d}}}{\seqv{M}} \cong \Fl[Q]{\seqv{\dvec{d}}}{M}$, where
$\seqv{M} = (M=M=\dots=M)$, and then use the previous
results to calculate the tangent space.
\begin{lmm}
    Let $\flvec{d}$ be a filtration and
    $M \in \Rep[Q]{\dvec{d}^\nu}(K)$.
    Let $\seqv{U} = (U^0, U^1, \dots, U^\nu) \in \Fl[Q]{\seqv{\dvec{d}}}{M}(F)$ for a field extension
    $F$ of $K$. Then we have that
    \[
T_\seqv{U} \Fl[Q]{\seqv{\dvec{d}}}{M} \cong \Hom_{\Lambda\tensor F}(\seqv{U}, (\seqv{M}\tensor F)/\seqv{U}),
    \]
    where $\Lambda = (KQ)A_{\nu + 1}$ and $\seqv{M} = (M=M=\dots=M)
    \in \MMod_\Lambda (\dvec{d}^\nu,\dvec{d}^\nu, \dots, \dvec{d}^\nu)(K)$.
\end{lmm}
\begin{proof}
    For a submodule
    \[
    \seqv{U}=(U^0 \rightarrow U^1 \rightarrow\dots\rightarrow U^\nu) \in \Gr[\Lambda]{\seqv{\dvec{d}}}{\seqv{M}}(R)
    \]
    we have automatically that the maps $U^i \rightarrow U^{i+1}$ are injections. Therefore,
    such a submodule $\seqv{U}$ gives, in a natural way, rise to a flag
    $\seqv{U} \in \Fl[Q]{\seqv{\dvec{d}}}{M}(R)$ and vice versa. This yields an isomorphism
    $\Gr[\Lambda]{\seqv{\dvec{d}}}{\seqv{M}} \cong \Fl[Q]{\seqv{\dvec{d}}}{M}$. Since
    $\Gr[\Lambda]{\seqv{\dvec{d}}}{\seqv{M}}$ is open in
    $\Gr[\Lambda]{d}{\seqv{M}}$, where $d = \sum_{k,i} d^k_i$, we have for a point
    $\seqv{U} \in \Fl[Q]{\seqv{\dvec{d}}}{M}(F)$ that
    \[ T_\seqv{U} \Fl[Q]{\seqv{\dvec{d}}}{M} \cong T_\seqv{U} \Gr[\Lambda]{\seqv{\dvec{d}}}{\seqv{M}}
    = T_\seqv{U} \Gr[\Lambda]{d}{\seqv{M}} = \Hom_{\Lambda\tensor F}(\seqv{U}, (\seqv{M}\tensor F)/\seqv{U}). \]
\end{proof}


We define the closed subscheme $\RepFl[Q]{\seqv{\dvec{d}}}$ of
$\Rep[Q]{\dvec{d}^\nu} \times \OFl(\seqv{\dvec{d}})$ by its functor of points
\[ \RepFl[Q]{\seqv{\dvec{d}}}(R) := \Set{ (M, \seqv{U}) \in
\Rep[Q]{\dvec{d}^\nu}(R) \times \OFl(\flvec{d})(R) | \seqv{U} \in \Fl[Q]{\seqv{\dvec{d}}}{M}}. \]

We have the following.

\begin{lmm}
    Let $\flvec{d}$ be a filtration.
Consider the two natural projections from the fibre product restricted to $\RepFl[Q]{\seqv{\dvec{d}}}$.
\[
\xymatrix{
\RepFl[Q]{\seqv{\dvec{d}}} \ar[r]^{\pi_1} \ar[d]_{\pi_2}& \Rep[Q]{\dvec{d}^\nu}\\
\OFl(\flvec{d})&
}
\]
    \label{lmm:geomflags:projtriv}
Then $\pi_1$ is projective and $\pi_2$ is a vector bundle of rank
\[
\sum_{k=1}^\nu \sum_{\alpha:i \rightarrow j} d^k_j (d_i^k - d_i^{k-1}).
\]
Therefore, $\RepFl[Q]{\seqv{\dvec{d}}}$ is smooth and irreducible of dimension
\[ \sum_{k=1}^{\nu-1} \euf{\dvec{d}^k, \dvec{d}^{k+1} - \dvec{d}^k}_Q + \dim \Rep[Q]{\dvec{d}^\nu}. \]
Finally, the (scheme-theoretic) image
$\mathcal{A}_\flvec{d}:=\pi_1(\RepFl[Q]{\seqv{\dvec{d}}})$ is a closed, irreducible
subvariety of $\Rep[Q]{\dvec{d}^\nu}$.
\end{lmm}
\newcommand{\myunit}{1.5ex}
\tikzset{%
    node style ge/.style={rectangle,minimum size=\myunit}
}
\begin{proof}
    $\pi_1$ is projective since it factors as a closed immersion into projective space times
    $\ORep_{Q}$ followed by the projection to $\ORep_{Q}$.

    For $\dvec{I}=(I_i)_{i \in Q_0}$, each $I_i \subset \{1,\dots,d^\nu_i\}$,
    we set $W_\dvec{I}$ to be the
    graded subspace
    of $K^{\dvec{d}^\nu}$ with basis $\{e_j\}_{j \in I_i}$ in the $i$-th graded part $K^{d^\nu_i}$ and
    $|\dvec{I}| := (|I_i|)_{i \in Q_0} \in \N^{Q_0}$.
    For a sequence $\seqv{\dvec{I}} = (\dvec{I}^0, \dvec{I}^1, \dots, \dvec{I}^\nu)$ such that
    $I_i^k \supset I_i^{k+1}$ and $|\dvec{I}^k|= \dvec{d}^\nu - \dvec{d}^k$ we set
    \[
    \seqv{W}_{\seqv{\dvec{I}}} := (W_{\dvec{I}^0}, \dots, W_{\dvec{I}^\nu})
    \]
    to be the decreasing sequence of
    subspaces associated to $\seqv{\dvec{I}}$.
    We show that $\pi_2$ is trivial over the open affine subset $U_{\seqv{\dvec{I}}}$ of
    $\OFl(\seqv{\dvec{d}})$
    given by
    \[
    U_{\seqv{\dvec{I}}} (R) := \Set{ \seqv{U} \in \OFl(\flvec{d})(R) |
    U^k \oplus (W_{\dvec{I}^k} \tensor R) = R^{\dvec{d}^\nu}}.
    \]
    Without loss
    of generality we assume $I^k_i = \{d^k_i+1, \dots, d^\nu_i\}$. Each element
    $\seqv{U} \in U_\seqv{\dvec{I}}(R)$ is given uniquely by some matrices
    $A_i^k \in \Mat_{(d^\nu_i-d^k_i)\times (d^k_i-d^k_{i-1})}$
    such that
 \[U^k_i = \Bild
\begin{array}{c@{}}\begin{tikzpicture}
\matrix (A) [matrix of math nodes,%
             nodes = {node style ge},%
             left delimiter  = (,%
             right delimiter = )]
             {%
         \id_{d^1_i} & 0 & \cdots & 0\\
        \ & \id_{d^2_i-d^1_i} & \ddots & \vdots\\
        \ & \hphantom{AB} & \ddots &0\\
        \ & A^2_i & \ddots & \id_{d^{k}_i-d^{k-1}_i}\\
        \hphantom{AB} & \hphantom{AB} & \cdots & A^k_i\\
         };
         \path ($ (A-2-1)!0.5!(A-4-1) $) node {$A^1_i$};
         \draw (A-1-1.south west) rectangle (A-5-1.south east);
         \draw ($ (A-2-2.south west)!0.3!(A-2-2.south) $) rectangle (A-5-2.south east);
         \draw (A-5-4.north west) rectangle (A-5-4.south east);
     \end{tikzpicture}\end{array}.\]

      Let $V^k:=W_{(\{1,\dots,d^k_i\})}$. Let $X$ be the closed subscheme of
      $\Rep[Q]{\dvec{d}^\nu}$
      given by the functor of points
      \[X(R) := \Set{ M \in \Rep[Q]{\dvec{d}^\nu} | V^k \tensor R \text{ is a subrepresentation of } M 
      \; \forall \; 0 \le k \le \nu}.\]
      Note that $X$ is an affine space of dimension
      \[\sum_{k=1}^\nu \sum_{\alpha:i \rightarrow j} d^k_j (d_i^k - d_i^{k-1}).\]
      Let $g_\seqv{U}:= (g_i)_{i\in Q_0}$ where
    \[g_i =
\begin{array}{c@{}}
\begin{tikzpicture}
\matrix (A) [matrix of math nodes,%
             nodes = {node style ge},%
             left delimiter  = (,%
             right delimiter = )]
             {%
      \id_{d^1_i} & 0 & \cdots & 0 & 0\\
      \hphantom{AB} & \id_{d^2_i-d^1_i} & \ddots & \vdots & 0\\
       & \hphantom{AB}& \ddots & 0 & \vdots\\
      \hphantom{AB} & A^2_i & \ddots & \id_{d^{\nu-1}_i-d^{\nu-2}_i} & 0\\
      \hphantom{AB} &\hphantom{AB}  & \cdots & A^{\nu-1}_i & \id_{d^{\nu}_i-d^{\nu-1}_i}\\
      };
         \path ($ (A-2-1)!0.5!(A-4-1) $) node {$A^1_i$};
         \draw (A-1-1.south west) rectangle (A-5-1.south east);
         \draw ($ (A-2-2.south west)!0.3!(A-2-2.south) $) rectangle (A-5-2.south east);
         \draw (A-5-4.north west) rectangle (A-5-4.south east);
  \end{tikzpicture}
  \end{array}
      \in \GL_{d^\nu_i}(K).\]
      Then, the map from $X \times U_{\seqv{\dvec{I}}}$ to
      $U_{\seqv{\dvec{I}}} \times_{\OFl} \RepFl[Q]{\flvec{d}}$ given by sending
      $(M, \seqv{U})$ to $(g_\seqv{U} \cdot M, \seqv{U})$ is an isomorphism which induces
      an isomorphism of vector spaces on the fibres.
      Therefore, we have that $\pi_2$ is a vector bundle.

      Finally, we prove the claim on dimension. Since $\OFl(\flvec{d})$ is smooth,
      we have that
      \begin{multline*}
          \dim \OFl(\flvec{d}) =
          \sum_{i \in Q_0} \sum_{k=1}^{\nu-1} \sum_{l=k+1}^\nu (d^k_i - d^{k-1}_{i}) (d^l_i - d^{l-1}_{i})\\
          =\sum_{i \in Q_0} \left( \sum_{k=1}^{\nu-1} \sum_{l=k+1}^\nu d^k_i (d^l_i - d^{l-1}_{i}) -
          \sum_{k=1}^{\nu-2} \sum_{l=k+2}^\nu d^{k}_{i} (d^l_i - d^{l-1}_{i}) \right)\\
          =\sum_{i \in Q_0} \left( d^{\nu-1} (d^\nu_i - d^{\nu-1}_i) + \sum_{k=1}^{\nu-2} d^k_i (d^{k+1}_i - d^k_i) \right) =
          \sum_{i \in Q_0} \left(\sum_{k=1}^{\nu-1} d^k_i (d^{k+1}_i - d^k_i) \right).
      \end{multline*}
      Since
      $\RepFl[Q]{\flvec{d}}$ is smooth and $\pi_2$ is a vector bundle we obtain
      \begin{multline*}
        \dim \RepFl[Q]{\flvec{d}} = \dim \OFl(\flvec{d}) +
        \sum_{\alpha \colon i \rightarrow j}\sum_{k=1}^\nu (d^k_i- d^{k-1}_{i}) d^{k}_j\\
          =\sum_{i \in Q_0} \left(\sum_{k=1}^{\nu-1} d^k_i (d^{k+1}_i - d^k_i) \right) +
          \sum_{\alpha \colon i \rightarrow j}\sum_{k=1}^\nu (d^k_i- d^{k-1}_{i}) d^{k}_j\\
      \end{multline*}
      \begin{multline*}
          =\sum_{k=1}^{\nu-1} \left( \sum_{i \in Q_0}d^k_i (d^{k+1}_i - d^k_i)  +
          \sum_{\alpha \colon i \rightarrow j} d^k_i (d^{k}_j-d^{k+1}_j)\right)
          + \sum_{\alpha \colon i \rightarrow j} d^\nu_i d^\nu_j\\
      =\sum_{k=1}^{\nu -1} \euf{\dvec{d}^k, \dvec{d}^{k+1} - \dvec{d}^k}_Q + \dim \Rep[Q]{\dvec{d}^\nu}.
  \end{multline*}
\end{proof}
\begin{remark}
    Note that if $M$ is a $K$-valued point of $\mathcal{A}_\flvec{d}$ for $K$ not algebraically closed,
    then $M$ does not necessarily have a flag of type $\flvec{d}$. This only becomes true after a finite
    field extension.
\end{remark}
We now can give an estimate for the codimension of $\mathcal{A}_\flvec{d}$ in $\Rep[Q]{\dvec{d}^\nu}$.
For this we use Chevalley's theorem.
\begin{theom}[Chevalley]
    Let $X, Y$ be irreducible schemes over a field $K$ and let $f \colon X \rightarrow Y$ be a dominant morphism.
    Then for every point $y \in Y$ and every point $x \in f^{-1}(y)$, the scheme theoretic fibre, we
    have that
    \[ \dim_x f^{-1}(y) \ge \dim X - \dim Y. \]
    Moreover, on an open, non-empty subset of $X$ we have equality.
\end{theom}
\begin{proof}
    See \cite[\S 5, Proposition 5.6.5]{EGA4}.
\end{proof}
\begin{theom}
  Let $\flvec{d}$ be a filtration, $F$ a field extension of $K$, $\Lambda := (FQ)A_{\nu+1}$ and
  $(M,\seqv{U}) \in \RepFl[Q]{\flvec{d}}(F)$. Let
  $\seqv{M}=(M=\dots=M)$ as a $\Lambda$-module. Then we have that
  \[ \codim \mathcal{A}_\flvec{d} \le \dim \Ext^1_{\Lambda}(\seqv{U}, \seqv{M}/\seqv{U})
  \le \dim \Ext^1_{\Lambda}(\seqv{U}, \seqv{M})
  \le \dim \Ext^1_{FQ}(M, M). \]
  \label{them:codimflag}
\end{theom}
\begin{proof}
  Since $\dim \mathcal{A}_{\flvec{d}}$ is stable under flat base change we can
  assume $K=F$.
  Let $\seqv{V} := \seqv{M}/\seqv{U}$.
  Then we have the following short exact sequence of $\Lambda$-modules:
  \[
  \begin{CD}
    0:@.\quad           @.  0   @>>>  0  @>>> \cdots @>>> 0\\
    @VVV @.  @VVV      @VVV         @.       @VVV\\
    \seqv{U}:@.\quad   @.  U^0 @>>> U^1 @>>> \cdots @>>> U^\nu\\
    @VVV @.  @VVV      @VVV         @.       @VVV\\
    \seqv{M}:@.\quad   @.   M @=     M  @=   \cdots @=    M\\
    @VVV @.  @VVV      @VVV         @.       @VVV\\
    \seqv{V}:@.\quad   @.  V^0 @>>> V^1 @>>> \cdots @>>> V^\nu\\
    @VVV @.  @VVV      @VVV         @.       @VVV\\
    0:@.\quad         @.    0   @>>>  0  @>>> \cdots @>>> 0.\\
  \end{CD}
  \]
  We already know that $\Hom_\Lambda(\seqv{U}, \seqv{V})$
  is the tangent space of $\Fl{\flvec{d}}{M}$ at the point $\seqv{U}$.
  Using Chevalley's theorem, we have that
  \[\dim \Hom_\Lambda(\seqv{U}, \seqv{V}) \ge \dim_\seqv{U}\Fl{\flvec{d}}{M} \ge
  \dim \RepFl[Q]{\flvec{d}} - \dim \mathcal{A}_\flvec{d} \]
  and therefore
  \[ \dim \mathcal{A}_\flvec{d} \ge \dim \RepFl[Q]{\flvec{d}} - \dim \Hom_\Lambda(\seqv{U}, \seqv{V}).\]
  We now calculate
  \begin{multline*}
    \euf{\seqv{U}, \seqv{V}}_\Lambda = \sum_{k=0}^\nu \euf{U^k, V^k}_{Q} - \sum_{k=0}^{\nu-1} \euf{U^k, V^{k+1}}_{Q}\\
    = \sum_{k=1}^{\nu-1} \euf{\dvec{d}^k, \dvec{d}^\nu - \dvec{d}^k}_{Q} - \sum_{k=1}^{\nu-1}
    \euf{\dvec{d}^k, \dvec{d}^\nu - \dvec{d}^{k+1}}_{Q}
    = \sum_{k=1}^{\nu-1} \euf{\dvec{d}^k, \dvec{d}^{k+1} - \dvec{d}^k}_{Q}.
  \end{multline*}
  Recall that
\[ \dim \RepFl[Q]{\seqv{\dvec{d}}} = \sum_{k=1}^{\nu-1}
\euf{\dvec{d}^k, \dvec{d}^{k+1} - \dvec{d}^k}_Q + \dim \Rep[Q]{\dvec{d}^\nu}. \]
  In total
  \begin{multline*}
    \codim \mathcal{A}_\flvec{d} \le \dim \Rep{\dvec{d}} + \dim \Hom_\Lambda (\seqv{U}, \seqv{V})
    - \dim \RepFl[Q]{\flvec{d}}\\
    = \dim \Hom_\Lambda (\seqv{U}, \seqv{V}) - \sum_{k=1}^{\nu-1} \euf{\dvec{d}^k, \dvec{d}^{k+1} - \dvec{d}^k}_{Q}
    =\dim \Hom_\Lambda (\seqv{U}, \seqv{V}) - \bform{\seqv{U}}{\seqv{V}}_\Lambda\\
    = \dim \Ext^1_\Lambda(\seqv{U}, \seqv{V}) - \dim \Ext^2_\Lambda(\seqv{U}, \seqv{V}).
  \end{multline*}
  Here we have the last equality since $\gldim \Lambda \le 2$.
  Since $\gldim KQ = 1$ and $P = P = \dots = P$ is projective in $\Mod \Lambda$
  for every projective $P$ in $\Mod KQ$ we see that
  $\pd_\Lambda \seqv{M} \le 1$ and similarly $\id_\Lambda \seqv{M} \le 1$. Consider, as before, the short exact sequence
  \[ 0 \rightarrow \seqv{U} \rightarrow \seqv{M} \rightarrow \seqv{V} \rightarrow 0. \]
  Applying $(-, \seqv{V})$ yields $(\seqv{U}, \seqv{V})^2 = 0$. Applying $(-,\seqv{M})$ gives
  a surjection $(\seqv{M},\seqv{M})^1 \rightarrow (\seqv{U}, \seqv{M})^1$. Applying $(-,\seqv{U})$
  yields that $(\seqv{U}, \seqv{U})^2 = 0$. Therefore,
  applying $(\seqv{U},-)$ yields a surjection $(\seqv{U}, \seqv{M})^1 \rightarrow (\seqv{U}, \seqv{V})^1$.
  Hence the above result simplifies to
  \[ \codim \mathcal{A}_\flvec{a} \le \dim \Ext^1(\seqv{U}, \seqv{V}) \le \dim \Ext^1(\seqv{U}, \seqv{M})
  \le \dim \Ext^1(\seqv{M}, \seqv{M}). \]

  Obviously, $\Ext^1_\Lambda(\seqv{M}, \seqv{M}) \cong \Ext^1_{KQ}(M,M)$ and the claim follows.
\end{proof}
\begin{remark}
    Note that if the characteristic of $K$ is $0$, then, by generic smoothness,
    there is a point $M \in \mathcal{A}_\flvec{d}$ and an $\seqv{U} \in \Fl{\flvec{d}}{M}$ such
    that $\Fl{\flvec{d}}{M}$ is smooth in $\seqv{U}$. In this case we have that
  \[ \codim \mathcal{A}_\flvec{d} = \dim \Ext^1_{\Lambda \tensor F}(\seqv{U}, (M\tensor F)/\seqv{U}).\]
\end{remark}

  We also construct an additional vector bundle.
  \begin{defn}
      Let $\flvec{d}$ be a filtration.
      Let $\Sch{Rep}_{Q,A_{\nu+1}}(\flvec{d})$ be the scheme given via its functor of points
      \begin{multline*}
          \Sch{Rep}_{Q,A_{\nu+1}}(\flvec{d})(R) := \\
          \Set{ (\seqv{U}, \seqv{f})  \in \prod_{i=0}^\nu \Rep[Q]{\dvec{d}^i}(R) \times
          \prod_{i=0}^{\nu-1} \Hom(\dvec{d}^i, \dvec{d}^{i+1})(R) | f^i \in \Hom_{RQ}(U^i, U^{i+1})}.
      \end{multline*}
      Let $\Sch{IRep}_{Q,A_{\nu+1}}(\flvec{d})$ be the open subscheme of
      $\Sch{Rep}_{Q,A_{\nu+1}}(\flvec{d})$ given by its functor of points
      \[
      \Sch{IRep}_{Q,A_{\nu+1}}(\flvec{d})(R) := \Set{ (\seqv{U}, \seqv{f}) \in \Sch{Rep}_{Q,A_{\nu+1}}(\flvec{d})(R) |
      f^i \in \Sch{Inj}(\dvec{d}^i, \dvec{d}^{i+1})(R)}.
      \]
  \end{defn}
  \begin{remark}
      Note that $\Sch{Rep}_{Q, A_{\nu +1}}(\flvec{d})(K)$ consists of sequences
      of representations of $Q$. Therefore, these are modules over $(KQ)A_{\nu+1}$. Vice
      versa, every $(KQ)A_{\nu+1}$-module of dimension vector $\flvec{d}$ is isomorphic to an element of
      $\Sch{Rep}_{Q, A_{\nu +1}}(\flvec{d})(K)$.

      For shortness, we will often write $\seqv{U}$ instead of $(\seqv{U}, \seqv{f})$ for
      an $(\seqv{U}, \seqv{f}) \in \Sch{Rep}_{Q, A_{\nu +1}}(\flvec{d})(R)$.
  \end{remark}
  \begin{lmm}
      Let $\flvec{d}$ be a filtration.
      Then the projection
      \begin{alignat*}{2}
          \pi&\colon\quad& \Sch{IRep}_{Q,A_{\nu+1}}(\flvec{d}) &\rightarrow
          \prod_{i=0}^{\nu-1} \Sch{Inj}(\dvec{d}^i, \dvec{d}^{i+1})\\
          \intertext{given by sending}
      &&(\seqv{U}, \seqv{f}) &\mapsto \seqv{f}
      \end{alignat*}
      is a vector bundle and therefore flat.
      In particular, $\Sch{IRep}_{Q,A_{\nu+1}}(\flvec{d})$ is smooth and irreducible.
  \end{lmm}
  \begin{proof}
      The first part is analogously to lemma \ref{lmm:geomflags:projtriv}. Irreducibility
      then follows by the fact that flat morphisms are open and proposition \ref{isch:propn:openirred}.
  \end{proof}
  Before we continue, we give the following easy lemma, stated by K. Bongartz
  in \cite{Bongartz_singularities}, which gives
  rise to a whole class of vector bundles.
  \begin{lmm}
      \label{lmm:geomflags:vecbun}
      Let $X$ be a variety over a ground ring $K$. Let $m, n \in \N$ and
      $f \colon X \rightarrow \Hom(m,n)_K$ a morphism. Then
      for any $r\in \N$, the variety $X(r)$ given by the functor of points
      \[ X(r)(R):= \Set{ x \in X(R) | f(x) \in\Hom(m,n)_{m-r}(R)} \]
      is a locally closed subvariety of $X$. Moreover, the closed subvariety
      \[
      U_r(R):=\Set{ (x, v) \in X(r)(R)\times R^m | f(x)(v) = 0 }
      \]
      of $X(r) \times K^m$
      is a sub
      vector bundle
      of rank $r$ over $X(r)$.
  \end{lmm}
  \begin{proof}
  The claim follows easily by using Fitting ideals.
  \end{proof}
  \begin{beispiel}
      Let $(\seqv{M},\seqv{g}) \in \Rep[Q,A_{\nu+1}]{\flvec{e}}(K)$.
      
      Let $\varphi \colon \Rep[Q,A_{\nu+1}]{\flvec{d}} \rightarrow \Hom(m,n)$
      be the morphism given by
      \[ (\seqv{U},\seqv{f}) \mapsto \left( \seqv{h}=(h_i^k) \mapsto
      \left( (h^k_j U_\alpha^k - M_\alpha^k h^k_i)_{\substack{\alpha\colon i \rightarrow j\\0 \le k \le \nu}}, 
      (h^{k+1}_i f^k_i - g^k_i h^k_i)_{\substack{i \in Q_0\\0 \le k \le \nu}} \right) \right)
      \]
      for every $(\seqv{U}, \seqv{f}) \in \Rep[Q,A_{\nu+1}]{\flvec{d}}(R)$,
      where
      \begin{align*}
          m &= \sum_{i \in Q_0} \sum_{k=0}^{\nu} d_i^k e_i^k &
          &\text{and}&
      n &= \sum_{k=0}^\nu \sum_{\alpha \colon i \rightarrow j} d_i^k e_j^k + 
      \sum_{k=0}^{\nu-1} \sum_{i \in Q_0} d_i^k e_i^{k+1}.
      \end{align*}
      Then $\seqv{h} \in \ker \varphi (\seqv{U}, \seqv{f})$ if and only if
      $\seqv{h} \in \Hom_{(R Q)A_{\nu+1}}(\seqv{U}, \seqv{M}\tensor R)$.

      Set
      \[
      \RepHom[Q,A_{\nu+1}]{\flvec{d}, \seqv{M}}_r := U_r(R)
      \]
      from the previous lemma.
      Note that elements $(\seqv{U}, \seqv{h}) \in \RepHom[Q,A_{\nu+1}]{\flvec{d}, \seqv{M}}_r(R)$
      are all pairs consisting of a representation $\seqv{U} \in \Rep[Q,A_{\nu+1}]{\flvec{d}}(R)$
      and a morphism $\seqv{h} \in \Hom_{(RQ)A_{\nu+1}}(\seqv{U}, \seqv{M})$ such that
      $\rank \Hom_{(RQ)A_{\nu+1}}(\seqv{U}, \seqv{M})=r$.

      The lemma yields that the projection
      \begin{align*}
          \RepHom[Q,A_{\nu+1}]{\flvec{d}, \seqv{M}}_r &\rightarrow \Rep[Q,A_{\nu+1}]{\flvec{d}}(r)\\
          (\seqv{U}, \seqv{h}) &\mapsto \seqv{U}
      \end{align*}
      is a vector bundle of rank $r$.

      It also stays a vector bundle if we restrict it to
      the open subset $\Sch{IRep}_{Q,A_{\nu+1}}(\flvec{d})(r)$ of $\Rep[Q,A_{\nu+1}]{\flvec{d}}(r)$.
      We denote the preimage under the projection to this variety by
      $\Sch{IRepHom}_{Q, A_{\nu+1}}(\flvec{d},\seqv{M})_r$.\label{ex:geomflags:hombundle}
%
  \end{beispiel}
  We obtain the following.
  \begin{theom}
      Let $\flvec{d}$ be a filtration and $F$ a field extension of $K$.
      \label{theom:geomflags:irred}
      Assume that there is an $M \in \mathcal{A}_\flvec{d}(F)$ such that $\dim \Ext^1_{FQ}(M,M) =
      \codim \mathcal{A}_\flvec{d}$. Then $\Fl[Q]{\flvec{d}}{M}$ is smooth over $F$ and
      geometrically irreducible.
  \end{theom}
  \begin{proof}
      Smoothness is immediate, since by the last theorem we have that
      $\dim T_{\seqv{U}}\Fl[Q]{\flvec{d}}{M}$ is constant and smaller or equal
      to the dimension at each irreducible component living in $\seqv{U}$. Since $F$
      is a field this implies smoothness. See \cite[I, \S 4, no 4]{DG}.

      Now we prove irreducibility. By base change we can assume that $F$ is algebraically closed.
      Consider all the following schemes as $F$-varieties.
      We construct the following diagram of varieties.
\[
\xymatrix{ & \Sch{IRepHom}_{Q, A_{\nu+1}}(\flvec{d},\seqv{M})_r \ar@{->>}[d]^{\text{vector bundle}} &
\Sch{IRepInj}_{Q, A_{\nu +1}}(\flvec{d},\seqv{M})_r \ar@{_{(}->}[l]_{\text{open}} \ar@{->>}[d]\\
\Sch{IRep}_{Q,A_{\nu+1}}(\flvec{d}) & \ar@{_{(}->}[l]_{\text{open}} \Sch{IRep}_{Q,A_{\nu+1}}(\flvec{d})(r)
& \Fl[Q]{\flvec{d}}{M},\\
}
\]
$r$ being equal to $\euf{\flvec{d}, \seqv{M}}_\Lambda + [M,M]^1$.
Since open subvarieties and images of irreducible varieties are again irreducible and by application
of proposition \ref{isch:propn:openirred} we then obtain that $\Fl[Q]{\flvec{d}}{M}$ also is irreducible.

  Consider the minimal value $r$ of $\dim\Hom(\seqv{U}, \seqv{M})$ for
  $\seqv{U} \in \Sch{IRep}(\flvec{d})(F)$. Denote by
  \begin{alignat*}{2}
      \pi \colon&& \Sch{IRep}(\flvec{d}) &\rightarrow \Rep[Q]{\dvec{d}^\nu}\\
                && \seqv{U} &\mapsto U^\nu.
  \end{alignat*}
  Since $\Orbit_M$ is open in $\mathcal{A}_\flvec{d}$ and $\Sch{IRep}$ is irreducible,
  the intersection of the two open sets $\pi^{-1}(\Orbit_M)$ 
  and
  $\Sch{IRep}(r)$ is non-empty. For all elements $\seqv{U}$ of $\pi^{-1}(\Orbit_M)$
  we have, by theorem \ref{them:codimflag}, that $[\seqv{U}, \seqv{M}]^1_\Lambda = [M,M]^1_Q$.
  We already saw that $[\seqv{U}, \seqv{M}]^2 = 0$, therefore
  \[
  [\seqv{U}, \seqv{M}]_\Lambda = \euf{\seqv{U}, \seqv{M}}_\Lambda + [\seqv{U}, \seqv{M}]^1_\Lambda
  = \euf{\flvec{d}, \seqv{M}}_\Lambda + [M,M]^1_Q.
  \]
  This means that the dimension of the homomorphism space is constant on $\pi^{-1}(\Orbit_M)$ and
  we obtain that $r=\euf{\flvec{d}, \seqv{M}}_\Lambda + [M,M]^1_Q$. Moreover, $\Sch{IRep}_r$
  is irreducible as an open subset of $\Sch{IRep}$.

  We then have that $\Sch{IRepHom}_{Q,A_{\nu+1}}(\flvec{d},\seqv{M})_r$
  is irreducible, since it is a vector bundle on $\Sch{IRep}_r$ by example \ref{ex:geomflags:hombundle}.
  Take the open subvariety $\Sch{IRepInj}
  (\flvec{d},\seqv{M})_r$
  of $\Sch{IRepHom}
  (\flvec{d},\seqv{M})_r$ where the morphism to $\seqv{M}$ is injective.
  It is irreducible as an open subset of an irreducible variety.
  The projection
  from this variety to $\Fl[Q]{\flvec{d}}{M}$ is surjective since
  $\pi^{-1}(\Orbit_M)$ is contained in $\Sch{IRep}(\flvec{d})(r)$, and therefore $\Fl[Q]{\flvec{d}}{M}$ is irreducible.
  \end{proof}

  We now want to interpret theorem \ref{theom:geomflags:irred} in terms of Hall numbers. Let $X_0$
  be a variety defined over a finite field $\F_q$, where $q=p^n$ for a prime $n$. Denote
  by $\overline{\F_q}$ the algebraic closure of $\F_q$ and by $X := X_0 \tensor \overline{\F_q}$
  the variety obtained from $X_0$
  by base change from $\F_q$ to $\overline{\F_q}$. Let $F$ be the Frobenius automorphism acting on $X$.
  Denote by $H^i(X, \Q_\ell)$ the $\ell$-adic cohomology group with compact support for a prime $\ell \neq p$,
  see for example \cite{Freitag_etaleweil}. Denote by $F^*$ the induced action of $F$ on
  cohomology $H^*(X, \Q_\ell)$.
  P. Deligne proved the following theorem.
\begin{theom}[P. Deligne \cite{Deligne_weil2}, 3.3.9]
Let $X_0$ be a proper and smooth variety over $\F_q$.
For every $i$, the characteristic polynomial $\det( T \id - F^*, H^i(X, \Q_\ell))$
is a polynomial with coefficients in $\Z$, independent of $\ell$ ( $\ell \neq p)$.
The complex roots $\alpha$ of this polynomial have absolute value $\lvert \alpha \rvert = q^{\frac{i}{2}}$.
\end{theom}

Moreover, the Lefschetz fixed point formula yields that
\[
\#X_0(\F_{q^n}) = \sum_{i \ge 0} (-1)^i \Tr( (F^*)^n, H^i(X, \Q_\ell)).
\]

Assume now that there is a polynomial $P \in \Q[t]$ such that, for each finite field extension
$L/\F_q$ we have that $\# X_0(L) = P(|L|)$, $|L|$ being the number of elements of the
finite field. We call $P$ the counting polynomial of $X_0$. Then we have the following.
\begin{theom}
	Let $X_0$ be proper smooth over $\F_q$ with counting polynomial $P$.
	Then odd cohomology of $X$ vanishes and
    \[ P(t) = \sum_{i=0}^{\dim X_0} \dim H^{2i}(X, \Q_\ell) t^i.
	\]
\end{theom}
\begin{proof}
	See \cite[Lemma A.1]{BillVandenBergh_absolutelyindec}.
\end{proof}

Assume now that $Y$ is a projective scheme over $\Z$ and set
$Y_K:=Y \tensor K$ for any field $K$. Note that for $Y$ $\ell$-adic
cohomology agrees with $\ell$-adic cohomology with compact support.
Assume furthermore that there is a counting polynomial $P \in \Q[t]$ such that,
for each finite field $K$, we have that $\# Y_K(K) = P(|K|)$.
By the base change theorem \cite[Theorem 1.6.1]{Freitag_etaleweil} we have
\[ H^i(Y_{\overline{\Q}}, \Q_\ell) \cong H^i (Y_\C, \Q_\ell). \]
By the comparison theorem \cite[Theorem 1.11.6]{Freitag_etaleweil} we have
\[ H^i(Y_\C, \Q_\ell) \cong H^i (Y_\C(\C), \Q_\ell),\]
where on the right hand side we consider the usual cohomology of the complex analytic manifold
attached to $Y_\C$.

Moreover, there is an open, non-empty subset $U$ of $\Spec \Z$ such that
$H^i (Y_{\overline{\kappa(v)}}, \Q_\ell) \cong H^i (Y_{\overline{\Q}}, \Q_\ell)$ for all
$v \in U$, where $\kappa(v)$ denotes the residue field at $v$. This means that for almost all primes $p$ we have that
\[H^i (Y_{\overline{\F_p}}, \Q_\ell) \cong H^i (Y_{\overline{\Q}}, \Q_\ell) \cong H^i(Y_\C(\C), \Q_\ell).\]

Therefore, if we know the Betti numbers of $Y_\C(\C)$, then we know the coefficients
of the counting polynomial. In order to apply this to our situation we use
the following theorem of W. Crawley-Boevey \cite{Bill_rigidintegral}.
\begin{theom}
    \label{theom:geomflags:rigidreps}
	Let $M$ be an $K$-representation without self-extensions. Then there is
	a $\Z$-representation $N$ such that $M= N \tensor K$ and for
	all fields $F$ we have that $\Ext (N \tensor F, N \tensor F) = 0$.
\end{theom}
Putting all this together, we obtain the following.
\begin{theom}
      Assume that there is an $M \in \mathcal{A}_\flvec{d}(\Q)$, being a direct sum of
      exceptional representations, such that
      \[
      \dim \Ext^1_{\Q Q}(M,M) =
      \codim \mathcal{A}_\flvec{d}.
      \]
      Let $N$ be a $\Z$-representation and $P \in \Q[t]$ a polynomial, such that
      $N \tensor \Q \cong M$ and
      $\# \Fl[Q]{\flvec{d}}{N \tensor \F_q} = P(q)$. Then
      $P(0) = 1$ and $P(1) = \chi\left( \Fl[Q]{\flvec{d}}{M \tensor \C}\right) > 0$.

      Moreover, if $Q$ is Dynkin or extended Dynkin, then there is a representation $N$ and
      a polynomial $P$ with the required properties.
\end{theom}
\begin{proof}
  Let $X := \Fl[Q]{\flvec{d}}{N}$ as a scheme
  over $\Z$.
  Using theorem \ref{theom:geomflags:irred} we
  obtain that $X_K$ is smooth and irreducible for every field $K$. By the previous discussion
  we have then that the
  $i$-th coefficient of $P$ is exactly
  $\dim H^{2i}(X_\C(\C), \Q_\ell)$ and that odd cohomology vanishes. Therefore,
  \[
  0<\sum\dim H^{2i}(X_\C(\C), \Q_\ell) = \chi\left( X_\C \right) = P(1).
  \]
  By irreducibility we have that
  \[
  P(0) = \dim H^{0}(X_\C(\C), \Q_\ell) = 1
  \]
  and this proves the first claims.

  If $Q$ is Dynkin or extended Dynkin,
  then let $N$ be the $\Z$-representation given by theorem \ref{theom:geomflags:rigidreps}. We have the
  polynomial $P$ since we have Hall polynomials and there is a decomposition symbol $\alpha = (\mu, \emptyset)$
  such that
  $N\tensor K$ is in $\SegreK(\alpha, K)$ for any field $K$, since it is a direct sum of exceptional
  representations and therefore discrete.
\end{proof}

\chapter{Reflections on Quiver Flags}
\dictum[Terry Pratchett]{Always be wary of any helpful item that weighs less than its operating manual.}
\label{chap:flags}
Let $K$ be an arbitrary field and $Q$ a quiver.
\section{Reflections and Quiver Flags}
Let $M$ be a $K$-representation of $Q$ and $\flvec{d}$ a filtration of $\dimve M$.
We want to define reflections on a flag $\seqv{U} \in \Fl[Q]{\seqv{\dvec{d}}}{M}$. Let $a$ be
a sink and let $\seqv{U} \in \Fl[Q]{\seqv{\dvec{d}}}{M}$. Then, for each $i$, we have the following commutative diagram with exact rows.
\[
\begin{CD}
    0 @>>>  (S^+_a U^{i-1})_a      @>>> \bigoplus U^{i-1}_j @>{\phi^{U^{i-1}}_a}>> \Bild \phi^{U^{i-1}}_a @>>> 0\\
    @.          @VVfV                        @VVV           @VVgV     @.\\
  0 @>>> (S^+_a U^i)_a @>>> \bigoplus U^i_j @>>\phi^{U^i}_a>
  \Bild \phi^{U^i}_a @>>> 0\\
\end{CD}
\]
By definition, $g$ and the map in the middle are injective. Therefore, $f$ is injective. This immediately yields that
$S^+_a \seqv{U}$ is a new quiver flag of $S^+_a M = S^+_a U^\nu$. The problem is that the
dimension of $S^+_a U^i$ is dependent on the rank of $\phi^{U^i}_a$.
This motivates the next definition. Recall that $d_a - \rank \phi^X_a = \dim \Hom(X, S_a)$
for a representation $X$ of dimension vector $\dvec{d}$.
\begin{defn}
  Let $a$ be a sink, $\dvec{d}$ a dimension vector and $s$ an integer. Define
    \[ \Rep[Q]{\dvec{d}}\lrang{a}^{s} :=
    \Set{ M \in \Rep[Q]{\dvec{d}} | \dim \Hom(M, S_a) = s}.\]

    Let $\seqv{\dvec{d}}$ be a $(\nu+1)$-tuple of dimension vectors
    and
    $\seqv{r} = (r^0, r^1, \dots , r^\nu)$ be a $(\nu+1)$-tuple of integers.
    For each representation $M$ define
\[
\Fl[Q]{\seqv{\dvec{d}}}{M}\lrang{a}^{\seqv{r}} :=
\Set{ \seqv{U} \in \Fl[Q]{\seqv{\dvec{d}}}{M} |
    \dim \Hom(U^i, S_a) = r^i}.
\]

     Moreover, let
     \[ \Rep[Q]{\dvec{d}}\lrang{a} := \Rep[Q]{\dvec{d}}\lrang{a}^0\] and
     \[ \Fl[Q]{\seqv{\dvec{d}}}{M}\lrang{a} := \Fl[Q]{\seqv{\dvec{d}}}{M}\lrang{a}^{\seqv{0}}.\]
\end{defn}
\begin{remark}
    Recall that a filtration $\flvec{d}$ of $\dimve M$ is a sequence of dimension
    vectors such that $\dvec{d}^0=0$, $\dvec{d}^\nu = \dimve M$ and that $\dvec{d}^i \le \dvec{d}^{i+1}$.
  In order to know a filtration $\flvec{d}$ it is enough to know the terms
  $\dvec{d}^1, \dots, \dvec{d}^{\nu-1}$, since $\dvec{d}^0$ is always $0$ and
  $\dvec{d}^\nu$ is always $\dimve M$. Therefore, we identify the $(\nu-1)$-tuple
  $(\dvec{d}^1, \dots, \dvec{d}^{\nu-1})$ with the $(\nu+1)$-tuple $\flvec{d}$.
\end{remark}
\begin{beispiel}
    Let
\[Q =
\makeatletter%
\let\ASYencoding\f@encoding%
\let\ASYfamily\f@family%
\let\ASYseries\f@series%
\let\ASYshape\f@shape%
\makeatother%
\setlength{\unitlength}{1pt}
\includegraphics{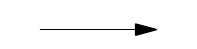}%
\definecolor{ASYcolor}{gray}{0.000000}\color{ASYcolor}
\fontsize{12.000000}{14.400000}\selectfont
\usefont{\ASYencoding}{\ASYfamily}{\ASYseries}{\ASYshape}%
\ASYalign(-52.462417,5.372285)(-0.500000,-0.500000){$1$}
\ASYalign(-4.443095,5.372285)(-0.500000,-0.500000){$2$}
\ASYalign(-28.452756,8.985785)(-0.500000,0.000000){$\scriptstyle\alpha$}
.\]
Consider the representation $M$ given by $M_1 = M_2 = K^2$ and
$M_\alpha =\left(\begin{smallmatrix}
1 & 0\\
0 & 0\\
\end{smallmatrix}\right)$. We have that
$M \in \Rep[Q]{(2,2)}\lrang{2}^1$. Now consider flags of type $((0,0), (1,1), (2,2))$, i.e.
subrepresentations $N$ of dimension vector $(1,1)$. We need two injective linear maps
$f_1, f_2 \colon K^1 \rightarrow K^2$ making the following diagram commutative.
\[
\xymatrix{
K \ar[r]^{N_\alpha} \ar[d]_{f_1} & K \ar[d]^{f_2}\\
K^2 \ar[r]_{\left(\begin{smallmatrix}
    1 & 0\\
    0 & 0
\end{smallmatrix}\right)} & K^2
}
\]
We have the following situations.
\label{bsp:grass}
\begin{itemize}
    \item $N \in \Gr[Q]{(1,1)}{M}\lrang{2}^1$: This means that $N_\alpha = 0$. Therefore,
        we need that the image of $f_1$ is in
        the kernel of $M_\alpha$, which is $1$-dimensional. Hence, a subrepresentation in
        $\Gr[Q]{(1,1)}{M}\lrang{2}^1$ is given by $f_1 = \left( \begin{smallmatrix}
            0\\
            1
        \end{smallmatrix}\right)$ and $f_2$ being an arbitrary inclusion. The point
        $f_1 = f_2 = \left( \begin{smallmatrix}
            0\\
            1
        \end{smallmatrix}\right)$ is special, since for this inclusion we have
        that $M/N \cong S_1 \oplus S_2$ and otherwise $M/N \cong K \overset{1}{\rightarrow} K$.
    \item $N \in \Gr[Q]{(1,1)}{M}\lrang{2}^0$: This means that $N_\alpha \neq 0$. Therefore,
        we need that the image of $f_1$ is not in
        the kernel of $M_\alpha$, which is $1$-dimensional. Hence, a subrepresentation in
        $\Gr[Q]{(1,1)}{M}\lrang{2}^0$ is given by $f_1=f_2= \left( \begin{smallmatrix}
            1\\
            x
        \end{smallmatrix}\right)$ for any $x \in K$.
\end{itemize}
The variety $\Gr[Q]{(1,1)}{M}$ consists therefore of two $\Proj^1_K$ glued together at one point.
Graphically,
\[
\Gr[Q]{(1,1)}{M} = {\color{red}\Gr[Q]{(1,1)}{M}\lrang{2}^1} \amalg {\color{blue}\Gr[Q]{(1,1)}{M}\lrang{2}^0} \cong
\begin{array}{c@{}}
    \begin{tikzpicture}
        \draw[blue] (2,0) circle (1);
        \draw[red] (0,0) circle (1);
        \filldraw[fill=red,red] (1,0) circle (0.05); 
    \end{tikzpicture}
\end{array}.
\]
Note that the Grassmannian is neither irreducible nor smooth.
\end{beispiel}

In order to get rid of $\seqv{r}$ we define the following maps and then look at the fibres.
\begin{defn}
    Let $a$ be a sink, $\dvec{d}$ a dimension vector, $s\in \N$ and $M \in \Rep[Q]{\dvec{d}}\lrang{a}^s$.
    We have that $M \cong M' \oplus S_a^s$ for some element $M' \in \Rep[Q]{\dvec{d} - s\epsilon_a}\lrang{a}$.
    Without loss of generality we can assume that $M=M' \oplus S_a^s$ and we set
    $\pi_a M := M'$. Obviously, $\pi_a M$ is unique up to isomorphism.

  Now let $\flvec{d}$ be a filtration
    and $\seqv{r}=(r^0, \dots, r^\nu)$ a $(\nu +1)$-tuple of integers.
    Define
    \begin{alignat*}{3}
    \pi_a^{\seqv{r}} &\colon\quad&  \Fl[Q]{\seqv{\dvec{d}}}{M}\lrang{a}^{\seqv{r}} &\rightarrow
      \Fl[Q]{\seqv{\dvec{d}} - \seqv{r} \epsilon_a}{\pi_a M}\lrang{a}\\
      && \seqv{U} &\mapsto \seqv{V} & & \text{where } V^i_j:=\begin{cases}
      U^i_j & \text{ if } j \neq a,\\
      \Bild \phi^{U^i}_a & \text{ if } j = a.
    \end{cases}
    \end{alignat*}
    \label{def:flagtoredfib}
\end{defn}
\begin{remark}
    Note that $\pi_a M$ is $\iota_{a,M} S^-_a S^+_a M$.
\end{remark}
\begin{beispiel}
    Coming back to example \ref{bsp:grass} we see that $\pi_2^1$ collapses $\Gr[Q]{(1,1)}{M}\lrang{2}^1$
    to the point
    \[\Gr[Q]{(1,0)}{K^2 \overset{\left(\begin{smallmatrix}
        1 & 0\\
    \end{smallmatrix}\right)}{\rightarrow} K}\lrang{2}.\]
    The fibre of $\pi_2^1$ over this point is the vector space Grassmannian $\Gr{1}{K^2}$, being
    isomorphic to $\Proj^1_K$.
\end{beispiel}

We now introduce a little bit more notation. If $\seqv{d}$ is a sequence, then denote by
$\overleftarrow{\seqv{d}}$ the sequence given by $(\overleftarrow{\seqv{d}})^i = d^{\nu-i}$. Moreover,
we define the sequence $\seqv{e}$ by $e^i := d^\nu - d^{\nu-i}$.
Therefore, if $\seqv{\dvec{d}}$ is a filtration of $\dvec{d}^\nu$,
then $\seqv{\dvec{e}}$ is a filtration of $\dvec{d}^\nu$.

The fibre of the map $\pi_a^{\seqv{r}}$ is a set of the following type.
\begin{defn}
  Let $\seqv{e} = (e^0, e^1, \dots, e^\nu)$ and $\seqv{r} = (r^0, r^1, \dots , r^\nu)$ be sequences of non-negative integers such that
  $\seqv{e} + \overleftarrow{\seqv{r}}$ is a filtration.
  Let $A_{\nu+1}$ be the quiver
  \[
  0 \rightarrow 1 \rightarrow 2 \rightarrow 3 \rightarrow \dotsm \rightarrow \nu.
  \]
  Then define
  \[
  \xymatrix{
  X^{\seqv{r}, \seqv{e}} := &
  K^{e^{\nu}+r^0} \ar@{->>}[r] &K^{e^{\nu-1}+r^1} \ar@{->>}[r] & \cdots \ar@{->>}[r] & K^{e^1 + r^{\nu-1}}\ar@{->>}[r]  & K^{e^0 + r^\nu}.
  }
  \]
\end{defn}
\begin{remark}
    The $K$-representation $X^{\seqv{r}, \seqv{e}} \in \repK{A_{\nu +1}}{K}$ is injective and its
    isomorphism class does not depend on the choice of the surjections.
\end{remark}
\begin{lmm}
  Let $\seqv{e} = (e^0, e^1, \dots, e^\nu)$
  and $\seqv{r} = (r^0, r^1, \dots , r^\nu)$ be sequences of non-negative integers such that
  $\seqv{e} + \overleftarrow{\seqv{r}}$ is a filtration.
  Then $X^{\seqv{r}, \seqv{e}}$ has
  a subrepresentation of dimension vector $\seqv{r}$ if and only if $\seqv{e}$ is a filtration of $e^\nu$
  \label{lmm:flagfib_not_empty}. Moreover, if $K$ is a finite field with $q$ elements, then
  the number of $K$-subrepresentations is given by
  \[
  \# \Gr[A_{\nu+1}]{\seqv{r}}{X^{\seqv{r},\seqv{e}}} = \prod_{i=0}^\nu \qbinom{e^{\nu-i} - e^{\nu-i-1} + r^i}{r^i}_q.
  \]
  In particular, this number is equal to $1$ modulo $q$ if and only if the set of subrepresentations is non-empty.
\end{lmm}
\begin{proof}
  We prove this by induction on $\nu$.
  \begin{description}
    \item[$\nu=0$]
      There is a subspace of dimension $r^0$ of $K^{r^0 + e^0}$ if and only if $e^0 \ge 0$ and, for
      $K$ a finite field of cardinality $q$, the number of
      those is obviously $\qbinom{e^0 + r^0}{r^0}_q$.
    \item[$\nu \ge 1$]
      If $(U^0,U^1,U^2, \dots, U^\nu)$ is a subrepresentation of dimension vector $\seqv{r}$ of $X^{\seqv{r}, \seqv{e}}$,
      then
      $(U^1,U^2, \dots, U^\nu)$ is a subrepresentation of
      $X^{(r^1, r^2, \dots, r^\nu), (e^0, e^1, \dots, e^{\nu-1})}$ of dimension vector $(r^1, r^2,\dots, r^\nu)$. Therefore, by induction,
      $e^i \le e^{i+1}$ for $0 \le i < \nu-1$ and $0 \le e^{0}$. The preimage $V$ of $U^1$ under the surjection from $U^0$ has
      dimension $r^1 + ( (e^{\nu} + r^0) - (e^{\nu-1} + r^1))= e^\nu - e^{\nu-1} + r^0$. Since $\seqv{U}$ is a subrepresentation, we must have that
      $U^0 \subset V$. Therefore, $r^0 \le e^\nu - e^{\nu-1} + r^0$ or equivalently $e^{\nu -1} \le e^\nu$.

      On the other hand, if $e^0 \ge 0$ and $e^i \le e^{i+1}$ for all $0 \le i < \nu$, then, by induction, there is a subrepresentation
      $(U^1,U^2, \dots, U^\nu)$ of $X^{(r^1, r^2, \dots, r^\nu), (e^0, e^1, \dots, e^{\nu-1})}$ of dimension vector $(r^1, r^2,\dots, r^\nu)$.
      As before, the dimension of the preimage $V$ of $U^1$ under the surjection from $U^0$ has dimension
      $e^\nu - e^{\nu-1} + r^0 \ge r^0$. If we choose now any subspace $U^0$ of dimension $r^0$ in $V$, then we obtain
      a subrepresentation
      of $X^{\seqv{r}, \seqv{e}}$ of dimension vector $\seqv{r}$.
      
      If $K$ is a finite field of cardinality $q$, then, by induction, we have that
      the number of subrepresentations
      of $X^{(r^1, r^2, \dots, r^\nu), (e^0, e^1, \dots, e^{\nu-1})}$
      of dimension vector $(r^1, r^2,\dots, r^\nu)$ is equal to
  \[
  \prod_{i=1}^\nu \qbinom{e^{\nu-i} - e^{\nu-i-1} + r^i}{r^i}_q.
  \]
  To complete such a subrepresentation to a subrepresentation of $X^{\seqv{r}, \seqv{e}}$ we
  have to choose an $r^0$-dimensional subspace of an $(e^{\nu} - e^{\nu-1} + r^0)$-dimensional
  space. Therefore, the number of subrepresentations is equal to
  \[
  \qbinom{e^{\nu} - e^{\nu-1} + r^0}{r^0}_q \prod_{i=1}^\nu \qbinom{e^{\nu-i} - e^{\nu-i-1} + r^i}{r^i}_q.
  \]
  This yields the claim.
  \end{description}
\end{proof}

\begin{theom}
    Let $a$ be a sink of $Q$, $\flvec{d}$ a filtration, $\seqv{r}$ a $(\nu+1)$-tuple of non-negative integers
    and $M \in \Rep[Q]{\dvec{d}^\nu}$. Then
    \[\pi_a^{\seqv{r}} \colon
      \Fl[Q]{\seqv{\dvec{d}}}{M}\lrang{a}^{\seqv{r}} \rightarrow
      \Fl[Q]{\seqv{\dvec{d}} - \seqv{r} \epsilon_a}{\pi_a M}\lrang{a}\]
      is surjective and
      the fibre $(\pi_a^{\seqv{r}})^{-1}(\seqv{U})$
      over any $\seqv{U} \in \Fl[Q]{\seqv{\dvec{d}} - \seqv{r} \epsilon_a}{\pi_a M}\lrang{a}$
      is isomorphic to $\Gr[A_\nu]{\seqv{r}}{X^{\seqv{r}, \seqv{\dvec{e}}_a}}$,
      where $\dvec{e}^{\nu -i} := \dvec{d}^\nu - \dvec{d}^i$.
      In particular the number
      of points in the fibre only depends on $\seqv{r}$ and $\seqv{\dvec{d}}$ and not on $\seqv{U}$.
     \label{theom:qgrass_to_red_fib}
\end{theom}
\begin{proof}
  Fix a flag $\seqv{V} \in \Fl[Q]{\seqv{\dvec{d}} - \seqv{r} \epsilon_a}{\pi_a M}\lrang{a}$. Let
  $\seqv{U} \in \Fl[Q]{\seqv{\dvec{d}}}{M}\lrang{a}^{\seqv{r}}$. The flag $\seqv{U}$ is
  in $(\pi_a^{\seqv{r}})^{-1}(\seqv{V})$ if and only if $U_j^i = V_j^i$ for all $j \neq a$
  in which case $\Bild(\Phi^{U^i}_a) = V^i_a$.
  Therefore, we only have to choose $U^i_a \subset M_a$ such that $V^i_a \subseteq U^i_a$,
  $U^{i-1}_a \subset U^i_a$ and $\dim U^i_a = d^i_a$ for all $i$. This is the same as choosing
  $\overline{U^{i}_a} \subset M_a/V^i_a$ such that $\theta^i(\overline{U^{i-1}_a}) \subset \overline{U^i_a}$ and
  $\dim \overline{U^i_a} = r^i$ if we denote by $\theta^i \colon M_a/V^{i-1}_a \rightarrow M_a/V^i_a$ the
  canonical projection. This is equivalent to finding a subrepresentation of
  \[ M_a/V^0_a \rightarrow M_a/V^1_a \rightarrow \cdots \rightarrow M_a/V^\nu_a \]
  of dimension vector $\seqv{r}$. All the maps in this representation are surjective since $\seqv{V}$
  is a flag, therefore this representation
  of $A_{\nu+1}$ is isomorphic to
  $X^{\seqv{r}, \seqv{e}}$. Since $\flvec{d}$ is a filtration of $M$ we have that $\flvec{e}$ is a filtration.
  Therefore, $\Gr[A_\nu]{\seqv{r}}{X^{\seqv{r}, \seqv{\dvec{e}}_a}}$ is non-empty by lemma \ref{lmm:flagfib_not_empty}
  and $\pi_a^\seqv{r}$ is surjective.
\end{proof}
Now we are nearly ready to do reflections. The only thing left to define is what happens on a source. If $b$ is a source in
$Q$, then $b$ is a sink in $Q^{op}$, so we just dualise everything.
\begin{defn}
  Let $\seqv{U} \in \Fl[Q]{\seqv{\dvec{d}}}{M}$ and let $\dvec{e}^{\nu-i} = \dvec{d}^\nu - \dvec{d}^i$.
  Then define
  \begin{alignat*}{2}
    \hat D &\colon\quad &
    \Fl[Q]{\seqv{\dvec{d}}}{M}&\rightarrow
    \Fl[Q^{op}]{\seqv{\dvec{e}}}{DM}\\
    &&\seqv{U} & \mapsto (\hat D (\seqv{U}))^i := \ker( DM \rightarrow D(U^i)) = D(M/U^i).
  \end{alignat*}
\end{defn}
\begin{remark}
    Obviously, $\hat D ^2 = \id$ and the map $\hat D$ is an isomorphism of varieties.
\end{remark}
\begin{defn}
  Let $b$ be a source, $\dvec{d}$ a dimension vector and $s$ an integer. Define
    \[ \Rep[Q]{\dvec{d}}\lrang{b}^{s} :=
    \Set{ M \in \Rep[Q]{\dvec{d}} | \dim \Hom(S_b, M) = s}.\]

    Let $\seqv{\dvec{d}}$ be a $(\nu+1)$-tuple of dimension vectors
    and
    $\seqv{r} = (r^0, r^1, \dots , r^\nu)$ be a $(\nu+1)$-tuple of integers.
    For each representation $M$ define
\[
\Fl[Q]{\seqv{\dvec{d}}}{M}\lrang{b}^{\seqv{r}} :=
\Set{ \seqv{U} \in \Fl[Q]{\seqv{\dvec{d}}}{M} |
    \dim \Hom(S_b, M/U^i) = r^i}.
\]

     Moreover, let
     \[ \Rep[Q]{\dvec{d}}\lrang{b} := \Rep[Q]{\dvec{d}}\lrang{b}^0\] and
     \[ \Fl[Q]{\seqv{\dvec{d}}}{M}\lrang{b} := \Fl[Q]{\seqv{\dvec{d}}}{M}\lrang{b}^{\seqv{0}}.\]
\end{defn}
%
%
\begin{remark}
    Note that $\seqv{U} \in \Fl[Q]{\seqv{\dvec{d}}}{M}\lrang{b}^{\seqv{r}}$ if and only if
    $\hat D \seqv{U} \in \Fl[Q^{op}]{\seqv{\dvec{e}}}{DM}\lrang{b}^{\seqv{r}}$.
\end{remark}
Now we state the main result on reflections.
\begin{theom}
    Let $a$ be a sink of $Q$, $\flvec{d}$ be a filtration and $M \in \Rep[Q]{\dvec{d}^\nu}\lrang{a}$.
  The map
  \begin{alignat*}{2}
    S^+_a &\colon\quad & \Fl[Q]{\seqv{\dvec{d}}}{M}\lrang{a} &\rightarrow
    \Fl[\sigma_a Q]{\sigma_a \seqv{\dvec{d}}}{S^+_a M}\lrang{a}\\
    && \seqv{U} & \mapsto S^+_a \seqv{U}
  \end{alignat*}
  is an isomorphism of varieties with inverse $\hat D \circ S^+_a \circ \hat D = S^-_a$.
  \label{theom:reflflagiso}
\end{theom}
\begin{proof}
  First, we show that $S^+_a \seqv{U}$ lies in the correct set. Let $\dvec{e}^{\nu-i} = \dvec{d}^\nu - \dvec{d}^i$.
For each $i$, we have the following commutative diagram with exact columns.
  \[
\begin{CD}
      @.          0        @. 0 @. 0 @. \\
    @.          @VVV                        @VVV           @VVV     @.\\
    0 @>>>  (S^+_a U^{i})_a      @>>> \bigoplus\limits_{j \rightarrow a} U^{i}_j @>{\phi^{U^{i}}_a}>> U^{i}_a @>>> 0\\
    @.          @VVV                        @VVV           @VVV     @.\\
  0 @>>> (S^+_a M)_a @>>> \bigoplus\limits_{j \rightarrow a} M_j @>\phi^{M}_a>> M_a @>>> 0\\
    @.          @VVV                        @VVV           @VVV     @.\\
   0 @>>> (S^+_a M / S^+_a U^{i})_a  @>>> \bigoplus\limits_{j \rightarrow a} (M/U^{i})_j @>\phi^{M/U^{i}}_a>> (M/U^{i})_a @>>> 0\\
    @.          @VVV                        @VVV           @VVV     @.\\
     @.          0        @. 0 @. 0 @.
\end{CD}
  \]
  The two top rows are exact since $\seqv{U} \in \Fl[Q]{\seqv{\dvec{d}}}{M}\lrang{a}$. By the snake lemma,
  we have that the
  bottom row is exact. Therefore, the map
  \[ (S^+_a M / S^+_a U^{i})_a \rightarrow \bigoplus_{j\rightarrow a} (M/U^{i})_j \]
  is injective and hence $S^+_a(M)/S^+_a(U^i) \in \Rep[\sigma_a Q]{\sigma_a(\dvec{e}^{\nu-i})}\lrang{a}$.
  The diagram also yields that $\hat D \circ S^+_a \circ \hat D \circ S^+_a = \id$.
  Since $S^+_a$ is a functor and all choices where natural, we have that $S^+_a$ gives
  a natural transformation between the functors of points of these two varieties. Therefore,
  it is a morphism of varieties.

  Dually, $S^+_a \circ \hat D \circ S^+_a \circ \hat D = \id$. This concludes the proof.
%
\end{proof}

\section{Reflections and Hall Numbers}
Let $\F_q$ be the finite field with $q$ elements and
$Q$ a quiver. Let $w=(i_r, \dots, i_1)$ be a word in vertices of $Q$. Recall that
\[
u_w = u_{i_r} \diamond \dots \diamond u_{i_2} \diamond u_{i_1} = \sum_{X} F_w^{X} u_X.
\]
Define a filtration $\flvec{d}(w)$ by letting
\[\dvec{d}(w)^k := \sum_{j=1}^k \epsilon_{i_k}.\]
Then we obviously have $F_w^X = \# \Fl[Q]{\flvec{d}(w)}{X}$. Therefore, coefficients
in the Hall algebra are closely related to counting points of quiver
flags over finite fields. In the following, we will use reflection functors
to simplify the problem of counting the number of points modulo $q$. As an
application, we will show that for a preprojective or preinjective representation $X$
we have that $\# \Fl[Q]{\flvec{d}}{X} = 1 \mod q$ if $\Fl[Q]{\flvec{d}}{X}$ is non-empty.

\begin{lmm}
    Let $a$ be a sink of $Q$, $K$ a field, $\flvec{d}$ a filtration and $M \in \RepK[Q]{\dvec{d}^\nu}{K}$.
    Then
	\[
    \# \Fl[Q]{\seqv{\dvec{d}}}{M}
	=\sum_{\seqv{r}\ge 0} \# \Gr[A_\nu]{\seqv{r}}{X^{\seqv{r}, \seqv{e}_a}}
    \# \Fl[Q]{\seqv{\dvec{d}} - \seqv{r} \epsilon_a}{\pi_a M}\lrang{a}
    \]
    (on both sides we possibly have $\infty$).

    We further note that, for each sequence of non-negative integers $\seqv{r}$, if
    $
    \Fl[Q]{\seqv{\dvec{d}} - \seqv{r} \epsilon_a}{\pi_a M}\lrang{a}
    $
    is non-empty,
    then so is
    $
    \Gr[A_\nu]{\seqv{r}}{X^{\seqv{r}, \seqv{e}_a}}.
    $
    \label{reflflag:lmm:decompflag}
\end{lmm}
\begin{proof}
    We have that
    \[
    \Fl[Q]{\seqv{\dvec{d}}}{M} = \coprod_{\seqv{r}\ge 0} \Fl[Q]{\seqv{\dvec{d}}}{M}\lrang{a}^{\seqv{r}}.
    \]
    By theorem \ref{theom:qgrass_to_red_fib}, we have for each sequence of non-negative integers $\seqv{r}$ that
    \[
	\# \Fl[Q]{\seqv{\dvec{d}}}{M}\lrang{a}^{\seqv{r}}
	=
	\# \Gr[A_\nu]{\seqv{r}}{X^{\seqv{r}, \seqv{e}_a}} \# \Fl[Q]{\seqv{\dvec{d}} - \seqv{r} \epsilon_a}{\pi_a M}\lrang{a}.
    \]
    The same theorem yields that, if 
    $\Fl[Q]{\seqv{\dvec{d}} - \seqv{r} \epsilon_a}{\pi_a M}\lrang{a}$ is non-empty, then so
    is $\Gr[A_\nu]{\seqv{r}}{X^{\seqv{r}, \seqv{e}_a}}$.
\end{proof}
\begin{lmm}
    Let $a$ be a sink of $Q$, $K$ a field, $\flvec{d}$ a filtration
    and $M \in \RepK[Q]{\dvec{d}^\nu}{K}\lrang{a}^s$. Let $\seqv{r}_+= \seqv{r}_+(\flvec{d})$ be given as follows:
    \begin{align*}
        r^0_+ &:= 0; &\\
        r^i_+ &:= \max\{0, (\sigma_a(\dvec{d}^{i-1} - \dvec{d}^i))_a + r^{i-1}_+\}& \text{for } 0 < i < \nu;\\
        r^\nu_+ &:= s.
    \end{align*}
    Now let $\seqv{r}$ be a sequence of integers. If
    $\Fl[Q]{\seqv{\dvec{d}}}{M}\lrang{a}^{\seqv{r}}$ is non-empty, then $\seqv{r} \ge \seqv{r}_+$.
    \label{reflflag:lmm:rplusmin}
\end{lmm}
\begin{proof}
    Let $\seqv{U} \in \Fl[Q]{\seqv{\dvec{d}}}{M}\lrang{a}^{\seqv{r}}$. By definition,
    $\Hom(U^i, S_a) = r^i$. We have
    \[r^i= \codim \Bild \Phi_a^{U^i} = \dim \ker \Phi_a^{U^i} + d^i_a - \sum_{j \rightarrow a} d^i_j
    = \dim \ker \Phi_a^{U^i} - (\sigma_a \dvec{d}^i)_a.\]
    We prove $r^i \ge r_+^i$ by induction on $i$. For $i=0$ the claim
    is obviously true. Now let $0\le i \le \nu-2$. Obviously, $\dim \ker \Phi_a^{U^i} \le \dim \ker \Phi_a^{U^{i+1}}$.
    Therefore,
    \[ r^i_+  + (\sigma_a \dvec{d}^i)_a 
    \le r^i  + (\sigma_a \dvec{d}^i)_a= \dim \ker \Phi_a^{U^i} \le \dim \ker \Phi_a^{U^{i+1}} =
    r^{i+1} + (\sigma_a \dvec{d}^{i+1})_a.\]
    Hence, $r^{i+1} \ge \max\{0,(\sigma_a(\dvec{d}^{i} - \dvec{d}^{i+1}))_a + r^i_+\} = r^{i+1}_+$.

    For $r^\nu_+$ note that, by definition, $U^\nu = M$ and therefore
    \[ r^\nu =\codim \Bild \Phi_a^{U^\nu} = \codim \Bild \Phi_a^{M} = s.\]
\end{proof}
\begin{remark}
    Note that $\seqv{r}_+(\flvec{d}-\seqv{r}_+(\flvec{d}) \epsilon_a)=0$ since
    \[
    \sigma_a (\dvec{d}^i - \dvec{d}^{i-1})_a + r^i_+(\flvec{d}) - r^{i-1}_+(\flvec{d})
    \ge r^i_+(\flvec{d}) - r^i_+(\flvec{d}) = 0.
    \]
    For any filtration $\flvec{d}$ of some representation $M$ it is enough to remember the terms
    \[
    (\dvec{d}^1, \dots, \dvec{d}^{\nu-1})
    \]
    since we always have $\dvec{d}^0 = 0$
    and $\dvec{d}^\nu = \dimve M$. Note that the rule to construct $r_+^i$ for
    $0<i<\nu$
    depends neither on $\dvec{d}^0$ nor on $\dvec{d}^\nu$. Therefore, we
    can define
    \[S_a^+ \flvec{d} = S_a^+ (\dvec{d}^1, \dots, \dvec{d}^{\nu-1}) :=
    (\sigma_a \dvec{d}^1 + r_+^1 \epsilon_a, \dots, \sigma_a \dvec{d}^{\nu-1} + r_+^{\nu-1}). \]
    If $\flvec{d}$ is
    a filtration of $M$, then $S_a^+ \flvec{d}$ is a filtration of $S_a^+ M$ if and only if
    $(S_a^+ \flvec{d})^{\nu-1} \le \dimve S_a^+ M$.
    \label{reflflag:rem:reflfilt}
\end{remark}
\begin{cory}
    Let $a$ be a sink, $K$ a field, $\flvec{d}$ a filtration and $M \in \RepK[Q]{\dvec{d}^\nu}{K}\lrang{a}^s$.
    Then
	\[
    \# \Fl[Q]{\seqv{\dvec{d}}}{M}
    =\sum_{\seqv{r}\ge 0} \# \Gr[A_\nu]{\seqv{r}+\seqv{r}_+}{X^{\seqv{r}+\seqv{r}_+, \seqv{e}_a}}
    \# \Fl[\sigma_a Q]{S_a^+\seqv{\dvec{d}}+ \seqv{r}\epsilon_a}{S_a^+ M}\lrang{a}.
    \]

    In particular, if $K$ is a finite field of cardinality $q$, we have
    \[
    \# \Fl[Q]{\seqv{\dvec{d}}}{M}
    =\sum_{\seqv{r}\ge 0}
    \# \Fl[\sigma_a Q]{S_a^+\seqv{\dvec{d}}+\seqv{r}\epsilon_a}{S_a^+ M}\lrang{a} \mod q.
    \]
    \label{reflflag:cory:decompreflmodq}
\end{cory}
\begin{proof}
    By lemmas \ref{reflflag:lmm:decompflag} and \ref{reflflag:lmm:rplusmin} we obtain
    that
	\[
    \# \Fl[Q]{\seqv{\dvec{d}}}{M}
    =\sum_{\seqv{r}\ge 0} \# \Gr[A_\nu]{\seqv{r}+\seqv{r}_+}{X^{\seqv{r}+\seqv{r}_+, \seqv{e}_a}}
    \# \Fl[Q]{\flvec{d} - (\seqv{r} + \seqv{r}_+)\epsilon_a}{\pi^s_a M}\lrang{a}.
    \]
    Note that $\sigma_a (\dvec{d}^i - (r^i+r^i_+)\epsilon_a) =
    (S_a^+\flvec{d})^i + r^i \epsilon_a$ for all $0 < i < \nu$. Therefore, theorem \ref{theom:reflflagiso}
    yields that
    \[
    \Fl[Q]{\flvec{d} - (\seqv{r} + \seqv{r}_+)\epsilon_a}{\pi^s_a M}\lrang{a} \cong
    \Fl[\sigma_a Q]{S_a^+\flvec{d} + \seqv{r}\epsilon_a}{S_a^+ M}\lrang{a}.
    \]
    This proves the first claim.

    Now let $K$ be a finite field of cardinality $q$. If
    $\Gr[A_\nu]{\seqv{r}+\seqv{r}_+}{X^{\seqv{r}+\seqv{r}_+, \seqv{e}_a}}$
    is non-empty, then its number is one modulo $q$ by lemma \ref{lmm:flagfib_not_empty}.
    The second part of lemma \ref{reflflag:lmm:decompflag} yields that whenever
    $\Fl[Q]{\flvec{d} - (\seqv{r} + \seqv{r}_+)\epsilon_a}{\pi^s_a M}\lrang{a}$
    is non-empty, then $\Gr[A_\nu]{\seqv{r}+\seqv{r}_+}{X^{\seqv{r}+\seqv{r}_+, \seqv{e}_a}}$
    is non-empty. This finishes the proof.
\end{proof}
We obtain the following.
\begin{theom}
    \label{theom:reflflagmodq}
    Let $a$ be a sink, $K$ a field, $\flvec{d}$ a filtration and $M \in \RepK[Q]{\dvec{d}^\nu}{K}\lrang{a}^s$.
	Then $\Fl[Q]{\seqv{\dvec{d}}}{M}$ is empty if and only if
    $\Fl[\sigma_a Q]{S^+_a \flvec{d}}{S^+_a M}$ is.

	Moreover, if $K=\F_q$ is a finite field, then
    \[
	\# \Fl[Q]{\seqv{\dvec{d}}}{M}=
	\# \Fl[\sigma_a Q]{S^+_a \flvec{d}}{S^+_a M} \mod q.
    \]
\end{theom}
\begin{proof}
    By corollary \ref{reflflag:cory:decompreflmodq} we obtain that
	\[
    \# \Fl[Q]{\seqv{\dvec{d}}}{M}
    =\sum_{\seqv{r}\ge 0} \# \Gr[A_\nu]{\seqv{r}+\seqv{r}_+}{X^{\seqv{r}+\seqv{r}_+, \seqv{e}_a}}
    \# \Fl[\sigma_a Q]{S_a^+\seqv{\dvec{d}}+ \seqv{r}\epsilon_a}{S_a^+ M}\lrang{a}.
    \]
    Note that $S_a^+ \flvec{d}$ is a filtration of $S_a^+M$ if and only if
    $(S_a^+ \flvec{d})^{\nu-1} \le \dimve S_a^+ M$. Therefore, if $S_a^+ \flvec{d}$
    is not a filtration of $S_a^+ M$, then 
    each $\Fl[\sigma_a Q]{S_a^+\seqv{\dvec{d}}+ \seqv{r}\epsilon_a}{S_a^+ M}\lrang{a}$
    is empty for all $\seqv{r} \ge 0$ and hence, so is 
    $\Fl[Q]{\seqv{\dvec{d}}}{M}$. In this case, we also have that
    $\Fl[\sigma_a Q]{S^+_a \flvec{d}}{S^+_a M}$ is empty. Both claims follow.

    Assume now that $S_a^+ \flvec{d}$ is a filtration of $S_a^+M$.
    We have that
    \[\Fl[\sigma_a Q]{S^+_a \flvec{d}}{S^+_a M} \cong \Fl[\sigma_a Q^{op}]{\dimve S^+_a M - \rever{S^+_a \flvec{d}}}{DS^+_a M}\]
    via $\hat D$. Let $\flvec{f}:=\dimve S^+_a M - \rever{S^+_a \flvec{d}}$. By using lemma
    \ref{reflflag:lmm:decompflag} we obtain that
    \[
    \# \Fl[\sigma_a Q^{op}]{\flvec{f}}{DS^+_a M} 
    =\sum_{\seqv{r}\ge 0} \# \Gr[A_\nu]{\rever{\seqv{r}}}{X^{\rever{\seqv{r}}, (S^+_a\flvec{d})_a}}
    \# \Fl[\sigma_a Q^{op}]{\flvec{f} - \rever{\seqv{r}} \epsilon_a}{D S_a^+ M}\lrang{a}.
    \]
    Moreover, the same lemma yields that for each $\seqv{r}\ge 0$ we have that if
    $\Fl[\sigma_a Q^{op}]{\flvec{f} - \rever{\seqv{r}} \epsilon_a}{D S_a^+ M}\lrang{a}$ is
    non-empty, then $\Gr[A_\nu]{\rever{\seqv{r}}}{X^{\rever{\seqv{r}}, (S^+_a\flvec{d})_a}}$ is non-empty.
    Using $\hat D$ yields
    \[
    \Fl[\sigma_a Q^{op}]{\flvec{f} - \rever{\seqv{r}} \epsilon_a}{D S_a^+ M}\lrang{a}
    \cong \Fl[\sigma_a Q]{S_a^+\flvec{d} + \seqv{r}\epsilon_a}{S_a^+ M}\lrang{a}.
    \]
    Combining these equalities, we have that
    \[
    \# \Fl[\sigma_a Q]{S^+_a \flvec{d}}{S^+_a M}
    =\sum_{\seqv{r}\ge 0} \# \Gr[A_\nu]{\rever{\seqv{r}}}{X^{\rever{\seqv{r}}, (S^+_a\flvec{d})_a}}
    \# \Fl[\sigma_a Q]{S_a^+\flvec{d} + \seqv{r}\epsilon_a}{S_a^+ M}\lrang{a}.
    \]
    Therefore, $\Fl[\sigma_a Q]{S^+_a \flvec{d}}{S^+_a M}$ is empty if and only if for all
    $\seqv{r} \ge 0$ the variety $\Fl[\sigma_a Q]{S_a^+\flvec{d} + \seqv{r}\epsilon_a}{S_a^+ M}\lrang{a}$
    is empty. The same is true for $\Fl[Q]{\seqv{\flvec{d}}}{M}$ and this proves the first claim.

    Now let $K$ be  a finite field with $q$ elements. We already saw that if
    $\Fl[\sigma_a Q]{S_a^+\flvec{d} + \seqv{r} \epsilon_a}{S_a^+ M}\lrang{a}$
    is non-empty, then
    $\Gr[A_\nu]{\rever{\seqv{r}}}{X^{\rever{\seqv{r}}, (S^+_a\flvec{d})_a}}$ is non-empty.
    Therefore, lemma \ref{lmm:flagfib_not_empty} yields that
    \[\# \Fl[\sigma_a Q]{S^+_a \flvec{d}}{S^+_a M}
    =\sum_{\seqv{r}\ge 0}
    \# \Fl[\sigma_a Q]{S_a^+\flvec{d} + \seqv{r}\epsilon_a}{S_a^+ M}\lrang{a} \mod q.\]
    By corollary \ref{reflflag:cory:decompreflmodq} this is equal to
    $\# \Fl[Q]{\flvec{d}}{M}$. This finishes the proof.
\end{proof}

\begin{remark}
    The Coxeter functor $C^+$ is by definition the composition of reflection functors
    associated to an admissible ordering
    $(a_1, \dots, a_n)$ of $Q$. The action on a filtration, which we also denote by $C^+$, is
    given by $C^+ \flvec{d} := S^+_{a_n} \dots S^+_{a_1} \flvec{d}$. It is not clear
    that $C^+$ on a filtration does not depend on the choice of the admissible ordering.
\end{remark}
We immediately obtain the following.
\begin{cory}
    \label{reflflags:cory:preprojone}
    Let $M$ be a preprojective $K$-representation and
    let $\flvec{d}$ be a filtration of $\dimve M$. Take $r\ge0$ such that
    $(C^+)^r M = 0$.
    
    Then
    $\Fl[Q]{\seqv{\dvec{d}}}{M}$ is non-empty if and only if
    we have that $(C^+)^r \flvec{d} =0$ and that for every intermediate sequence $w$ of admissible sink reflections
    $S^+_w \flvec{d}$ is a filtration of $S^+_w M$. In particular, this depends only
    on the isomorphism class of $M$ and the filtration $\flvec{d}$, but not on the choice
    of $M$ or the field $K$.

    Moreover, if $K$ is a finite field with $q$ elements, then
    $\Fl[Q]{\seqv{\dvec{d}}}{M}$ non-empty implies that
    \[\# \Fl[Q]{\seqv{\dvec{d}}}{M} = 1 \mod q.\]
\end{cory}
\begin{proof}
    Using remark \ref{reflflag:rem:reflfilt} we obtain that for each reflection at a sink
    $a$ of $Q$ we have that
    $S^+_a \flvec{d}$ is again a filtration of $S^+_a M$ if and only if
    $(S^+_a \flvec{d})^{\nu-1} \le \dimve S^+_a M$.
    If this is not the case, then
    the quiver flag is empty by theorem \ref{theom:reflflagmodq}.
    Therefore, if the quiver flag is non-empty, then for every
    intermediate sequence $w$ of admissible sink reflections
    we have that $S^+_w \flvec{d}$ is a filtration of $S^+_w M$. We call this condition (*).

    Assume that (*) holds. Iteratively applying
    theorem \ref{theom:reflflagmodq} we have that
    $\Fl[Q]{\flvec{d}}{M}$ is empty if and only if $\Fl[Q]{(C^+)^r\flvec{d}}{(C^+)^r M} = \Fl[Q]{(C^+)^r\flvec{d}}{0}$
    is empty.
    There is only
    one filtration of the $0$ representation, namely $(0,0,\dots,0)$. Therefore,
    $\Fl[Q]{\flvec{d}}{M}$ is non-empty if and only if $(C^+)^r\flvec{d} = 0$.
    This proves the first part since we already have seen that if (*) does not hold, then
    $\Fl[Q]{\flvec{d}}{M}$ is empty.
    
    Assume now that $K$ is a finite field with $q$ elements. If (*) does not
    hold, then the quiver flag is empty and the claim holds.
    Assume therefore that (*) holds. As before, applying
    theorem \ref{theom:reflflagmodq} yields that
    \[ \#\Fl[Q]{\flvec{d}}{M} = \#\Fl[Q]{(C^+)^r\flvec{d}}{(C^+)^r M} = \#\Fl[Q]{(C^+)^r\flvec{d}}{0} \mod q.\]
    There is only
    one filtration of the zero representation, namely $(0,0,\dots,0)$, and
    the number of flags of this type is obviously equal to one.
    This concludes the proof.
\end{proof}

\section{Dynkin Case}
In this section let $Q$ be
a Dynkin quiver. Then every representation is preprojective, and
there are Hall polynomials with respect to isomorphism classes. We can use the machinery we just developed to prove
that, for $Q$, the generic composition algebra
specialised at $q=0$
and the composition monoid are isomorphic.

\begin{propn}
    \label{reflflags:propn:dynkinallone}
    Let $X$ be a $K$-representation of $Q$ and $w$ a word in vertices of $Q$.
    Then the condition that $X$
    has a filtration of type $w$ only depends on $w$ and $[X]$ and not on the choice of $X$ or the field
    $K$.

    Moreover, we have that
    \[ u_w = \sum_{[X] \in [\mathcal{A}_w]} u_{[X]} \in \mathcal{H}_0(Q). \]
\end{propn}
\begin{proof}
    Since all representations of $Q$ are preprojective, the first part of the
    statement follows directly from corollary \ref{reflflags:cory:preprojone}.
    Therefore, the sum in the second part is well-defined (i.e. the set
    $[\mathcal{A}_w]$ does not depend on the field).

    If $K$ is a finite field with $q$ elements, corollary \ref{reflflags:cory:preprojone} also yields that
    \[ F_w^X = \# \Fl[Q]{\flvec{d}(w)}{X}  = \begin{cases}
        1 \mod q & \text{if } X \in \mathcal{A}_w,\\
        0  & \text{else.}
    \end{cases}
    \]
    Since $F_w^X= f_w^{[X]}(q)$ and we just showed that this is one modulo $q$ for all prime powers $q$
    we have that
    \[f_w^{[X]}(0) =\begin{cases}
        1 & \text{if } [X] \in [\mathcal{A}_w],\\
        0 & \text{else.}
    \end{cases}
    \]
    This yields the claim.
\end{proof}
We obtain the main theorem for the Dynkin case.
\begin{theom}
  The map
  \begin{alignat*}{2}
    \Psi &\colon\quad & \Q\mathcal{M}(Q) & \rightarrow \mathcal{H}_0(Q)\\
    && \mathcal{A} & \mapsto \sum\limits_{[M] \in [\mathcal{A}]} u_{[M]}
  \end{alignat*}
  is an isomorphism of $\Q$-algebras.
  \label{reflflags:theom:dynkin}
\end{theom}
\begin{proof}
    Note that for $Q$ Dynkin we have that $\mathcal{M}(Q) \cong \mathcal{CM}(Q)$
    and $\mathcal{H}_q(Q) \cong \mathcal{C}_q(Q)$. Therefore,
    for each $\mathcal{A} \in \mathcal{M}(Q)$ there is a word $w$ in
    vertices of $Q$ such that $\mathcal{A} = \mathcal{A}_w$. In
    the previous proposition we showed that the map sending $\mathcal{A}_w$ to
    \[\Psi(\mathcal{A}_w) = \sum_{[M] \in [\mathcal{A}_w]} u_{[M]} = u_w\]
    is well-defined. Therefore, $\Psi$ is a homomorphism, since
    \[\Psi(\mathcal{A}_{w} * \mathcal{A}_{v}) = \Psi(\mathcal{A}_{wv}) =
    u_{wv}=u_{w} \diamond u_v = \Psi(\mathcal{A}_{w}) \diamond \Psi(\mathcal{A}_{v}).  \]
    
%
    $\Psi$ is surjective since it is a homomorphism, and the generators
    $u_i$ of $\mathcal{H}_0(Q)$ are in the image of $\Psi$. More precisely,
    $\Psi(\Orbit_{S_i}) = u_i$. Obviously, $\Psi$ is a graded morphism
    of graded algebras. The dimension of the $\dvec{d}$-th graded part
    of $\Q\mathcal{M}(Q)$ is the same as the dimension
    of the $\dvec{d}$-th graded part of $\mathcal{H}_0(Q)$, namely the number of isomorphism classes of representations
    of dimension vector $\dvec{d}$. Since each graded part is finite dimensional and
    $\Psi$ is surjective, we have
    that $\Psi$ is an isomorphism.
\end{proof}

\addtocontents{toc}{\protect\newpage}
\chapter{Extended Dynkin Case}
\dictum[Proverb]{When all you have is a hammer, everything starts to look like a nail.}
\bigskip
In this chapter we examine the relation between the generic composition algebra
and the composition monoid of an extended Dynkin quiver.
In the following, let $Q$ be a connected, acyclic, extended Dynkin quiver.
Fix again a total order $\prec_t$ on the preprojective
and preinjective Schur roots refining the order $\prec$ on $\mathcal{P} \cup \mathcal{I}$.
\label{chap:extdynkin}
\section{Basic Results}
In this section we prove some general facts on the number of points
of quiver flag varieties. First, we want to prove the following.
\begin{theom}
  Let $M$ be any $\F_q$-representation and let $\seqv{\dvec{d}}$ be a filtration of $\dimve M$.
  Then the number of points
  of $\Fl[Q]{\seqv{\dvec{d}}}{M}$ is equal to the number of points of $\Fl[Q]{\seqv{\dvec{d'}}}{M_R}$ modulo $q$ for a
  sequence $\seqv{\dvec{d'}}$ such that the defect of each $\dvec{d'}^i$ is equal to zero. Moreover, $\Fl[Q]{\seqv{\dvec{d}}}{M}$
  is empty if and only if  $\Fl[Q]{\seqv{\dvec{d'}}}{M_R}$ is empty.
\end{theom}
\begin{proof}
By using the Coxeter functors we can
reduce to the case where $M=M_R$. More precisely, eliminate $M_P$ by using $C^+$ $r\ge 0$ times and then apply $C^-$ $r$
times to return to $M$ without
the preprojective part.
In the same way eliminate $M_I$. If, at some point, the reflection of
the filtration $\flvec{d}$ is not a filtration any more, then the flag variety is empty. In this case
set $\flvec{d'}$ to any sequence which is not a filtration of $\dimve M_R$.

Therefore, we can assume that $M$ is purely regular.
If the defect $\partial \dvec{d}^i$ is positive for some $0<i<\nu$, then
$\Fl[Q]{\flvec{d}}{M}$ is empty, since every representation of dimension vector $\dvec{d}^i$ has a preinjective
summand and can therefore not be a subrepresentation of $M$ which is purely regular.
Therefore, $\partial \dvec{d}^i \le 0$ for all $i$.

Now prove the claim by induction on $h=-\sum \partial \dvec{d}^i$.
If $h$=0 we are done. Assume that $h>0$ and $j$ is minimal
with the property that $\partial \dvec{d}^j < 0$. Again, if $\Fl[Q]{\seqv{\dvec{d}}}{M}$ is
non-empty, then there is a representation of
dimension vector $\dvec{d}^j$ having at least one preprojective and possibly some regular summands.
Therefore, the $C^+$-orbit
(here just reflections on a dimension vector) of
$\dvec{d}^j$ tends to $-\infty$ and, applying our machinery of reflections on flags, at some point we have to add
some $r_+ \epsilon_a$, $r_+> 0$,
where $a$ is a source of $Q$, increasing the defect. Therefore, by induction and inverting all reflections, we are done.
\end{proof}
We obtain the following.
\begin{cory}
	Let $M$ be an $\F_q$-representation such that $\dimve M_x \ge \dvec{\delta}$ for at most one $x \in \Proj^1_{\F_q}$
  and let $\seqv{\dvec{d}}$ be a filtration of $\dimve M$. Then $\#\Fl[Q]{\seqv{\dvec{d}}}{M}$ is
  equal to one modulo $q$ if and only if it is non-empty.
  \label{extdyn:onetube_equal_one}
\end{cory}
\begin{proof}
    Using the last theorem we can assume that the defect of each $\dvec{d}^i$ is
    $0$ and $M$ is purely regular with the condition.
    If $U^i$ is a subrepresentation of $M$ of dimension vector $\dvec{d}^i$, then $U^i$
    is purely regular, too. Assume that this is
  not the case, then it must have a non-zero preinjective summand eliminating the defect of a non-zero preprojective
  summand. But this cannot be a subrepresentation of $M$, since $M$ is purely regular and therefore has no preinjective
  subrepresentations. If we decompose $\dvec{d}^i = \sum_{x \in \Proj^1} \dvec{d}^{x, i}$, where
  $\dvec{d}^{x, i} = \dimve U^i_x$,
  then one has that $\ext(\dvec{d}^{x, i}, \dvec{d}^{y, i}) = 0$ for $x \neq y$. Therefore, the
  canonical decomposition of $\dvec{d}^i$ is a refinement of $\sum \dvec{d}^{x, i}$ by lemma \ref{lmm:candecrefine}.
  Since the canonical decomposition
  is unique, and, if $\dvec{\delta}$ appears, the condition on $M$ yields that there is
  only one possibility in which tube it can
  occur, we know that for each subrepresentation $V^i$ of $M$ with $\dimve V^i = \dvec{d}^i$ that
  $\dimve V^i_x = \dvec{d}^{x, i}$. Hence, choosing a flag of type $\seqv{\dvec{d}}$ of $M$ is the
  same as choosing flags of type $\seqv{\dvec{d}^x}$ of $M_x$, i.e.
  \[
  \Fl[Q]{\flvec{d}}{M} = \prod_{x} \Fl[Q]{\flvec{d}^x}{M_x}.
  \]
  By the result on the cyclic quiver the number of points of each
  of these is equal to one modulo $q$, and therefore, so is the number of points of the product.
\end{proof}
The next few lemmas deal with quiver Grassmannians.
\begin{lmm}
  Let $\dvec{d}$ be the dimension vector of a preprojective representation
  and let $M$ be an arbitrary
  $\F_q$-representation
  of dimension $\dvec{d}+\dvec{e}$ for some dimension vector $\dvec{e}$. If
  $\Gr{\dvec{d}}{M}$ is non-empty, then $\#\Gr{\dvec{d}}{M} = 1 \mod q$.

  Dually, if $\dvec{e}$ is the dimension vector of
  an indecomposable preinjective, $M$ an arbitrary $\F_q$-representation of dimension $\dvec{d}+\dvec{e}$ and
  $\Gr{\dvec{d}}{M}$ is non-empty, then $\#\Gr{\dvec{d}}{M} = 1 \mod q$.

  Moreover, in both cases $\Gr{\dvec{d}}{M}$ is non-empty
  if and only if $\Gr{\dvec{d}}{M\tensor F}$ is non-empty, where $F$ is the algebraic
  closure of $\F_q$.  \label{extdyn:preproj_equal_one}
\end{lmm}
\begin{proof}
  Let $w=(a_1, \dots, a_n)$ be an admissible ordering of the vertices of $Q$.
  Recall the map $\sigma$ defined in section \ref{sec:intro:reflfun}.
  For each indecomposable preprojective representation $P$ there is a natural number
  $r=kn+s$ for some $k\ge0$ and $0 \le s < n$ such that
  \[ S^+_{a_s} \dotsm S^+_{a_1}(C^+)^k P =0. \]
  Let $\sigma(P)$ be the minimal such number. The number $\sigma(P)$ depends on the
  choice of the admissible ordering.
  
  The basic idea is now to use $S_a^+$ to reduce $\sigma(P)$ and then use induction.

  Let $\dvec{d}$ be the dimension vector of a preprojective representation $P= \bigoplus P_i$ for
  some indecomposable preprojective representations $P_i$.
  For each $P_i$ let $t_i:=\sigma(P_i)$ for a fixed admissible ordering $(a_1, \dots, a_n)$.

  We proceed by induction on $t=\max \{ t_i \}$.
  If $t = 0$, then we have that each $P_i$ is $0$. Therefore, $P=0$ and $\dvec{d}=0$.
  Obviously, $\# \Gr{\dvec{d}}{M} = \# \Gr{\dvec{0}}{M} = 1$ and
  $\Gr{\dvec{0}}{M}$ is non-empty if and only if $\Gr{\dvec{0}}{M \tensor F}$
  is non-empty.

  Now let $t > 0$. If $\dvec{d} > \dimve M$, then $\Gr{\dvec{d}}{M}$ is empty
  and so is $\Gr{\dvec{d}}{M\tensor F}$. We are done in this case. Assume therefore
  that $\dvec{d} \le \dimve M$.

  Let $a=a_1$. By definition, $\sigma(S_a^+ P_i) = \sigma(P_i)-1$ for
  each $P_i \ncong S_a$ if we
  calculate $\sigma$ with respect to the admissible ordering $(a_2, \dots, a_{n}, a)$
  on $\sigma_a Q$. We also have
  \[\sigma_a \dvec{d} = \sum \sigma_a \dimve P_i = \sum_{P_i \ncong
  S_a} \dimve S_a^+ P_i - s \epsilon_{a}, \]
  where $s$ is the number of $P_i$ isomorphic to $S_{a}$.

  By our algorithm we have that
  \[ \# \Gr{\dvec{d}}{M} = \# \Gr{S_{a}^+\dvec{d}}{S_{a}^+M} \mod q \]
  and the right Grassmannian is non-empty if and only if the
  left Grassmannian is. The same is true for
  $\Gr{\dvec{d}}{M\tensor F}$ and $\Gr{S_{a}^+\dvec{d}}{S_{a}^+M \tensor F}$
  since reflection functors commute with field extension. Therefore,
  it is enough to prove the claim for
  $\Gr{S_{a}^+\dvec{d}}{S_{a}^+M}$.

  We have that
  \[
  S_{a}^+ \dvec{d} = \sigma_{a} \dvec{d} + r_+ \epsilon_a =
  \sum_{P_i \ncong S_{a}} \dimve S_{a}^+ P_i + (r_+- s) \epsilon_{a} \ge 0,\]
  where $r_+ = \max\{0, - (\sigma_{a} \dvec{d})_a\}$. Since
  $\sum_{P_i \ncong S_{a}} \dimve S_{a}^+ P_i \ge 0$, we have that
  $r_+ \le s$. Since $a$ is a source in $\sigma_a Q$, we obtain that
  $S_a$ is a simple injective and $X:=\bigoplus S_a^+ P_i$
  is a representation having $(\dimve X)_a \ge s-r_+$. Therefore, there
  is a surjection $\pi \colon X \rightarrow S_a^{s-r_+}$ and the
  kernel is again preprojective, say isomorphic to
  $\bigoplus P'_i$. Since
  for each $P'_i$ there is at least one $P_j$ such that
  $P'_i \preceq S_a^+ P_j$, we have that $\sigma(P'_i) \le \sigma(S_a^+ P_j)$ by
  lemma \ref{lmm:intro_refl:order}.
Therefore, $\max\{ \sigma(P'_i) \} \le \max\{\sigma(S_a^+ P_i)\} < t$.
  We are done by induction.
%

  If $\dvec{e}$ is the dimension vector of an indecomposable preinjective, then
  $\dvec{e}$ is the dimension vector of an indecomposable preprojective on $Q^{op}$,
  and the claim follows from the preprojective case since
  \[ \Gr[Q]{\dvec{d}}{M} \cong \Gr[Q^{op}]{\dvec{e}}{DM}. \]
\end{proof}
\begin{lmm}
  Let $\dvec{d}$ be a dimension vector having $d_a = 0$ for some $a \in Q_0$ and $M$ an
  $\F_q$-representation of $Q$.
  If $\Gr{\dvec{d}}{M} \neq \emptyset$, then its number is one modulo $q$.

  Moreover, $\Gr{\dvec{d}}{M}$ is non-empty if and only if $\Gr{\dvec{d}}{M \tensor F}$ is, where
  $F$ is the algebraic closure of $\F_q$.
  \label{extdyn:supp_equal_one}
\end{lmm}
\begin{proof}
  Note that there is an admissible series of source reflections
  $S^-_{b_1}\dotsm S^-_{b_r}$ such that each $b_i \neq a$ and
  $a$ is a source in $\sigma_{b_1} \dotsm \sigma_{b_r} Q$.
  Since $b_i \neq a$, we have that $(S^-_{b_1} \dotsm S^-_{b_r} \dvec{d})_a=0$ and neither the number
  of points modulo $q$ nor
  whether it is empty changed. Therefore, we can assume that $a$ is a source in $Q$.

  Let $Q'$ be the quiver obtained from $Q$ by deleting the vertex $a$.
  Subrepresentations $U$ of dimension vector $\dvec{d}$ of $M$ have to make the following diagram commute for
  each $\alpha \colon a \rightarrow i$.
  \[
   \begin{CD}
     0 @>U_\alpha>> U_i \\
     @VVV @VVV\\
     M_a @>M_\alpha>> M_i \\
   \end{CD}
  \]
  Therefore, $U$ is a subrepresentation of $M$ if and only if $U|_{Q'}$ is a subrepresentation
  of $M|_{Q'}$.
   The quiver
  $Q'$ is obviously a union of Dynkin quivers.
  The restriction $M|_{Q'}$ is therefore preprojective, and whether it is empty is independent of the field.
  Hence the claims follow.
\end{proof}

\section{Basis of PBW-Type}
In this section, we will construct a basis of PBW-type for $\mathcal{C}_0(Q)$ consisting
of monomial elements, for an acyclic, extended Dynkin quiver $Q$. This basis then lifts to a basis of
$\mathcal{C}_q(Q)$.
We prove this by showing first that in $\mathcal{C}_0(Q)$
the relations of $\EMon(Q)$ hold and then that the partial normal form in $\EMon(Q)$ is
a normal form in $\mathcal{C}_0(Q)$. This will be the desired basis.

Let $\mathcal{A}$ be a closed irreducible $\GL_\dvec{d}$-stable subvariety of $\Rep{\dvec{d}}$ defined over $\Z$.
For example, take a word $w$ in vertices of $Q$ and consider $\mathcal{A}_w$. Then we define
\[u_\mathcal{A} := \sum_{\alpha \in \mathcal{A}} u_\alpha \in \mathcal{C}_0(Q),\]
where we say that $\alpha \in \mathcal{A}$ for an $\alpha \in \Segre$ if for all
finite fields $K$ and all $X \in \SegreK(\alpha, K)$ we have
that $X \tensor \overline{K} \in \mathcal{A}$, $\overline{K}$ being the algebraic closure of $K$.
Note that for some finite fields $K$ the set $\SegreK(\alpha, K)$ may be empty.

Assume that the vertices $Q_0=\{1,\dots,n\}$
of $Q$ are ordered in such a way that we have $\Ext(S_i, S_j) =0$ for $i\ge j$.
\begin{defn}
    Let $\dvec{d}$ be a dimension vector.
    Define $u_\dvec{d} := u_1^{d_1} \diamond u_2^{d_2} \diamond \dots \diamond u_n^{d_n} \in \mathcal{C}_q(Q)$.
\end{defn}
Then we have in the composition algebra of $Q$:
\[
u_\dvec{d}= \left(\prod_{i\in Q_0} \qnum{d_i}_q!\right) \sum_{\alpha\in \Rep{\dvec{d}}}  \ u_\alpha.
\]
Therefore, specialising $q$ to zero yields
\[
u_\dvec{d} = \sum_{\alpha\in \Rep{\dvec{d}}}  \ u_\alpha = u_{\Rep{\dvec{d}}} \in \mathcal{C}_0(Q).
\]

Let $\seqv{\dvec{d}} = (\dvec{d}^0=0, \dvec{d}^1, \dots , \dvec{d}^\nu)$ be a filtration $\dvec{d}^\nu$. Then
\[
    u_{\dvec{d}^\nu- \dvec{d}^{\nu-1}} \diamond \dots \diamond u_{\dvec{d}^2- \dvec{d}^1}
    \diamond u_{\dvec{d}^1- \dvec{d}^0}
    = \prod\qnum{d^i_j-d^{i-1}_j}_q! \sum_{\alpha\in\mathcal{A}_\flvec{d}} f_{\seqv{\dvec{d}}}^{\alpha}(q) u_\alpha
    =\sum_{\alpha \in \mathcal{A}_\flvec{d}} f_{\seqv{\dvec{d}}}^{\alpha}(0) u_\alpha \in \mathcal{C}_0(Q)
\]
for some polynomials $f_{\seqv{\dvec{d}}}^{\alpha}$ for each $\alpha \in \mathcal{A}_\flvec{d}$ such
that for each finite field $\F_q$ and each $X \in \SegreK(\alpha, \F_q)$
we have that $\# \Fl[Q]{\flvec{d}}{X} = f_\flvec{d}^\alpha(q)$. In the
sum only $\alpha \in \mathcal{A}_\flvec{d}$ appear since if the polynomial is non-zero, then the
flag variety will be non-empty over the algebraic closure. On the other hand,
if $\alpha \in \mathcal{A}_{\flvec{d}}$ we do not necessarily have that $f^\alpha_\flvec{d} \neq 0$.

In this section we will often deal with quiver Grassmannians. Therefore, we define
\[ f^\alpha_{\dvec{e}\ \dvec{d}}:= f^\alpha_{(0, \dvec{d}, \dvec{d} + \dvec{e})}.\]

We want to say something about $f(0)$ for a Hall polynomial $f$. For this we often use the following.
\begin{lmm}
	Let $f \in \Q[x]$ be a polynomial. Assume that there is a $c \in \Z$ such that
    $f(q) \in \Z$ and $f(q) = c \mod q$, for infinitely many integers
	$q \in \Z$. Then $f(0) = c$.
	\label{lmm:extdyn:modq_eqconst}
\end{lmm}
\begin{proof}
    We have that $f = \frac{1}{N}g$ for a polynomial $g \in \Z[x]$ and a positive integer $N \in \N$.
	Therefore,
	\[ c \cdot N= f(q) \cdot N = g(q) = g(0) \mod q
	\]
	for infinitely many integers $q$. Hence, $g(0) = c \cdot N$ and therefore
    $f(0) = \frac{1}{N} g(0) = c$.
\end{proof}

In the composition algebra $\mathcal{C}_0(Q)$ we consider the following elements:
\begin{itemize}
  \item $u_\dvec{d}$ for $\dvec{d}$ a Schur root,
  \item $u_{s\dvec{\delta}}$ for $s > 0$.
\end{itemize}

Note that in the composition monoid we have that
\[\rep{s\dvec{\delta}} * \rep{t\dvec{\delta}} = \rep{(s+t) \dvec{\delta}} = \rep{t \dvec{\delta}} * \rep{s \dvec{\delta}}.\]
In the composition algebra this will not be true anymore, but at least we will prove that one has
\[u_{s\dvec{\delta}} \diamond u_{t\dvec{\delta}} = u_{t\dvec{\delta}} \diamond u_{s\dvec{\delta}}.\]

In the following we will denote by $\ext(\dvec{d}, \dvec{e}) := \ext_{FQ}(\dvec{d}, \dvec{e})$, where
$F$ is any algebraically closed field. This number is independent of $F$ by work
of A. Schofield and W. Crawley-Boevey.

We will first prove that all relations of $\EMon(Q)$ hold in $\mathcal{C}_0(Q)$,
replacing $\SRoot{s\dvec{d}}$ by $u_{s\dvec{d}}$. In the following all elements
$u_\dvec{d}$ will be in $\mathcal{C}_0(Q)$ if not otherwise stated.

%
\begin{lmm}
   Let $\dvec{d}$ be a preprojective Schur root, $\dvec{e}$ an arbitrary
   dimension vector such that $\ext(\dvec{d}, \dvec{e}) = 0$ and $r, s \in \N$.
   Then
   \[
   u_{s\dvec{e}} \diamond u_{r\dvec{d}} = u_{r\dvec{d} + s\dvec{e}}.
   \]

   Dually, let $\dvec{e}$ be a preinjective Schur root,
   $\dvec{d}$ an arbitrary dimension vector such that
   $\ext(\dvec{d}, \dvec{e}) = 0$ and $r, s \in \N$.
   Then
   \[
   u_{s\dvec{e}} \diamond u_{r\dvec{d}} = u_{r\dvec{d} + s\dvec{e}}.
   \]
   \label{lmm:extdyn:preproj_full_dimvec}
\end{lmm}
\begin{proof}
    Since $\ext(r\dvec{d}, s\dvec{e}) = 0$, we have, for any
    $\F_q$-representation $M$ of dimension vector
    $r\dvec{d} + s\dvec{e}$, that $\Gr{r\dvec{d}}{M\tensor F}$ is non-empty,
    where $F$ is the algebraic closure of $\F_q$.

    By lemma \ref{extdyn:preproj_equal_one} we have that
    $f^\alpha_{s\dvec{e}\ r\dvec{d}}(q) = 1 \mod q$
    for all prime powers $q$ and all $\alpha \in \Rep{r\dvec{d}+s\dvec{e}}$.
    Therefore, $f^\alpha_{s\dvec{e}\ r\dvec{d}}(0) =1$ and the claim follows.
\end{proof}
\begin{lmm}
   Let $\dvec{d}$ be a preinjective Schur root, $\dvec{e}$ a preprojective Schur
   root such that $\ext(\dvec{d}, \dvec{e}) = 0$ and $r, s \in \N$.
   Then
   \[
   u_{s\dvec{e}} \diamond u_{r\dvec{d}} = u_{r\dvec{d} + s\dvec{e}}.
   \]
   \label{lmm:extdyn:preinj_preproj_full_dimvec}
\end{lmm}
\begin{proof}
      We first prove the claim for $s=t=1$.

      Since there is an indecomposable preprojective representation
      $P$ of dimension vector $\dvec{e}$
      and an indecomposable preinjective representation $I$ of dimension vector $\dvec{d}$, the minimal value
      of $\Ext$ is taken on those. Moreover, since $\Hom(I, P) = 0$, we have that
      \[
      0 = \ext(\dvec{d}, \dvec{e}) = \Ext(I, P) = \Ext(I, P) - \Hom(I, P) = -\bform{\dvec{d}}{\dvec{e}}.
      \]
      Therefore, $\dvec{d} < \dvec{\delta}$ since, otherwise,
      $\dvec{d} = \dvec{d'} + \dvec{\delta}$ with $\dvec{d'}$ being a positive preinjective root and
      \[\bform{\dvec{d}}{\dvec{e}} = \bform{\dvec{d'}}{\dvec{e}} + \partial{\dvec{e}} < 0,\]
      in contradiction to $\bform{\dvec{d}}{\dvec{e}} = 0$. Dually, we have that $\dvec{e} < \dvec{\delta}$.

      Summing up, we have that $\dvec{d} + \dvec{e} < 2 \dvec{\delta}$. Therefore, for each
      $\F_q$-representation $M$ of dimension
      vector $\dvec{d} + \dvec{e}$ we have that $\dimve M_x \ge \dvec{\delta}$ for at most one
      $x \in \Proj^1_{\F_q}$ and for each tube $\mathcal{T}_x$ with
      $M_x \neq 0$ we have that $\deg x = 1$. Over the algebraic closure $F$ of
      $\F_q$ we have that every representation of dimension vector $\dvec{d}+\dvec{e}$ has a subrepresentation of
      dimension vector $\dvec{d}$. Let $\alpha$ be the decomposition symbol of $M$.
      Since $M$ lives only in tubes of degree one by the restriction on the dimension,
      we have that the decomposition symbol
      of $M \tensor F$ is also $\alpha$. Therefore, the polynomial $f^\alpha_{\dvec{e}\ \dvec{d}}$
      is non-zero. We conclude that
      $f^\alpha_{\dvec{e}\ \dvec{d}}(q)$
      modulo $q$ is equal to one for an infinite number of prime powers $q$, namely all $q$ which
      are not a zero of $f^\alpha_{\dvec{e}\ \dvec{d}}$, by
      corollary \ref{extdyn:onetube_equal_one}. Hence,
      the constant coefficient is equal to one and the claim follows.

      Therefore, we have proved that $u_\dvec{e} \diamond u_\dvec{d} = u_{\dvec{d} + \dvec{e}}$.
      We automatically have that $\ext(\dvec{e}, \dvec{d})=0$, since $\dvec{d}$ is
      preinjective and $\dvec{e}$ is preprojective. By the previous lemma, we have
      that $u_\dvec{d} \diamond u_\dvec{e} = u_{\dvec{d}+\dvec{e}}$. Therefore,
      $u_\dvec{d} \diamond u_{\dvec{e}} = u_\dvec{e} \diamond u_\dvec{d}$.
      Again by the previous lemma, we have that $(u_\dvec{d})^r = u_{r \dvec{d}}$ and
      $(u_{\dvec{e}})^r = u_{r\dvec{e}}$. Using the lemma once more, we obtain
      \[
      u_{s\dvec{e}} \diamond u_{r\dvec{d}} = u_{\dvec{e}}^s \diamond u_{\dvec{d}}^r =
      u_{\dvec{d}}^r \diamond u_{\dvec{e}}^s = u_{r\dvec{d}} \diamond u_{s\dvec{e}} =
      u_{r\dvec{d}+s\dvec{e}}.
      \]
\end{proof}
\begin{lmm}
   Let $\dvec{d}$, $\dvec{e}$ be
   dimension vectors, one being a real Schur root of defect $0$, such that
   $\ext(\dvec{d}, \dvec{e}) = 0$ and let $r, s \in \N$.
   Then
   \[
   u_{s\dvec{e}} \diamond u_{r\dvec{d}} = u_{r\dvec{d} + s\dvec{e}}.
   \]
   \label{lmm:extdyn:realreg_full_dimvec}
\end{lmm}
\begin{proof}
      First assume that $\dvec{d}$ is a real Schur root and has defect $0$. Then
      $\dvec{d} < \dvec{\delta}$ and there is a regular simple representation $T$ of dimension vector
      $\dvec{f}$ such that $\dvec{d} + \dvec{f} \le \dvec{\delta}$. There is a series of sink
      reflections $S^+_{a_1}\dotsm S^+_{a_k}$ such that
      $(S^+_{a_1}\dotsm S^+_{a_k}(r\dvec{f}))_j = r\delta_j$ for an extending vertex $j$ of $Q$.
      Therefore, $(S^+_{a_1}\dotsm S^+_{a_k} (r\dvec{d}))_j = 0$.

      Let $M$ be an $\F_q$-representation of $Q$ of dimension vector
      $r\dvec{d}+s\dvec{e}$ and let $F$ be the algebraic closure of $\F_q$.
      We have that $M\tensor F$ has a subrepresentation of dimension vector $r\dvec{d}$, since
      $\ext(\dvec{d}, \dvec{e}) = 0$. Therefore, for each $j\ge 1$
      \[
      \dimve S^+_{a_j}\dotsm S^+_{a_k} (M\tensor F) \ge S^+_{a_j}\dotsm S^+_{a_k} (r\dvec{d}).
      \]
      Since
      \[
      \dimve S^+_{a_j}\dotsm S^+_{a_k} (M\tensor F) = \dimve S^+_{a_j}\dotsm S^+_{a_k} M,
      \] the same is true for $M$. Applying lemma
      \ref{extdyn:supp_equal_one} yields that the number modulo $q$ is one since
      it is non-empty over the algebraic closure.

      If $\dvec{e}$ is a real Schur root with defect zero,
      then apply $\hat D$ to $\Gr{r\dvec{d}}{M}$ for $M \in \Rep{r\dvec{d} + s\dvec{e}}$
      to reduce to first case.
\end{proof}

\begin{lmm}
  For $s, t \in \N$ we have that
  \[ u_{s \dvec{\delta}} \diamond u_{t \dvec{\delta}} = u_{t \dvec{\delta}} \diamond u_{s \dvec{\delta}}
  \in \mathcal{C}_0(Q). \]
   \label{lmm:extdyn:imag_commute}
\end{lmm}
\begin{proof}
	We need to show that for all finite fields $\F_q$ and each $M \in \RepK{(s+t) \dvec{\delta}}{\F_q}$ we have that
  $\#\Gr{t \dvec{\delta}}{M} = \#\Gr{s \dvec{\delta}}{M} \mod q$. Since both Grassmannians
  are non-empty over the algebraic closure, any sequence of reflections applied to $s\dvec{\delta}$
  or $t\dvec{\delta}$ will
  yield a filtration of the reflection of $M$.

  Assume that
  $M$ has a preprojective summand. We want to show that
  \[
  \#\Gr{t \dvec{\delta}}{M} = \#\Gr{s \dvec{\delta}}{M} =1 \mod q.
  \]
  Since the defect of $M$ is zero and it has a preprojective summand, it will automatically have a
  preinjective summand. Let $S^+_{a_1} \dotsm S^+_{a_k}$ be a minimal sequence of admissible sink reflections
  eliminating a preprojective summand of $M$. Let $M':= \dimve S^+_{a_1} \dotsm S^+_{a_k}M$
  and $Q' := \sigma_{a_1} \dotsm \sigma_{a_k} Q$. Then,
  at this point, we have that
  $\dimve M' = \dimve M + r \epsilon_{a_1}$, where $r>0$ is the number of times
  $M$ has the eliminated preprojective as a summand. We have that
  $\sigma_{a_1} \dotsm \sigma_{a_k} t \dvec{\delta} = t \dvec{\delta}$,
  and therefore the cokernel of $M'$ by a subrepresentation of dimension vector $t \dvec{\delta}$
  has dimension vector
  $\dim M'  - t\dvec{\delta} = s \dvec{\delta} + r \epsilon_{a_1}$.
  Moreover, the vertex $a_1$ is a source in $Q'$. If we apply $\hat D$, we
  obtain that
  \[
  \#\Gr[Q']{t \dvec{\delta}}{M'}=\#\Gr[Q'^{op}]{s \dvec{\delta} + r \epsilon_{a_1}}{DM'}.
  \]
  Note that $a_1$ is a sink in $Q'^{op}$. Applying sink reflections we obtain that each entry of $(C^+)^i \epsilon_{a_1}$
  tends to
  $-\infty$ for $i \rightarrow \infty$. Therefore, at some point while using the
  algorithm of lemma \ref{theom:reflflagmodq} to do sink reflections,
  we will have that one component of the reflected dimension vector
  will become zero. But then,
  applying lemma \ref{extdyn:supp_equal_one} yields that $\#\Gr{t \dvec{\delta}}{M} = 1 \mod q$ since
  the Grassmannian is non-empty over the algebraic closure of $\F_q$.
  The same argument yields that $\#\Gr{s \dvec{\delta}}{M} = 1 \mod q$.

  Assume now that $M$ is purely regular. Consider a subrepresentation $U$ of $M$ of dimension vector
  $t\dvec{\delta}$. Then $U=U_R$. Since $\Ext(U_x, U_y) = 0$ for $x \neq y \in \Proj^1_{\F_q}$
  we have that
  $\ext(\dimve U_x, \dimve U_y) = 0$. Therefore, the canonical decomposition of
  $t \dvec{\delta}$ is a refinement of $\sum_{x\in\Proj^1} \dimve U_x$ by
  lemma \ref{lmm:candecrefine}. But the
  canonical decomposition of $t \dvec{\delta}$ is $t$ times $\dvec{\delta}$. Therefore, there
  are integers $t_x \in \N$ such that
  $\dimve U_x = t_x (\deg x)\dvec{\delta}$ with 
  $\sum (\deg x) t_x = t$. The regular representation $M$ has its regular
  part only in finitely many tubes, say $\mathcal{T}_{x_1}, \dots, \mathcal{T}_{x_k}$. By a similar argument using the
  canonical decomposition
  we have that $\dimve M_{x_i} = m_i (\deg x_i) \dvec{\delta}$ for some integers $m_i > 0$.
  Hence,
  \[\#\Gr{t \dvec{\delta}}{M} = \sum_{\substack{(t_1,\dots,t_k)\\ \sum t_i(\deg x_i)=t}}
  \prod_{i=1}^k \#\Gr{t_i(\deg x_i) \dvec{\delta}}{M_{x_i}}\overset{\mod q}{=}
  \sum_{\mathclap{\substack{(t_1,\dots,t_k)\\ \sum t_i(\deg x_i)=t\\
  t_i \le m_i}}} 1. \]
  It is easy to see that
  $M_{x_i}$ has a subrepresentation of dimension
  vector $t_i (\deg x_i) \dvec{\delta}$ if and only if $t_i \le m_i$. The last
  equality follows by
  theorem \ref{coeff-compmonoid-areone} and using that $\mathcal{T}_x$ is equivalent
  to $\Mod \kappa(x) C_{r-1}$, $r$ being the rank of $\mathcal{T}_x$. Note that $\kappa(x)$
  has $q^{\deg x}$ elements.

  Repeating the same argument for $s \dvec{\delta}$, we obtain that
  \[\#\Gr{s \dvec{\delta}}{M} = \sum_{\substack{(s_1,\dots,s_k)\\ \sum s_i(\deg x_i)=s}}
  \prod_{i=1}^k \#\Gr{s_i(\deg x_i) \dvec{\delta}}{M_{x_i}}\overset{\mod q}{=}
  \sum_{\mathclap{\substack{(s_1,\dots,s_k)\\ \sum s_i(\deg x_i)=s\\
  s_i \le m_i}}} 1. \]
  There is an obvious bijection between the two sets we are summing over, namely
  sending $(t_1, \dots, t_k) \mapsto (m_1 - t_1, \dots, m_k - t_k)$, and
  therefore the two numbers agree.
\end{proof}
We can now use the previous lemmas to prove the following.
\begin{theom}
  The following relations hold in $\mathcal{C}_0(Q)$.
  \begin{align*}
  u_{s\dvec{e}} \diamond u_{r\dvec{d}}&= u_{r\dvec{d}} \diamond u_{s\dvec{e}}&&
  \text{for all }
  s,\ t \ge 0,\ \dvec{d},\ \dvec{e} \text{ Schur roots such that}\\
      &&&\ext(\dvec{d}, \dvec{e}) = \ext(\dvec{e}, \dvec{d})= 0;\\
  u_{s\dvec{e}} \diamond u_{r\dvec{d}} &= u_{r\dvec{d} +s\dvec{e}}&&
  \text{for all }
  s,\ t \ge 0,\ \dvec{d},\ \dvec{e} \text{ Schur roots such that}\\
  &&&\ext(\dvec{d}, \dvec{e}) = 0 \text{ and not both roots are imaginary};\\
  (u_\dvec{d})^r &= u_{r\dvec{d}}&& \text{ for all } r >0,\ \dvec{d} \text{ a real Schur root}.
  \end{align*}
  \label{theom:extdyn:schurmult_commutes}
\end{theom}
\begin{proof}
    If both roots are imaginary, this is lemma \ref{lmm:extdyn:imag_commute}. In the
    other cases at least one root is real, and it is enough to prove the second statement
    since then the first statement follows by
  \[ u_{s\dvec{e}} \diamond u_{r\dvec{d}} = u_{r\dvec{d} +s\dvec{e}}=
  u_{r\dvec{d}} \diamond u_{s\dvec{e}}\in \mathcal{C}_0(Q). \]
  The third statement follows since for a real Schur root
  $\dvec{d}$ one always has $\ext(\dvec{d}, \dvec{d}) = 0$.

  Now we prove the second statement. If $\dvec{d} = \dvec{\delta}$, then
  $\dvec{e}$ is not preprojective. Therefore, lemmas \ref{lmm:extdyn:preproj_full_dimvec}
  and \ref{lmm:extdyn:realreg_full_dimvec} yield the claim. If $\dvec{e} = \dvec{\delta}$,
  then $\dvec{d}$ is not preinjective, and again lemmas \ref{lmm:extdyn:preproj_full_dimvec}
  and \ref{lmm:extdyn:realreg_full_dimvec} yield the claim. If one
  of them is real of defect zero, then apply lemma \ref{lmm:extdyn:realreg_full_dimvec}.
  Finally, if both are real not of defect zero, then lemmas \ref{lmm:extdyn:preproj_full_dimvec}
  and \ref{lmm:extdyn:preinj_preproj_full_dimvec} finish the proof.
\end{proof}
Now we obtain the following.
\begin{theom}
    Let $\dvec{d}$ and $\dvec{e}$ be Schur roots, at least one real, such that
    $\ext(\dvec{d}, \dvec{e}) = 0$ and $r, s \in \N$. Let $\sum_{i=1}^k t_i\dvec{f}^i$ be the canonical
  decomposition of $r\dvec{d}+s\dvec{e}$, i.e. the $\dvec{f}^i$ are pairwise different Schur roots
  with $\ext(\dvec{f}^i, \dvec{f}^j) =0$ and $t_i >0$.

  Then
  \[ u_{s\dvec{e}} \diamond u_{r\dvec{d}} = u_{t_1\dvec{f}^1} \diamond \dots \diamond u_{t_k\dvec{f}^k}. \]
  In particular, the product on the right hand side does not depend on the order.
\end{theom}
\begin{proof}
  We already proved in theorem \ref{theom:extdyn:schurmult_commutes} that
  \[ u_{s\dvec{e}} \diamond u_{r\dvec{d}} = u_{s\dvec{d}+r\dvec{e}}. \]
  We have to show that the right hand side equals the same.

  By the previous theorem, the product on the right hand side does not depend on the order.
  We can therefore
  assume that we have the preinjective roots on the left and the preprojective roots on the right.

  At most one $\dvec{f}^i$ is equal to the isotropic Schur root $\dvec{\delta}$.
  We start multiplying the expression together, starting with $u_{t\dvec{\delta}}$. Then the left
  factor will be a multiple of a preinjective
  or a real regular Schur root, or the right factor will be a multiple of a preprojective
  or a real regular Schur root. The claim then follows by lemmas \ref{lmm:extdyn:preproj_full_dimvec}
  and \ref{lmm:extdyn:realreg_full_dimvec}.
\end{proof}
As the last step before proving the main theorem we need to cope with Schur roots living in only one
inhomogeneous tube.
\begin{lmm}
  Let $\dvec{c}^1, \dots, \dvec{c}^k$ be real Schur roots such that the general representations of $\dvec{c}^i$
  live in a single inhomogeneous tube, say $\mathcal{T}_x$. Let $M$ be such that
  \[ \overline{\Orbit_M} = \rep{\dvec{c}^1} * \dots * \rep{\dvec{c}^k}. \]
  This representation exists and we have that
  \[ u_{\dvec{c}^1} \diamond \dots \diamond u_{\dvec{c}^k} = u_{\overline{\Orbit_M}}. \]
\end{lmm}
\begin{proof}
    We can assume that each $\dvec{c}^k$ is the dimension vector of a regular simple representation
    in $\mathcal{T}_x$ by
  lemma \ref{lmm:extdyn:realreg_full_dimvec}.
  We already showed in the part on the composition monoid that there is a generic extension
  $M$ with this property living in the inhomogeneous tube $\mathcal{T}_x$. Therefore,
  $M$ is given by a decomposition symbol $\alpha = (\mu, \emptyset)$. A representation over
  an algebraically closed field $F$ has a
  filtration of type $\dvec{c}^1, \dots, \dvec{c}^k$ if and only if it is a degeneration of $M$.
  Therefore, we only have
  to consider decomposition symbols $\beta \in \overline{\Orbit_M}$.

We use the same method as in lemma \ref{extdyn:preproj_equal_one}, but now
we deal with a flag and not only with a Grassmannian.
In order to make the proof more readable, for
a flag $\flvec{d}$ we set \[\flvec{a}=(\dvec{a}^1, \dots, \dvec{a}^{\nu-1}) :=
  (\dvec{d}^1 - 0, \dvec{d}^2-\dvec{d}^1, \dots, \dvec{d}^{\nu-1} - \dvec{d}^{\nu-2})\] and define
  $\HFl{\flvec{a}}{M}:= \Fl{\flvec{d}}{M}$. Note that we do not loose any information since
  $\dvec{d}^\nu = \dimve M$. Moreover, $S^+_a$ applied to a flag in the new notation
  is given by
  \[
  S_a^+ \flvec{a} = (\sigma_a \dvec{a}^1 + r_+^1 \epsilon_a, \sigma_a \dvec{a}^2 + (r_+^2 - r_+^1)\epsilon_a, \dots,
  \sigma_a \dvec{a}^{\nu-1} +(r_+^{\nu -1} -r_+^{\nu -2})\epsilon_a).
  \]

   We show that for any decomposition symbol $\beta \in \overline{\Orbit_M}$ and any
   $\F_q$-representation $N \in \SegreK(\beta, \F_q)$ we have that
  \[ \#\HFl{(\dvec{c}^k, \dvec{c}^{k-1}, \dots, \dvec{c}^2)}{N} =
  1 \mod q. \]
   This then yields the claim.

  Now let $N$ be such an $\F_q$-representation.
  First, apply sink reflections in the admissible ordering to the flag until the reflected
  $N$ has no preprojective summand left. Note that $\sigma_w \dvec{c}^i \ge 0$ for
  each admissible word in vertices of $Q$. Therefore, in each step $\seqv{r}_+=0$.
  This means that we end
  up with a flag of type $(\sigma_w \dvec{c}^k, \dots, \sigma_w \dvec{c}^2)$ and $\dimve S^+_w N =
  \sigma_w \dimve N + \dimve I$, where $I$ is some preinjective representation. Applying $\hat D$, we
  have to consider a flag of the following type of $D S^+_w N$:
  \[ (\sigma_w\dvec{c}^1+\dimve I, \sigma_w\dvec{c}^2, \dots, \sigma_w\dvec{c}^{k-1}). \]
  $D S^+_w N$ has no preinjective summands and $\dimve I$ is the dimension vector of a preprojective representation
  of $\sigma_w Q^{op}$.
  Note that $D S^+_w N \tensor \overline{\F_q}$ has a flag of the given type, since this is invariant under
  under $\hat D$ and, by theorem \ref{theom:reflflagmodq},
  sink reflections.
  We are left with showing that in this situation
  \[ \#\HFl{(\sigma_w\dvec{c}^1+\dimve I, \sigma_w\dvec{c}^2, \dots, \sigma_w\dvec{c}^{k-1})}{D S^+_w N} = 1 \mod q.\]

  The flag we have at this point has the following properties, formulated
  for a filtration $\flvec{a}$ and an $\F_q$-representation $N'$:
  \[ \HFl{(\dvec{a}^1, \dvec{a}^{2}, \dots, \dvec{a}^{k-1})}{N'\tensor \overline{\F_q}}\]
  is non-empty, $N'$ has no preinjective direct summand,
  each
  $\dvec{a}^i = \dimve R_i + \dimve P_i$ where $R_i$ is $0$ or a regular simple representation
  living in one inhomogeneous tube $\mathcal{T}_x$ and $P_i$ is $0$ or a
  preprojective representation.
  If we can show that, in this situation, we have
  \[ \#\HFl{(\dvec{a}^1, \dvec{a}^{2}, \dots, \dvec{a}^{k-1})}{N'} = 1 \mod q,\]
  then the claim follows.
  We want to prove this more general statement by induction.

  For a preprojective representation $P$ we define
  \[\sigma (P) := \max\{\sigma(P_1), \dots, \sigma(P_r)\},\]
  where $P_1, \dots, P_r$ are the indecomposable direct summands of $P$ and $\sigma(0) := 0$.

  Let $r$ be the number of $R_i$ not equal to $0$
  and $t:=\max\{ \sigma(P_i)\}$ for
  a fixed admissible ordering $(b_1,\dots, b_n)$ of the vertices of $Q$.
  To prove the claim we do induction on $(r,t)$, ordered lexicographically, i.e.
  $(r,t) < (r',t')$ if $r<r'$ or $r=r'$ and $t<t'$.

  We start the induction at $t=0$. In this case $\dimve P_i = 0$ for all $i$,
  therefore we can proceed as in the beginning and eliminate
  the preprojective summand of $N'$ by doing sink reflections. While doing this, $t$ remains $0$ and
  $r$ constant.
  We can therefore assume that $N'$ is purely regular since $N'$ did not have a preinjective summand.

  The only regular representation of dimension vector $\dimve R_i$ is $R_i$ itself,
  living in the inhomogeneous tube $\mathcal{T}_x$, and therefore all
  regular representations having a flag of type $(\dimve R_1, \dimve R_2, \dots, \dimve R_{k-1})$
  of dimension vector $\sum \dimve R_i$
  live in $\mathcal{T}_x$, too. Since the tubes are orthogonal, we have that
  \[\#\HFl{(\dimve R_1, \dimve R_2, \dots, \dimve R_{k-1})}{N'} =
\#\HFl{(\dimve R_1, \dimve R_2, \dots, \dimve R_{k-1})}{N'_x}.\]
  The representation $N'_x$ lives in an inhomogeneous tube which is equivalent to representations
  of the cyclic quiver $C_n$.
  By lemma \ref{cyclic:lmm:filtrindept},
  we have that
  \[
  \HFl{(\dimve R_1, \dimve R_2, \dots, \dimve R_{k-1})}{N'_x} \neq \emptyset.
  \]
  By the result on the cyclic quiver we finally conclude that
  \[
\#\HFl{(\dimve R_1, \dimve R_2, \dots, \dimve R_{k-1})}{N'_x} = 1 \mod q.
\]

  Now let $t>0$.
  The vertex $b=b_1$ is a sink of $Q$ and $\sigma_b \dvec{a}^i = \sigma_b \dimve R_i + \sigma_b \dimve P_i$.
  We know that $\sigma_b \dimve R_i = \dimve S_b^+ R_i \ge 0$. Let us assume that $P_i$
  has the simple $S_b$ $s_i$-times
  as a direct summand.
  Then $\sigma_b \dimve P_i = \dimve S^+_b P_i - s_i\epsilon_b$.
  Therefore, $r_+^i - r_+^{i-1} \le s_i$. Note that $\sigma(S^+_b P_i) < \sigma(P_i)$ calculated with respect
  to the admissible ordering $(b_2, \dots, b_n, b)$ of $\sigma_b Q$.
  Let $u, v \in \N$ such that $0\le s_i-r_+^i+r_+^{i-1}=u+v$, $\sigma_b \dimve R_i - u \epsilon_b \ge 0$ and
  $\dimve S^+_b P_i  - v \epsilon_b \ge 0$. Since $S_b$ is simple injective in $\sigma_b Q$, there is
  a surjection from the regular simple representation $S_b^+ R_i$ to $S_b^u$. If $u>0$, then the
  kernel will be preprojective since $S_b^+ R_i$ is regular simple, and therefore there is a preprojective
  representation $X_i$ such that $\dimve X_i = \dimve S_b^+ R_i - u\epsilon_b$. If $u=0$, set
  $X_i := R_i$.

  With a similar argument we
  have that there is a preprojective representation $P_i'$ such that $\dimve P_i' = \dimve S_b^+P_i -v \epsilon_b$.
  By lemma \ref{lmm:intro_refl:order} we also have that $\sigma(P_i') \le \sigma(S_b^+P_i) < \sigma(P_i)$.
Therefore,
\[
S_b^+ \flvec{a} = ( \dimve X_1 + \dimve P'_1, \dimve X_2 + \dimve P'_2, \dots,
\dimve X_{k-1} + \dimve P'_{k-1})
\]
is of the same form as $\flvec{a}$ and if $r$ did not decrease, then $t$ did. Therefore, we are done by induction.
%
%
\end{proof}
\begin{propn}
    There is a surjection of $\Q$-algebras
    \begin{alignat*}{2}
        \Xi &\colon\quad& \Q\EMon(Q) &\rightarrow \mathcal{C}_0(Q)\\
        && \SRoot{s\dvec{d}} & \mapsto u_{s\dvec{d}}.
    \end{alignat*}
\end{propn}
\begin{proof}
    We just showed that the defining relations (\ref{eq:rootrel1}), (\ref{eq:rootrel2}) and (\ref{eq:rootrel3})
    of $\EMon(Q)$ are
    satisfied by the elements $u_{s\dvec{d}}$ in $\mathcal{C}_0(Q)$.
    Therefore, the given map is well-defined. It is surjective since
    the generators $u_i$ are in the image of $\Xi$.
\end{proof}
Now we can finally prove the main result.
\begin{theom}
  Let $Q$ be a connected, acyclic, extended Dynkin quiver
  with $r$ inhomogeneous tubes indexed by $x_1, \dots, x_r \in \Proj^1_\Z$.
  Let $w$ be a word in vertices of $Q$. Then $u_w$ can be uniquely written as
  $\mathcal{P}\mathcal{C}_1\dotsm \mathcal{C}_r\mathcal{R}\mathcal{I} \in \mathcal{C}_0(Q)$, where
    \begin{align*}
        \mathcal{P}&=u_{s_1\dvec{p}^1} \diamond \dots \diamond u_{s_k\dvec{p}^k} && \text{with } s_i > 0,\ \dvec{p}^i
        \text{ preprojective Schur roots and }\\
        &&&\dvec{p}^i \prec_t \dvec{p}^j \text{ for all } i < j;\\
        \mathcal{C}_i&=u_{\overline{\Orbit_{M_i}}} && \text{for a separated } [M_i] \in [\mathcal{T}_{x_i}];\\
        \mathcal{R}&=u_{\lambda_1\dvec{\delta}} \diamond \dots \diamond u_{\lambda_l\dvec{\delta}} && \text{with }
        \lambda=(\lambda_1\ge \dots \ge \lambda_l)
        \text{ a partition};\\
        \mathcal{I}&=u_{t_1\dvec{q}^1} \diamond \dots \diamond u_{t_m\dvec{q}^m} && \text{with } t_i > 0,\ \dvec{q}^i
        \text{ preinjective Schur roots and }\\
        &&&\dvec{q}^i \prec_t \dvec{q}^j \text{ for all } i < j.
    \end{align*}
  Moreover, the set of elements of this form gives a basis of $\mathcal{C}_0(Q)$ and
  $\Xi\colon \Q\EMon(Q) \rightarrow \mathcal{C}_0(Q)$ is an isomorphism.
\end{theom}
\begin{proof}
    By lemma \ref{compmon:lmm:normformemon} every element of $\EMon(Q)$ can be written as above.
    Since $\Xi$ is a homomorphism we have that every monomial element $u_w \in \mathcal{C}_0(Q)$ can
    be written in this form.

  The composition algebra at $q=0$, $\mathcal{C}_0(Q)$,
  is naturally graded by dimension vector and each graded part has the same dimension
  as the corresponding graded part of the positive part $U^+(\mathfrak{g})$
  of the universal enveloping algebra of the Kac-Moody
  Lie algebra $\mathfrak{g}$ given by the Cartan datum associated to $Q$.
  Therefore, $\dim\mathcal{C}_0(Q)_\dvec{d}$ is the number of ways of
  writing $\dvec{d}$ as a sum of positive roots with multiplicities, in the sense that we
  count $m \dvec{\delta}$ for $m > 0$ with multiplicity $n-1$, where $n$ is the number of vertices of $Q$, and every
  other root with multiplicity one.

  For each inhomogeneous tube $\mathcal{T}_x$ of rank $l$ and each dimension vector $\dvec{d}$
  we have that $\dim \mathcal{C}_0(\mathcal{T}_x)_\dvec{d}$
  is equal to
  $\dim U^+(\hat{\mathfrak{sl}}_l)_\dvec{d}$,
  the dimension of the $\dvec{d}$-th graded part of the positive part of the universal enveloping algebra of
  the Lie algebra $\hat{\mathfrak{sl}}_l$, by \cite{Ringel_compalgofcyclic}.
  The roots of $\hat{\mathfrak{sl}}_l$ correspond to the dimension vectors of indecomposables
  in $\mathcal{T}_x$. Moreover, the multiplicity of $m \dvec{\delta}$ in $\hat{\mathfrak{sl}}_l$
  is $l-1$. Note
  that $\dim (U^+(\hat{\mathfrak{sl}}_l))_\dvec{d}$ is equal to the number of separated isomorphism classes
  in $\Pi_x$
  of dimension vector $\dvec{d}$.

  Theorem 4.1 in \cite{DlabRingel_repsgraphsalgs} yields that the sum $\sum (l_i - 1)$
  of the ranks of the inhomogeneous tubes minus one is equal to $n-2$, the number of vertices of
  $Q$ minus two.

  Since the elements in the above form which are homogeneous of degree $\dvec{d}$
  generate $\mathcal{C}_0(Q)_\dvec{d}$, being of the right number, they have to be linearly independent.
  Hence, they are a basis and every element $u_w$ can be written uniquely in such
  a form. Since every element in $\EMon(Q)$ can be written
  in this partial normal form, we have that the $\dvec{d}$-th graded part of $\Q\EMon(Q)$, 
  $\Q\EMon(Q)_\dvec{d}$, has dimension at most $\dim \mathcal{C}_0(Q)_\dvec{d}$. Since
  $\Xi$ is surjective, this yields that it is an isomorphism and that the partial normal
  form is a normal form.
\end{proof}
\begin{cory}
  The map
  \[\Phi \colon\quad \mathcal{C}_0(Q) \rightarrow \Q \mathcal{CM}(Q),
  \]
  sending
  $u_{S_i}$ to $\Orbit_{S_i}$ for each $i \in Q_0$, is a surjective $\Q$-algebra homomorphism with kernel
  generated by $u_{r \dvec{\delta}} = (u_{\dvec{\delta}})^r$.
  \label{extdyn:cory:morph}
\end{cory}
\begin{proof}
  Since $\mathcal{C}_0(Q) \cong \Q\EMon(Q)$, $\mathcal{CM}(Q) \cong \TEMon(Q)$ and
  $\TEMon(Q)$ arises from $\EMon(Q)$ by dividing out the relation $\SRoot{r\dvec{\delta}} = \SRoot{\dvec{\delta}}^r$,
  the claim follows.
\end{proof}
\begin{cory}
  Let $Q$ be a connected, acyclic, extended Dynkin quiver
  with $r$ inhomogeneous tubes indexed by $x_1, \dots, x_r \in \Proj^1_\Z$.
  For every separated partition $\pi_i \in \Pi_{x_i}^s$ let $u_{(\pi_i, x_i)} \in \mathcal{C}_q(Q)$ be
  a lift of $u_{\overline{\Orbit_{M_i}}}\in \mathcal{C}_0(Q)$ for a representation $M_i \in \mathcal{T}_{x_i}$
  of isomorphism class $\pi_i$.
  Then the elements $\mathcal{X}\tensor_{\Q[q]} \Q(q) \in \mathcal{C}_q(Q) \tensor_{\Q[q]} \Q(q)$ of the form
    $\mathcal{X} = \mathcal{P} \mathcal{C}_1\dotsm \mathcal{C}_r\mathcal{R}\mathcal{I} \in \mathcal{C}_q(Q)$, where
    \begin{align*}
        \mathcal{P}&=u_{s_1\dvec{p}^1} \diamond \dots \diamond u_{s_k\dvec{p}^k} && \text{with } s_i > 0,\ \dvec{p}^i
        \text{ preprojective Schur roots and }\\
        &&&\dvec{p}^i \prec_t \dvec{p}^j \text{ for all } i < j;\\
        \mathcal{C}_i&=u_{(\pi_i,x_i)} && \text{for a separated } \pi_i \in \Pi_{x_i}^s;\\
        \mathcal{R}&=u_{\lambda_1\dvec{\delta}} \diamond \dots \diamond u_{\lambda_l\dvec{\delta}} && \text{with }
        \lambda=(\lambda_1\ge \dots \ge \lambda_l)
        \text{ a partition};\\
        \mathcal{I}&=u_{t_1\dvec{q}^1} \diamond \dots \diamond u_{t_m\dvec{q}^m} && \text{with } t_i > 0,\ \dvec{q}^i
        \text{ preinjective Schur roots and }\\
        &&&\dvec{q}^i \prec_t \dvec{q}^j \text{ for all } i < j.
    \end{align*}
    are a basis of $\mathcal{C}_q(Q) \tensor_{\Q[q]} \Q(q)$.
\end{cory}
\begin{proof}
    Since the images of the set of elements of the form above under the specialisation to $q=0$ are
    a basis, they are linearly independent. Since each graded part $\mathcal{C}_q(Q)\tensor_{\Q[q]} \Q(q)_\dvec{d}$ is
    a $\Q(q)$-vector space of dimension $\dim \mathcal{C}_0(Q)_\dvec{d}$, we have that the elements
    are a basis.
\end{proof}

\section{Example}
\begin{sidewaysfigure}
\caption{Degeneration order on the quiver $\widetilde {A}_2$}
\[
\makeatletter%
\let\ASYencoding\f@encoding%
\let\ASYfamily\f@family%
\let\ASYseries\f@series%
\let\ASYshape\f@shape%
\makeatother%
\setlength{\unitlength}{1pt}
\includegraphics{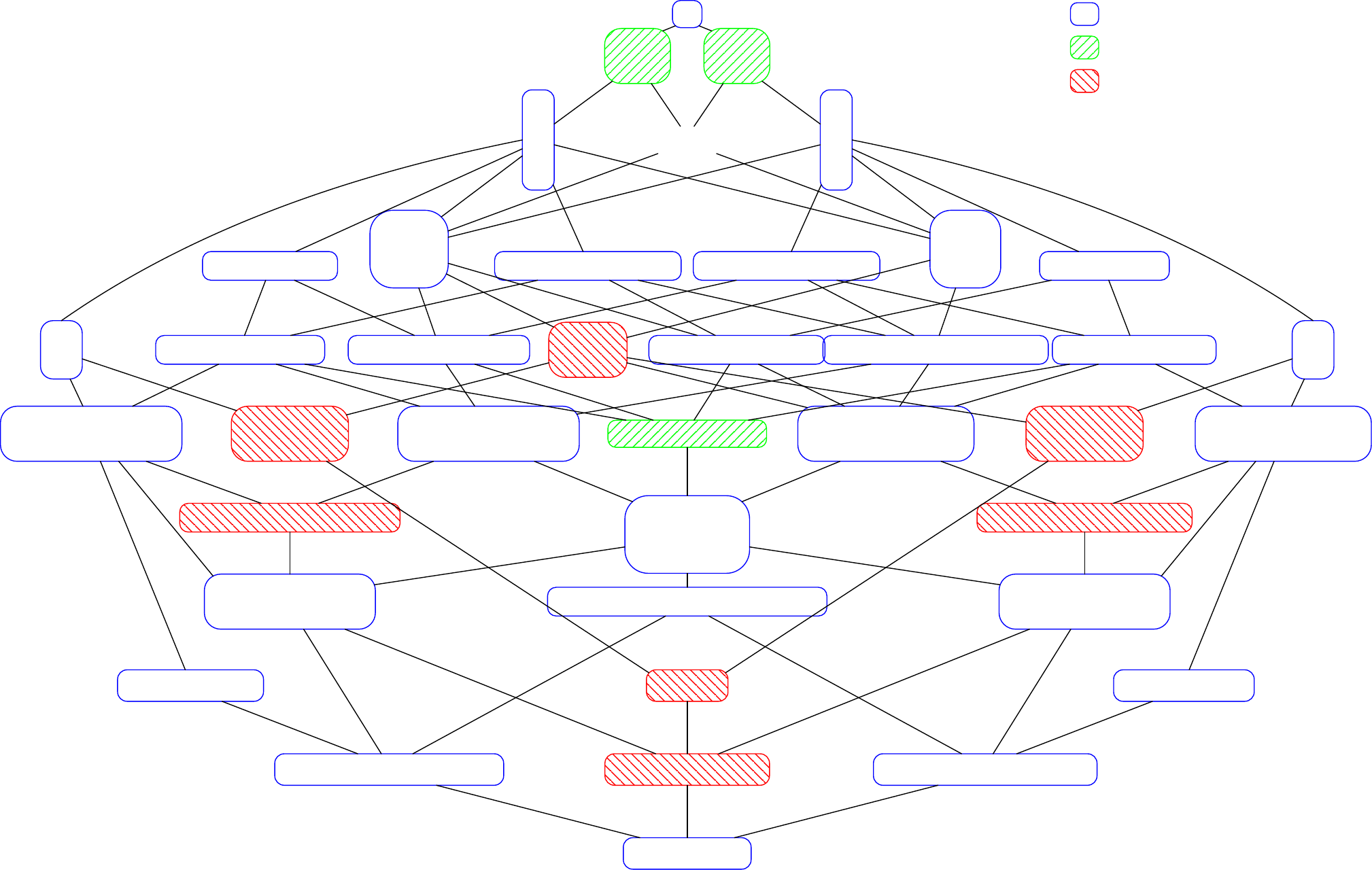}%
\definecolor{ASYcolor}{gray}{0.000000}\color{ASYcolor}
\fontsize{8.000000}{9.600000}\selectfont
\usefont{\ASYencoding}{\ASYfamily}{\ASYseries}{\ASYshape}%
\ASYalign(-122.123400,363.899422)(-0.500000,-0.500000){M}
\ASYalign(-109.676900,363.899422)(0.000000,-0.250000){$\in \mathcal{CM}(Q) \cap \{$test symbols$\}$}
\ASYalign(-122.123400,349.673044)(-0.500000,-0.500000){M}
\ASYalign(-109.676900,349.673044)(0.000000,-0.250000){$\in \mathcal{CM}(Q) \backslash \{$test symbols$\}$}
\ASYalign(-122.123400,335.446666)(-0.500000,-0.500000){M}
\ASYalign(-109.676900,335.446666)(0.000000,-0.250000){$\in \{$test symbols$\} \backslash \mathcal{CM}(Q)$}
\ASYalign(-291.075811,363.899422)(-0.500000,-0.255287){$2\underline{\delta}$}
\ASYalign(-312.194862,346.049215)(-0.500000,-0.104167){$\underline{\delta} \oplus \begin{matrix}S_2\\T\end{matrix}$}
\ASYalign(-269.956759,346.049215)(-0.500000,-0.104167){$\underline{\delta} \oplus \begin{matrix}T\\S_2\end{matrix}$}
\ASYalign(-354.432964,310.348801)(-0.500000,-0.052083){$\begin{matrix}S_2\\T\\S_2\\T\end{matrix}$}
\ASYalign(-227.718657,310.348801)(-0.500000,-0.052083){$\begin{matrix}T\\S_2\\T\\S_2\end{matrix}$}
\ASYalign(-291.075811,310.348801)(-0.500000,-0.255287){$\underline{\delta} \oplus S_2 \oplus T$}
\ASYalign(-172.809123,263.938262)(-0.500000,-0.069445){$\begin{matrix}S_2\\T\\S_2\end{matrix} \oplus T$}
\ASYalign(-409.342498,263.938262)(-0.500000,-0.069445){$\begin{matrix}T\\S_2\\T\end{matrix} \oplus S_2$}
\ASYalign(-113.675780,256.798179)(-0.500000,-0.250000){$I(2,2,1) \oplus S_3$}
\ASYalign(-333.313913,256.798179)(-0.500000,-0.250000){$I(2,1,1) \oplus P(0,1,1)$}
\ASYalign(-468.475841,256.798179)(-0.500000,-0.250000){$P(1,2,2) \oplus S_1$}
\ASYalign(-248.837708,256.798179)(-0.500000,-0.250000){$P(1,1,2) \oplus I(1,1,0)$}
\ASYalign(-269.956759,221.097765)(-0.500000,-0.250000){$I(2,1,1) \oplus S_2 \oplus S_3$}
\ASYalign(-396.671067,221.097765)(-0.500000,-0.250000){$P(1,1,2) \oplus S_2 \oplus S_1$}
\ASYalign(-101.004349,221.097765)(-0.500000,-0.250000){$\underline{\delta} \oplus I(1,1,0) \oplus S_3$}
\ASYalign(-481.147272,221.097765)(-0.500000,-0.250000){$\underline{\delta} \oplus S_1 \oplus P(0,1,1)$}
\ASYalign(-333.313913,221.097765)(-0.500000,-0.104167){$\begin{matrix}S_2\\T\end{matrix} \oplus \begin{matrix}T\\S_2\end{matrix}$}
\ASYalign(-557.175857,221.097765)(-0.500000,-0.130969){${\begin{matrix}S_2\\T\end{matrix}}^2$}
\ASYalign(-24.975764,221.097765)(-0.500000,-0.130969){${\begin{matrix}T\\S_2\end{matrix}}^2$}
\ASYalign(-185.480554,221.097765)(-0.500000,-0.250000){$I(1,1,0) \oplus T \oplus P(0,1,1)$}
\ASYalign(-291.075811,185.397351)(-0.500000,-0.255287){$\underline{\delta} \oplus S_1 \oplus S_2 \oplus S_3$}
\ASYalign(-206.599605,185.397351)(-0.500000,-0.104167){$\begin{matrix}S_2\\T\end{matrix} \oplus I(1,1,0) \oplus S_3$}
\ASYalign(-460.028221,185.397351)(-0.500000,-0.104167){$\begin{matrix}S_2\\T\end{matrix} \oplus S_2 \oplus T$}
\ASYalign(-544.504426,185.397351)(-0.500000,-0.104167){$\begin{matrix}S_2\\T\end{matrix} \oplus S_1 \oplus P(0,1,1)$}
\ASYalign(-37.647195,185.397351)(-0.500000,-0.104167){$\begin{matrix}T\\S_2\end{matrix} \oplus I(1,1,0) \oplus S_3$}
\ASYalign(-122.123400,185.397351)(-0.500000,-0.104167){$\begin{matrix}T\\S_2\end{matrix} \oplus S_2 \oplus T$}
\ASYalign(-375.552016,185.397351)(-0.500000,-0.104167){$\begin{matrix}T\\S_2\end{matrix} \oplus S_1 \oplus P(0,1,1)$}
\ASYalign(-291.075811,142.556854)(-0.500000,-0.069445){$\begin{matrix}S_2\\T\\S_2\end{matrix} \oplus S_1 \oplus S_3$}
\ASYalign(-122.123400,149.696937)(-0.500000,-0.250000){$I(1,1,0) \oplus S_2 \oplus T \oplus S_3$}
\ASYalign(-460.028221,149.696937)(-0.500000,-0.250000){$S_2 \oplus T \oplus S_1 \oplus P(0,1,1)$}
\ASYalign(-460.028221,113.996523)(-0.500000,-0.104167){$\begin{matrix}S_2\\T\end{matrix} \oplus S_1 \oplus S_2 \oplus S_3$}
\ASYalign(-122.123400,113.996523)(-0.500000,-0.104167){$\begin{matrix}T\\S_2\end{matrix} \oplus S_1 \oplus S_2 \oplus S_3$}
\ASYalign(-291.075811,113.996523)(-0.500000,-0.250000){$I(1,1,0) \oplus S_1 \oplus P(0,1,1) \oplus S_3$}
\ASYalign(-79.885298,78.296108)(-0.500000,-0.229289){$I(1,1,0)^2 \oplus S_3^2$}
\ASYalign(-291.075811,78.296108)(-0.500000,-0.229289){$S_2^2 \oplus T^2$}
\ASYalign(-502.266323,78.296108)(-0.500000,-0.229289){$S_1^2 \oplus P(0,1,1)^2$}
\ASYalign(-164.361503,42.595694)(-0.500000,-0.229289){$I(1,1,0) \oplus S_1 \oplus S_2 \oplus S_3^2$}
\ASYalign(-417.790118,42.595694)(-0.500000,-0.229289){$P(0,1,1) \oplus S_1^2 \oplus S_2 \oplus S_3$}
\ASYalign(-291.075811,42.595694)(-0.500000,-0.229289){$T \oplus S_1 \oplus S_2^2 \oplus S_3$}
\ASYalign(-291.075811,6.895280)(-0.500000,-0.229289){$S_1^2 \oplus S_2^2 \oplus S_3^2$}
\]
\label{fig:extdyn:degens}
\end{sidewaysfigure}
If $Q$ is a Dynkin quiver, we have shown that
\[
  f_w^\alpha (0) =
\begin{cases}
  1 & f_w^\alpha \neq 0,\\
  0 & \text{otherwise},
\end{cases}
\]
for an arbitrary decomposition symbol $\alpha$ and all words $w$ in vertices of $Q$.
We already saw that this is not true anymore in the extended Dynkin case.
\begin{defn}
	A decomposition symbol $\alpha$ is called a \emph{test symbol} if, for all
    words $w$ in vertices of $Q$,
\[
  f_w^\alpha (0) =
\begin{cases}
  1 & f_w^\alpha \neq 0,\\
  0 & \text{otherwise.}
\end{cases}
\]
\end{defn}

If it would we possible to distinguish between the sets $\mathcal{A}_w$ of the composition monoid
only by using test symbols, then it would be possible to prove that there is a homomorphism
\[ \Phi \colon \mathcal{C}_0(Q) \rightarrow \mathcal{CM}(Q) \]
by a similar method as we used in the Dynkin case. We will now illustrate that this is not possible
in general, namely if we have at least one inhomogeneous tube. The smallest example of this is the
quiver $\widetilde{A}_2$.
\[ \widetilde{A}_2 :=
\begin{array}{@{}c@{}}
\makeatletter%
\let\ASYencoding\f@encoding%
\let\ASYfamily\f@family%
\let\ASYseries\f@series%
\let\ASYshape\f@shape%
\makeatother%
\setlength{\unitlength}{1pt}
\includegraphics{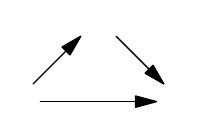}%
\definecolor{ASYcolor}{gray}{0.000000}\color{ASYcolor}
\fontsize{12.000000}{14.400000}\selectfont
\usefont{\ASYencoding}{\ASYfamily}{\ASYseries}{\ASYshape}%
\ASYalign(-52.462417,5.372285)(-0.500000,-0.500000){$1$}
\ASYalign(-28.452756,29.381946)(-0.500000,-0.500000){$2$}
\ASYalign(-4.443095,5.372285)(-0.500000,-0.500000){$3$}
\end{array}
\]
It has one inhomogeneous tube with regular simples $S_2$ and
$T$, where $\dimve T = (1,0,1)$.

We wrote down the degeneration graph of representations of
dimension vector $2\dvec{\delta}$ in figure \ref{fig:extdyn:degens}.
In there, we write classes of representations on the vertices. A connecting line
means, that the lower class is a subset of the closure of the $\GL_{2\dvec{\delta}}$-saturation
of the upper class. We write $P(d_1, d_2, d_3)$ for the unique indecomposable
preprojective representation of dimension vector $(d_1, d_2, d_3)$ and similarly
$I(d_1, d_2, d_3)$.

One easily sees that only
decomposition symbols $\alpha = (\mu, \sigma)$ having their regular part in one single tube
are test symbols. Moreover, just using test symbols, we cannot distinguish between
the following elements
\begin{align*}
    &\begin{matrix}T\\S_2\\T\\S_2\end{matrix}
        \text{ and }
        \Rep{\dvec{\delta}} \oplus \begin{matrix}T\\S_2\end{matrix},
            & \text{or between}&
&\begin{matrix}S_2\\T\\S_2\\T\end{matrix}
  \text{ and }
  \Rep{\dvec{\delta}} \oplus \begin{matrix}S_2\\T\end{matrix}
\end{align*}
since every test symbol which occurs in one class also occurs in the other class.
Similarly for
\[ \Rep{\dvec{\delta}} \oplus S_1 \oplus S_2 \oplus S_3
  \text{ and }
  \begin{matrix}S_2\\T\\S_2\end{matrix} \oplus S_1 \oplus S_3.\]

This immediately yields that this method cannot work in the extended
Dynkin case, therefore our hard work was necessary.

In general, it would be interesting to obtain the
morphism $\mathcal{C}_0(Q) \rightarrow \Q\mathcal{CM}(Q)$
in a global way, i.e. without using generators and relations. This could then possibly
generalise to the wild case, where it is hopeless to apply our method.

\appendix
\chapter{Fitting Ideals}
\label{app:fitting}
We define \emph{Fitting ideals}\index{Fitting ideals}
for morphisms and finitely presented modules and give some results for them.
For proofs see \cite{Northcott_finfreeres}.
Let $R$ be any commutative ring. A projective $R$-module $P$ is called of constant
rank $r$ if for every prime ideal $\mathfrak{p} \in \Spec R$ the localisation
$P_\mathfrak{p}$ is a free $R_\mathfrak{p}$-module of rank $r$.

\begin{defn}
	Let $f \colon R^d \rightarrow R^e$ be a homomorphism, which is given by a matrix
	$A \in R^{e \times d}$. For each $r \in \N$ define $\mathcal{F}_r (f)$ to be the ideal of
    $R$ generated by the $e-r$ minors of $A$. For $r \ge e$ set $\mathcal{F}_r(f) := R$.
\end{defn}
\begin{defn}
	Let $M$ be a finitely presented $R$-module and
	\[R^d \overset{f}{\rightarrow} R^e \rightarrow M \rightarrow 0 \]
	a presentation. Define $\mathcal{F}_r (M) := \mathcal{F}_r (f)$.
\end{defn}
\begin{remark}
	This definition does not depend on the choice of the presentation $f$.
\end{remark}
Here are the first basic properties of the Fitting ideals $\mathcal{F}_r$.
\begin{propn}
	Let $M$ be a finitely presented module. Then
	\begin{enumerate}
		\item $\mathcal{F}_r(M) \subset \mathcal{F}_{r+1}(M)$ for all $r \in \N$.
		\item If $M$ is a free module of rank $r$, then
			\[ (0) = \mathcal{F}_0 (M) = \mathcal{F}_1 (M) = \dots = \mathcal{F}_{r-1}(M)
			\subset \mathcal{F}_{r}(M) = R.\]
		\item If $S$ is an $R$-algebra, then $\mathcal{F}_r(M\tensor S) = \mathcal{F}_r(M) \tensor S$
			for all $r$.
	\end{enumerate}
\end{propn}

More generally, one has the following theorem (\cite{Campillo_projfittideal}).
\begin{theom}
	Let $M$ be a finitely presented $R$-module. Then $M$ is projective of constant rank $r$ if and only if
	$\mathcal{F}_{r-1}(M) = (0)$ and $\mathcal{F}_{r}(M) = R$.
\end{theom}

One can generalise this to finitely generated modules, but we do not need this.

\chapter{Tensor Algebras}
\label{app:tensalg}
Let $\Lambda_0$ be a ring and $\Lambda_1$ a $\Lambda_0$-bimodule. Define the tensor ring
$T(\Lambda_0, \Lambda_1)$ to be the $\N$-graded $\Lambda_0$-module
\[ \Lambda := \bigoplus_{r\ge 0} \Lambda_r, \quad \Lambda_r := \Lambda_1 \tensor_{\Lambda_0} \dotsm \tensor_{\Lambda_0} \Lambda_1 \: (r \text{ times}),\]
with multiplication given via the natural isomorphism $\Lambda_r \tensor_{\Lambda_0} \Lambda_s \cong \Lambda_{r+s}$. If $\lambda \in \Lambda_r$
is homogeneous, we write $\degr{ \lambda} = r$ for its degree.

The graded radical of $\Lambda$ is the ideal $\Lambda_+ := \bigoplus_{r \ge 1} \Lambda_r$. Note that
$\Lambda_+ \cong \Lambda_1 \tensor_{\Lambda_0} \Lambda$ as right $\Lambda$-modules.

\begin{lmm}
	Let $R$ and $S$ be rings. Let $A_R$, $_R B_S$ and $C_S$ be modules over the corresponding rings. Assume that $_R B_S$ is $S$-projective
	and $R$-flat. Then
	\[ \Ext^n_S (A \tensor_R B, C) \cong \Ext^n_R (A, \Hom_S(B,C))
	\]
\end{lmm}
\begin{proof}
	Choose a projective resolution $P_\bullet$ of $A_R$. The functor
    $- \tensor_R B$ is exact, and for any projective module $P_R$ the functor
	$\Hom_S(P \tensor_R B, -) \cong \Hom_R(P, -) \circ \Hom_S(B_S, -)$ is exact since $B_S$ is projective. Therefore,
	$P_\bullet \tensor_R B$ gives a projective resolution of $A \tensor_R B_S$ as an $S$-module. We obtain
	\begin{multline*}
		\Ext^n_S (A \tensor_R B, C) = H^n \Hom_S(P_\bullet \tensor_R B, C) \cong \\
		\cong H^n \Hom_R(P_\bullet, \Hom_S(B,C)) = \Ext^n_R (A, \Hom_S(B,C)).
	\end{multline*}
\end{proof}
\begin{theom}
	Let $\Lambda$ be a tensor ring and $M \in \MMod \Lambda$.
Then there is a short exact sequence
\[ 0 \rightarrow M \tensor_{\Lambda_0} \Lambda_+ \xrightarrow{\delta_M} M \tensor_{\Lambda_0} \Lambda \xrightarrow{\epsilon_M} M \rightarrow 0,\]
where, for $m \in M$, $\lambda \in \Lambda$ and
$\mu \in \Lambda_1$,
\begin{align*}
	\epsilon_M ( m \tensor \lambda)& := m \cdot \lambda\\
	\delta_M (m \tensor ( \mu \tensor \lambda) &:= m \tensor (\mu \tensor \lambda) - m \cdot \mu \tensor \lambda.
\end{align*}
Moreover, if $_{\Lambda_0} \Lambda_1$ is flat, then $\pd M \tensor_{\Lambda_0} \Lambda_+ \le \gldim \Lambda_0$ and
$\pd M \tensor_{\Lambda_0} \Lambda \le \gldim \Lambda_0$.
\end{theom}
\begin{proof}
	It is clear that $\epsilon_M$ is an epimorphism and that $\epsilon_M \delta_M = 0$. To see that $\delta_M$ is a monomorphism, we decompose
	$M \tensor \Lambda = \bigoplus_{r\ge 0} M \tensor \Lambda_r$ and similarly $M \tensor \Lambda_+$ as $\Lambda_0$-modules.
	Then $\delta_M$ restricts to
	maps $M \tensor \Lambda_r \rightarrow (M \tensor \Lambda_r) \oplus (M \tensor \Lambda_{r-1})$ for each $r \ge 1$ and moreover acts as
	the identity on the first component. In particular, if $\sum_{r=1}^t x_r \in \Ker(\delta_M)$ with $x_r \in M \tensor \Lambda_r$, then $x_t =0$.
	Thus $\delta_M$ is injective.

	Next we show that $M \tensor \Lambda = (M \tensor \Lambda_0) \oplus \Bild(\delta_M)$. Let $x = \sum_{r=0}^{t} x_r \in M \tensor \Lambda$.
	We show that $x \in (M \tensor \Lambda_0) + \Bild(\delta_M)$ by induction on $t$. For $t=0$ this is trivial. Let $t \ge 1$. Then
	$x_t \in M \tensor \Lambda_+$ and $x - \delta_M (x_t) = \sum_{r=0}^{t-1} x_r '$. So we are done by induction. If
	$x \in \Bild(\delta_M) \cap M \tensor \Lambda_0$, then there is a $y \in M \tensor \Lambda_+$ such that $\delta_M (y) =x \in M \tensor \Lambda_0$.
	But this is only the case if $y =0$ and therefore $x=0$.

	Now let $_{\Lambda_0} \Lambda_1$ be flat. Then $_{\Lambda_0}\Lambda$ is flat. If $N$ is any $\Lambda_0$-module,
	then $\pd_\Lambda N \tensor \Lambda \le \pd_{\Lambda_0} N$ since
	\[\Ext^n_\Lambda (N \tensor \Lambda, L) \cong \Ext^n_{\Lambda_0} (N, \Hom_{\Lambda}(\Lambda, L)) = \Ext^n_{\Lambda_0} (N, L)\] for any
	$L \in \MMod \Lambda$ by the lemma. Therefore, the second part follows since $\Lambda_+ = \Lambda_1 \tensor \Lambda$.
\end{proof}
Now let $Q$ be a quiver and $\Lambda_0 = R^{Q_0}$ for a fixed ring $R$. Let $\Lambda_1 := R^{Q_1}$ be the free $\Lambda_0$-bimodule
given by the arrows. The tensor ring $T(\Lambda_0, \Lambda_1)$ is then equal to $RQ$. For an $RQ$-module $M$ denote
$M\epsilon_i$ by
$M_i$. As $R$-modules we have that $M = \bigoplus M_i$.
\begin{theom}
	\label{euf_quiverquiver}
	Let $R$ be a ring with $\gldim R = n < \infty$. Let $M$ and $N$ be two modules of finite length over $RQ$. Then the Euler form is given by
	\[ \euf{M,N}_{RQ} = \sum_{i\in Q_0} \euf{M_i, N_i}_R - \sum_{\alpha\colon i \rightarrow j} \euf{M_i, N_j}_R. \]
\end{theom}
\begin{proof}
	Let $\Lambda_0 = R^{Q_0}$. There is a natural $\Lambda_0$-bimodule structure on $\Lambda_1 = R^{Q_1}$.
    Let $\Lambda = T(\Lambda_0, \Lambda_1) = RQ$.
	Let $0 \rightarrow M \tensor_{\Lambda_0} \Lambda_1 \tensor_{\Lambda_0} \Lambda \rightarrow M \tensor_{\Lambda_0} \Lambda \rightarrow M \rightarrow 0$
	be the short exact sequence of the previous theorem. Apply $\Hom(-, N)$ to it and consider the long exact sequence
	\[
	\xymatrix{
		0 \ar[r] & \Hom_\Lambda(M, N) \ar[r] & \Hom_\Lambda(M \tensor_{\Lambda_0} \Lambda , N) \ar[r]&
		    \Hom_\Lambda(M \tensor_{\Lambda_0} \Lambda_1 \tensor_{\Lambda_0} \Lambda, N) \ar `r[d] `[lll] `[dlll] `[dll] [dll]\\
		 & \Ext^1_\Lambda(M, N) \ar[r] &\Ext^1_\Lambda(M \tensor_{\Lambda_0} \Lambda , N) \ar[r]&
		    \Ext^1_\Lambda(M \tensor_{\Lambda_0} \Lambda_1 \tensor_{\Lambda_0} \Lambda, N) \\
			&\cdots & &  \\
		 **[r]\quad\ar[r]& \Ext^n_\Lambda(M, N) \ar[r] &\Ext^n_\Lambda(M \tensor_{\Lambda_0} \Lambda , N) \ar[r]&
		    \Ext^n_\Lambda(M \tensor_{\Lambda_0} \Lambda_1 \tensor_{\Lambda_0} \Lambda, N) \\
		**[r]\quad\ar[r]& \Ext^{n+1}_\Lambda(M, N) \ar[r] & **[l]0.\qquad\qquad &}
	\]
    Here we obtain the last $0$ since the $\pd M \tensor \Lambda \le \gldim R$
    by the previous theorem.

    Now, since $ _{\Lambda_0} \Lambda_\Lambda$ is projective as a $\Lambda_0$-module and as a $\Lambda$-module, we obtain
	\[\Ext^i_\Lambda(M \tensor_{\Lambda_0} \Lambda , N) = \Ext^i_{\Lambda_0}(M,N) = \bigoplus_{i \in Q_0} \Ext^i_R(M_i, N_i)\]
	and
	\[\Ext^i_{\Lambda}(M \tensor_{\Lambda_0} \Lambda_1 \tensor_{\Lambda_0} \Lambda, N) = \Ext^i_{\Lambda_0} (M \tensor_{\Lambda_0} \Lambda_1, N)=
	\bigoplus _{\alpha \colon i \rightarrow j} \Ext^i_R(M_i, N_j). \]
	This yields the claim.
\end{proof}
\begin{remark}
    Note that, for a field $K$, we recover C. M. Ringel's result, namely that
    \[
    \sum_{i\in Q_0} \dim M_i \dim N_i - \sum_{\alpha\colon i \rightarrow j} \dim M_i \dim N_j = [M,N]_{KQ} - [M,N]^1_{KQ}
    \]
    for two $K$-representations $M$ and $N$ of a quiver $Q$.
\end{remark}

\backmatter
\bibliography{repsofalgs}
\printindex
\end{document}